\documentclass[11pt]{amsart}
\usepackage{amsmath,amssymb,amsthm,mathtools}
\usepackage{mathrsfs}
\usepackage{enumitem}
\usepackage{caption}
\usepackage{graphicx}
\usepackage{tikz}
\usetikzlibrary{arrows.meta,calc,decorations.markings,decorations.pathreplacing,fit,patterns,positioning}
\usepackage[colorlinks=true,hyperindex,pagebackref]{hyperref}
\newcommand\myshade{85}
\colorlet{mylinkcolor}{violet}
\colorlet{mycitecolor}{red}
\colorlet{myurlcolor}{cyan}
\hypersetup{
  linkcolor  = mylinkcolor!\myshade!black,
  citecolor  = mycitecolor!\myshade!black,
  urlcolor   = myurlcolor!\myshade!black,
  colorlinks = true,
}
\renewcommand*{\backref}[1]{}
\renewcommand*{\backrefalt}[4]{%
  \ifcase#1\relax
  \or\space(Cited on p.~#2.)%
  \else\space(Cited on pp.~#2.)%
  \fi
}
\usepackage{aliascnt}
\usepackage[nameinlink,capitalise,noabbrev]{cleveref}

\hypersetup{
  pdftitle={Global-dimension spectra of gentle algebras},
  pdfauthor={Wen Chang and Quanyu Tang}
}

\theoremstyle{definition}
\newtheorem{problem}{Problem}
\theoremstyle{plain}
\newtheorem{theorem}{Theorem}[section]

\newaliascnt{proposition}{theorem}
\newtheorem{proposition}[proposition]{Proposition}
\aliascntresetthe{proposition}

\newaliascnt{lemma}{theorem}
\newtheorem{lemma}[lemma]{Lemma}
\aliascntresetthe{lemma}

\newaliascnt{corollary}{theorem}
\newtheorem{corollary}[corollary]{Corollary}
\aliascntresetthe{corollary}

\theoremstyle{definition}
\newaliascnt{definition}{theorem}
\newtheorem{definition}[definition]{Definition}
\aliascntresetthe{definition}

\newaliascnt{example}{theorem}
\newtheorem{example}[example]{Example}
\aliascntresetthe{example}

\theoremstyle{definition}
\newaliascnt{remark}{theorem}
\newtheorem{remark}[remark]{Remark}
\aliascntresetthe{remark}

\crefname{problem}{Problem}{Problems}
\Crefname{problem}{Problem}{Problems}
\crefname{theorem}{Theorem}{Theorems}
\Crefname{theorem}{Theorem}{Theorems}
\crefname{proposition}{Proposition}{Propositions}
\Crefname{proposition}{Proposition}{Propositions}
\crefname{lemma}{Lemma}{Lemmas}
\Crefname{lemma}{Lemma}{Lemmas}
\crefname{corollary}{Corollary}{Corollaries}
\Crefname{corollary}{Corollary}{Corollaries}
\crefname{definition}{Definition}{Definitions}
\Crefname{definition}{Definition}{Definitions}
\crefname{example}{Example}{Examples}
\Crefname{example}{Example}{Examples}
\crefname{remark}{Remark}{Remarks}
\Crefname{remark}{Remark}{Remarks}
\crefname{figure}{Figure}{Figures}
\Crefname{figure}{Figure}{Figures}

\newcommand{\Db}{\mathrm D^{\mathrm b}}

\newcommand{\Kb}{\mathrm K^{\mathrm b}}
\newcommand{\per}{\operatorname{per}}
\newcommand{\modu}{\operatorname{mod}}
\newcommand{\proj}{\operatorname{proj}}
\newcommand{\thick}{\operatorname{thick}}
\newcommand{\gldim}{\operatorname{gl.dim}}
\newcommand{\degmax}{\deg_{\max}}

\newcommand{\REnd}{\mathbf R\!\operatorname{End}}
\newcommand{\Hom}{\operatorname{Hom}}
\newcommand{\AG}{\operatorname{AG}}
\newcommand{\GSpec}{\mathcal G}
\newcommand{\AGSpec}{\mathcal G_{\mathrm{AAG}}}
\newcommand{\ZZ}{\mathbb Z}
\newcommand{\NN}{\mathbb N}

\newcommand{\cH}{\mathcal H}
\newcommand{\cT}{\mathcal T}

\newcommand{\cQ}{\mathcal Q}
\newcommand{\Arf}{\operatorname{Arf}}
\newcommand{\Ext}{\operatorname{Ext}}

\tikzset{
  arc/.style={line width=.85pt,line cap=round,line join=round},
  boundary/.style={line width=1.15pt,line cap=round,line join=round},
  construction/.style={draw=black!48,line width=.55pt},
  auxiliary/.style={draw=black!55,densely dashed,line width=.55pt},
  seam/.style={draw=black!60,dash pattern=on 3pt off 2pt,line width=.6pt},
  diagonal/.style={line width=1.2pt,line cap=round,line join=round},
  whitevertex/.style={circle,draw=black,line width=.65pt,fill=white,
                     inner sep=0pt,minimum size=4.7pt},
  blackvertex/.style={circle,fill=black,inner sep=0pt,minimum size=4.1pt},
  spinevertex/.style={blackvertex},
  panel/.style={font=\small,align=center},
  figlabel/.style={font=\small,fill=white,inner sep=1.2pt},
  flow/.style={line width=.7pt,-{Latex[length=2.1mm,width=1.45mm]}},
  identification/.style={draw=black!65,line width=.65pt,
                         -{Latex[length=1.8mm,width=1.25mm]}},
  rotation/.style={draw=black!60,line width=.55pt,-{Latex[length=1.65mm]}},
  pairedge/.style={line width=.75pt,postaction={decorate},
    decoration={markings,mark=at position .54 with {\arrow{Latex[length=1.8mm,width=1.2mm]}}}},
  leader/.style={draw=black!65,line width=.5pt},
  heartpiece/.style={draw=black!50,fill=black!4,rounded corners=2pt,
                    line width=.55pt,minimum width=1.2cm,minimum height=.55cm}
}

\title[Global-dimension spectra of gentle algebras]
{Global-dimension spectra of gentle algebras}
\author{Wen Chang}
\address{School of Mathematics and Statistics, Shaanxi Normal University,
Xi'an 710062, P.~R.~China}
\email{changwen161@163.com}
\author{Quanyu Tang}
\address{School of Mathematical Sciences, University of Science and
Technology of China, Hefei 230026, P.~R.~China}
\email{tangquanyu827@gmail.com}
\date{}
\begin{document}

\begin{abstract}
We classify the connected finite-dimensional gentle algebras
of finite global dimension whose global-dimension spectra under
derived equivalence are not intervals of integers.
These are precisely the algebras with AAG invariant $[(2,4g)]$,
where $g\geq2$, whose canonical line field is induced by a
nowhere-vanishing vector field.
Their spectrum is $\{2g,2g+2,\ldots,4g\}$.
For each $g\geq2$, there are exactly two such derived-equivalence
classes, distinguished by their Arf invariant.

The geometric foundation is a realisation theorem for prescribed
polygon weights in arbitrary genus.
For an AAG invariant $[(m_i,n_i)]_{i=1}^b$ arising from such an
algebra, the global dimensions attained in its AAG fibre form
the integer interval
$[\max_i\lceil n_i/m_i\rceil,\max_i n_i]$.
By controlling the complete line-field invariant through local
ribbon-graph constructions, we determine the spectrum in every
derived-equivalence class.
\end{abstract}

\maketitle

\tableofcontents
\section{Introduction}
\label{sec:introduction}

Tilting theory supplies one of the basic mechanisms for comparing
module categories as well as their derived categories
\cite{HappelRingel,Happel1987}.
For example, if a derived equivalence between two finite-dimensional algebras
$A$ and $B$ of finite global dimension is induced by a classical
tilting module of projective dimension at most one, then
$|\gldim A-\gldim B|\leq1$
\cite[Proposition~III.3.4]{Happel1988}. In general, Rickard's derived Morita theory
shows that an arbitrary derived equivalence between finite-dimensional
algebras is controlled by a tilting complex \cite{Rickard}. In this case, it is well known that finiteness of global dimension is preserved \cite[12.5]{GabrielRoiter}.  

More generally, for a tilting object in the derived category of an abelian category, a bound for the global dimension of the endomorphism algebra is given in \cite{KellerKrause}. However, note that such bounds do
not by themselves give uniform information about every algebra in one fixed
derived-equivalence class.
Thus, the following two questions are natural, which are raised by Kalck in 2016
and recorded by Schr\"oer in \cite[Questions~13.40 and~13.41]{SchroerAtlas}.

\begin{problem}
\label{prob:interval}
Let $A$ and $B$ be finite-dimensional algebras of finite global dimension such
that $\Db(\modu A)\simeq\Db(\modu B)$.  If $d$ is an integer satisfying
$\gldim A\leq d\leq\gldim B$, must there exist a finite-dimensional algebra
$C$ with $\Db(\modu C)\simeq\Db(\modu A)$ and $\gldim C=d$?
\end{problem}

\begin{problem}
\label{prob:boundedness}
Let $A$ be a finite-dimensional algebra of finite global dimension.  Does
there exist $b_A\geq0$ such that every finite-dimensional algebra $B$ with
$\Db(\modu B)\simeq\Db(\modu A)$ satisfies
$|\gldim A-\gldim B|\leq b_A$?
\end{problem}

For a finite-dimensional algebra $A$ of finite global dimension, define the \emph{global dimension spectrum} of $A$ under derived equivalence:
\begin{equation}
\label{eq:GSpec-definition}
 \GSpec(A):=
 \{\gldim B\mid B\text{ is finite-dimensional and }
                   \Db(\modu B)\simeq\Db(\modu A)\}.
\end{equation}
Thus \cref{prob:interval} asks whether $\GSpec(A)$ is always an interval of
integers, while \cref{prob:boundedness} asks whether it is always bounded.
The latter is known in important special settings, including piecewise
hereditary categories \cite{HappelReitenSmalo,Ladkani}.  The former follows
whenever the relevant representatives can be joined by a chain of tilting
steps induced by modules of projective dimension at most one, because each
such step changes global dimension by at most one.

This paper uses a special kind of algebras--gentle algebras--to test these problems. 
Gentle algebras
arose in the study of tilted and iterated tilted algebras of types $A$
and $\widetilde A$ \cite{AssemHappel,AssemSkowronski}, and the class of
gentle algebras is closed under derived equivalence
\cite{SchroerZimmermann}. At the module level, it is shown that the endomorphism algebra
of a module without self-extensions over a gentle algebra is again a 
gentle algebra \cite{Schroer1999Modules}. Avella-Alaminos and Geiss introduced the so-called AAG invariant for a gentle algebra \cite{AvellaAlaminosGeiss}, which is a
combinatorial derived invariant defined by using forbidden threads.  For a gentle algebra of finite global dimension, the global
dimension is the maximum length of a forbidden thread.

More recently, a geometric model for the derived category of any finite-dimensional gentle algebra in
  terms of graded curves on a marked surface is constructed in \cite{OpperPlamondonSchroll}.
  For graded gentle algebras, it is proved that the derived categories are equivalent to the partially wrapped Fukaya categories
  of graded marked surfaces \cite{HaidenKatzarkovKontsevich, LekiliPolishchuk}, and this correspondence is used to interpret derived
  equivalences through the associated line fields in \cite{LekiliPolishchuk}.
Building on \cite{LekiliPolishchuk}, a complete derived invariant for any finite-dimensional gentle algebra is constructed in 
\cite{AmiotPlamondonSchroll} by using boundary winding numbers
and, according to the genus, divisibility, parity, and Arf data. 

By deforming the surface model for the module category developed
in \cite{BaurCoelhoSimoes}, the first author realised each
indecomposable module as an object in the surface model of the
derived category and described its minimal projective resolution
and projective dimension \cite{Chang}. Taking the maximum over all indecomposable modules gives the
global dimension of the corresponding gentle algebra.  An equivalent formulation in which the global dimension is the
largest polygon weight of the dissection is given in
\cite{LiuGaoHuang}.  

Based on the above description, we aim to determine the set of
values obtained when the admissible dissection varies subject to a fixed
surface.
Let $A$ be a finite-dimensional gentle algebra of finite global dimension.  Display its
AAG invariant, including multiplicities, as
\begin{equation}
\label{eq:intro-AAG-list}
 \AG(A)=[(m_1,n_1),\ldots,(m_b,n_b)],
\end{equation}
where $m_i\geq 1$ since the global dimension of $A$ is finite. 

Put
\begin{equation}
\label{eq:LU-definition}
 L_A:=\max_i\left\lceil\frac{n_i}{m_i}\right\rceil,
 \qquad U_A:=\max_i n_i,
\end{equation}
and define the \emph{AAG-global dimension spectrum} obtained by fixing the AAG invariant by
\begin{equation}
\label{eq:AAG-fibre-definition}
\begin{aligned}
 \AGSpec(A):=\{\,\gldim C\mid {}&C\text{ is a connected gentle algebra of finite }\\[-2pt]& \text{ global dimension such that} \AG(C)=\AG(A)\,\}.
\end{aligned}
\end{equation}

Our first theorem shows that the global dimension spectrum of the AAG invariant is an interval.

\begin{theorem}
\label{thm:intro-AAG-realization}
Let $A$ be a finite-dimensional gentle algebra of finite global dimension.  Then
\begin{equation}
\label{eq:intro-AAG-realization}
 \AGSpec(A)=[L_A,U_A].
\end{equation}
In particular, every finite-dimensional algebra $B$ derived equivalent to $A$
satisfies
\begin{equation}
\label{eq:intro-enclosure}
 L_A\leq\gldim B\leq U_A.
\end{equation}
\end{theorem}

Therefore \cref{prob:boundedness} has an affirmative answer for every
gentle algebra of finite global dimension. 
We mention that for gentle algebras, boundedness in \cref{prob:boundedness} already follows
from \cite[Corollary~1.2]{SchroerZimmermann}, which shows that a fixed
derived-equivalence class contains only finitely many Morita-equivalence
classes. The enclosure \eqref{eq:intro-enclosure} also follows directly
from the definition of the AAG invariant and the forbidden-thread formula
for global dimension \cite{AvellaAlaminosGeiss}.
The realisation assertion in \cref{thm:intro-AAG-realization} shows that
every integer in the enclosure occurs when only the AAG data are fixed.
The distinction between $\AGSpec(A)$ and $\GSpec(A)$ is essential:
the former fixes only the boundary part of the surface invariant,
whereas the latter also fixes the complete line-field orbit.
When the AAG data already determine the line-field orbit, the finite formula
simplifies to a closed interval, see \cref{thm:intro-derived-interval}.
The complete invariant nevertheless gives an exact finite description of the
ordinary spectrum.

\begin{theorem}
\label{thm:intro-finite-formula}
For every finite-dimensional gentle algebra $A$ of finite global dimension, there is a
finite, explicitly defined set $\mathfrak R(A)$ of ribbon graphs with one
distinguished corner at each vertex, determined by the AAG data and the
complete numerical line-field invariant of $A$, such that
\begin{equation}
\label{eq:intro-finite-formula}
 \GSpec(A)=\{\degmax(G)-1\mid G\in\mathfrak R(A)\}.
\end{equation}
Here $\degmax(G)$ is the maximum vertex degree.  In particular, $\GSpec(A)$ is
computable by a finite enumeration.
\end{theorem}

The finite description can be sharpened to explicit formulas for every
derived-equivalence class. The following theorem collects the formulas
proved in \cref{sec:spectrum-classification}, together with the genus-zero
case of \cref{thm:intro-derived-interval}.
For a canonical line field $\eta_A$, write $d(\eta_A)$ for the genus-one
winding-divisibility invariant, see \eqref{eq:genus-one-gcd}, and let $\sigma(\eta_A)=0$ mean that every
closed curve has even winding number.
The following theorem gives the complete global-dimension spectrum for any gentle algebra.
\begin{theorem}
\label{thm:intro-complete-spectrum}
Let $A$ be a connected finite-dimensional gentle algebra of finite
global dimension, and let $g$ be the genus of its associated surface.
Write
\[
 \AG(A)=[(m_i,n_i)]_{i=1}^b,
 \qquad m=\sum_{i=1}^b m_i,
 \qquad w_i=m_i-n_i,
 \qquad L=L_A,\quad U=U_A.
\]
For $g=1$, put $D=d(\eta_A)$.
Then $\GSpec(A)$ is given by the following table.

\begin{center}
\begingroup
\normalfont\small
\setlength{\tabcolsep}{4pt}
\renewcommand{\arraystretch}{1.3}
\begin{tabular}{@{}p{0.17\linewidth}p{0.48\linewidth}p{0.27\linewidth}@{}}
\hline
\textbf{Case}
&
\textbf{Additional conditions}
&
\textbf{Spectrum $\GSpec(A)$}
\\
\hline
\textup{(i)} $g=0$
&
None
&
$[L,U]$
\\
\hline
\textup{(ii)} $g=1$
&
Some $w_i\neq-2$
&
$[L,U]$
\\
&
All $w_i=-2$, $D=0$, and $b=1$
&
$[2,m]$
\\
&
All $w_i=-2$, $D=0$, and $b\geq2$
&
$[L,U]$
\\
&
All $w_i=-2$ and $D>0$
&
$[L_D,U]$
\\
\hline
\textup{(iii)} $g\geq2$
&
$\sigma(\eta_A)=1$ or $b\geq2$
&
$[L,U]$
\\
&
$b=1$, $\sigma(\eta_A)=0$, and $m=2$
&
$\{2g,2g+2,\ldots,4g\}$
\\
&
$b=1$, $\sigma(\eta_A)=0$, and $m\geq3$
&
$[L+\varepsilon,U]$
\\
\hline
\end{tabular}
\endgroup
\end{center}

In case~\textup{(ii)}, the row with $D=0$ and $b=1$ requires
$m\geq2$. In the positive-$D$ row, one has $1\leq D\leq m$, and
\begin{equation}
\label{eq:intro-spectrum-LD}
 L_D:=\min\left\{r\in[L,U]\;\middle|\;
   \sum_{i=1}^b\left\lceil\frac{m_i+2}{r}\right\rceil
       \leq m-D+b\right\}.
\end{equation}
The set defining $L_D$ is nonempty. In particular, if all $w_i=-2$,
$b=1$, and $D>0$, then
\[
 \GSpec(A)=
 \left[
   \left\lceil\frac{m+2}{m-D+1}\right\rceil,
   m+2
 \right].
\]

In case~\textup{(iii)}, the conditions $b=1$ and
$\sigma(\eta_A)=0$ require $m\geq2$.
In the last row, $U=n_1=m+4g-2$ and
\[
 \varepsilon=
 \begin{cases}
  1,&n_1=mL\text{ and }L\text{ is odd},\\
  0,&\text{otherwise}.
 \end{cases}
\]
Whenever the Arf invariant is defined, its value does not affect
the spectrum.
\end{theorem}

Thus every spectrum is an integer interval except when
$g\geq2$, $b=1$, $m=2$, and $\sigma(\eta_A)=0$.
Equivalently, we obtain the following classification of all internal gaps.

\begin{theorem}
\label{thm:noninterval-classification}
Let $A$ be a finite-dimensional gentle algebra of finite global
dimension.  Then $\GSpec(A)$ is not an interval if and only if
\[
 \AG(A)=[(2,4g)]\quad\text{for some }g\geq2, 
 \quad \text{and}  \quad  \sigma(\eta_A)=0.
\]
In this case
\[
 \GSpec(A)=\{2g,2g+2,\ldots,4g\}.
\]
For each $g\geq2$ there are exactly two such derived-equivalence classes,
distinguished by their Arf invariant, and both have this spectrum.
\end{theorem}
The above theorem gives a family of counterexamples to
\cref{prob:interval}.  An explicit example with five simple modules
is given in \cref{thm:intro-counterexample}. The same classification also yields derived-equivalent gentle algebras
which are not tilting--cotilting equivalent, recovering the existence
phenomenon established by Ni~\cite[Theorem~1.2]{Ni2026} for each odd
number $2g+1\geq5$ of simple modules.
In fact, each exceptional derived-equivalence class of genus $g$
contains at least $g+1$ tilting--cotilting equivalence classes;
see \cref{rem:tilting-cotilting-existence}.

We conclude with consequences for tilting and silting
mutation. Combining the classification with the
global-dimension estimate of \cite[Lemma~4.6(4)]{ChenLinZhou}, we show that the
exceptional gentle classes have at least $g+1$ connected
components in the subgraph induced by tilting objects;
see \cref{cor:classified-tilting-components}.
The same estimate gives an integer interval of relative
global dimensions on each silting component. Thus
silting objects provide interpolation between tilting
objects lying in the same component, as explained in
\cref{rem:silting-dg-interpolation}.




The paper is organized as follows.
Section~\ref{sec:surface-data} recalls the necessary
background and fixes our surface conventions and notation.
Sections~\ref{sec:quadrangulations} and~\ref{sec:arc-realization}
establish the realisation theorem.
Section~\ref{sec:finite-formula} proves the finite formula
and illustrates it in a genus-one example.
Section~\ref{sec:spectrum-classification} determines the
ordinary spectra and classifies all noninterval cases.
Appendix~\ref{sec:HRS-components} records consequences
for tilting and silting mutation. In particular, the
exceptional gentle classes have at least $g+1$ connected
components in the subgraph induced by tilting objects.

Throughout, $k$ is an algebraically closed field, algebras
are nonzero unital $k$-algebras, and modules are right modules.
Ordinary finite-dimensional algebras are assumed to be basic.
Since Morita equivalence preserves global dimension and
bounded derived categories, this convention does not change
the spectra considered here. Derived equivalences are
$k$-linear triangle equivalences. Unless stated otherwise,
the gentle algebras under consideration are connected.
We write $\NN=\ZZ_{\geq0}$ and use intervals to mean intervals
of integers. Quiver paths are composed from left to right,
whereas maps are composed from right to left.

\section{Gentle algebras and surface models}
\label{sec:surface-data}
There are several references for the geometric models of various categories associated with gentle algebras.
For our purpose, we follow the notations as used in \cite{ChangJinSchroll}. The $\circ$-dissection $\Delta$ defines the algebra $A(\Delta)$.
Its dual $\bullet$-dissection $D:=\Delta^*$ is used to read the
homological data and to construct the canonical line field
associated with $A(\Delta)$, following \cite{AmiotPlamondonSchroll}.

\noindent\begin{minipage}{\linewidth}
\subsection{Gentle algebras and admissible \texorpdfstring{$\circ$}{o}-dissections}

A finite-dimensional algebra $A=kQ/I$ is \emph{gentle} if the following
conditions hold:
\begin{enumerate}[label=\textup{(\roman*)},leftmargin=2.6em]
\item each vertex of $Q$ is the source of at most two arrows and the target
      of at most two arrows;
\item the ideal $I$ is generated by paths of length two;
\item for every arrow $\alpha$, there is at most one arrow $\beta$ such
      that $\alpha\beta$ is composable and $\alpha\beta\notin I$, and at
      most one arrow $\gamma$ such that $\gamma\alpha$ is composable and
      $\gamma\alpha\notin I$;
\item for every arrow $\alpha$, there is at most one arrow $\beta$ such
      that $\alpha\beta$ is composable and $\alpha\beta\in I$, and at
      most one arrow $\gamma$ such that $\gamma\alpha$ is composable and
      $\gamma\alpha\in I$.
\end{enumerate}
\end{minipage}\par

A \emph{marked surface} is a pair $(S,M)$, where $S$ is a connected compact
oriented surface with nonempty boundary and $M$ is a finite set of marked
points on its boundary $\partial S$.  Write
\[
 \partial S=\coprod_{i=1}^{b}\partial_iS,
 \qquad M=M^{\circ}\amalg M^{\bullet}\subseteq\partial S.
\]
The elements of $M^{\circ}$ and $M^{\bullet}$ are called $\circ$-points and
$\bullet$-points, respectively.  These points alternate on each boundary
component, with at least one point of each colour.  In the general notation
of \cite{ChangJinSchroll}, our assumption is
$P^{\circ}=P^{\bullet}=\varnothing$, that is, the surface has no punctures.
We also exclude the zero surface, namely the disc with only one point of
each colour.

An \emph{$\circ$-arc} has endpoints in $M^{\circ}$ and is a non-contractible
curve; $\bullet$-arcs have endpoints in $M^{\bullet}$ instead.  Curves are
taken up to homotopy relative to the marked points, with transverse
interior intersections.

\begin{definition}
\label{def:admissible-dissection}
An \emph{admissible dissection} $\Delta$ of $(S,M)$ consists
of distinct $\circ$-arcs without interior self-intersections
or mutual interior intersections, which cut the surface into polygons, each containing precisely one
$\bullet$-point.
\end{definition}

We regard the polygons of $\Delta$ as complementary regions
in $S$. Their sides and corners are counted with multiplicity
along their boundary walks, so the same arc or marked point
may occur more than once.
Write $\mathcal P(\Delta)$ for these polygons, and let $P_p$
denote the unique polygon containing $p\in M^\bullet$.
For a polygon $P$, let $q(P)$ be the number of occurrences
of dissection sides in its boundary walk, excluding its
unique boundary side.

We use the following orientation conventions. Along the
boundary walk of a polygon, clockwise traversal keeps its
interior on the right, whereas anticlockwise traversal
keeps it on the left. Each boundary component of $S$ is
oriented so that $S$ lies on the right.
At a boundary endpoint $p$, use an orientation-preserving local
picture in which $p$ is the origin, the boundary is horizontal,
and the surface lies above it.  Draw the initial segments of
the incident arcs as distinct rays pointing away from $p$ into
the upper half-plane.  Their clockwise order is the order in
which they meet a small upper semicircle traversed from left
to right; their anticlockwise order is obtained by traversing
that semicircle from right to left.  

Let $q$ be a common endpoint of arcs $\ell_i,\ell_j$ in $\Delta$.
An \emph{oriented intersection} from $\ell_i$ to $\ell_j$ at $q$ is a local
anticlockwise angle from $\ell_i$ to $\ell_j$, based at $q$ and contained
in the interior of the surface.
It is \emph{minimal} when no arc
of $\Delta$ divides that angle.  The algebra
\[
 A(\Delta)=kQ(\Delta)/I(\Delta)
\]
has the arcs of $\Delta$ as vertices and their minimal oriented
intersections as arrows.  For composable arrows
$\alpha:\ell_i\to\ell_j$ and $\beta:\ell_j\to\ell_k$, the product
$\alpha\beta$ is a generating relation precisely when the two angles use
opposite ends of $\ell_j$.  If the endpoints of $\ell_j$ coincide, its two
end occurrences are still distinguished.

A grading $\mathfrak g$ assigns an integer $\mathfrak g(\alpha)$ to each
minimal oriented intersection $\alpha$.  It gives a \emph{graded admissible
dissection} $(\Delta,\mathfrak g)$ and a \emph{graded marked surface}
$(S,M,\Delta,\mathfrak g)$.  Giving each arrow $\alpha$ of $Q(\Delta)$ the
degree $|\alpha|=\mathfrak g(\alpha)$ defines the graded gentle algebra
$A(\Delta,\mathfrak g)$, whose underlying ordinary algebra is $A(\Delta)$.

Unless a grading is explicitly specified, $A(\Delta)$ denotes
the ordinary algebra $A(\Delta,0)$, with all arrow degrees zero.

Write $D:=\Delta^*$ for the dual $\bullet$-dissection: its arc
$\ell^*$ crosses the corresponding $\ell\in\Delta$ once and
crosses no other arc of $\Delta$.  The dual dissection is unique
up to homotopy.

Conversely, every connected finite-dimensional gentle algebra of
finite global dimension is isomorphic to $A(\Delta)$ for an admissible
$\circ$-dissection on a marked surface of the above kind
\cite{OpperPlamondonSchroll}.
Indeed finite
dimensionality excludes $\circ$-punctures, while finite global dimension
excludes $\bullet$-punctures. In fact, a gentle algebra has finite global dimension if and only if
it has no oriented cycles with full relations; this also follows
from the description of its singularity category
\cite[Theorem~2.5]{Kalck2015}.

For $p\in M^\bullet$, let $\nu_D(p)$ be the number of endpoint
occurrences of $D$ at $p$.  An arc in $D$ whose endpoints both equal $p$
contributes two. Since $p$ lies on the boundary, these occurrences form a
linearly ordered fan, with $\nu_D(p)-1$ successive minimal angles.
When identifying this fan with a forbidden thread of $A(\Delta)$,
we read its endpoint occurrences clockwise.  Indeed, listing the
dissection-side occurrences $\ell_1,\ldots,\ell_d$ of $P_p$
clockwise after its boundary side gives arrows
$\ell_k\to\ell_{k+1}$ at the successive $\circ$-corners:
the corresponding angles at these corners are anticlockwise.
The dual endpoint occurrences at $p$ have the same
clockwise order $\ell_1^*,\ldots,\ell_d^*$.

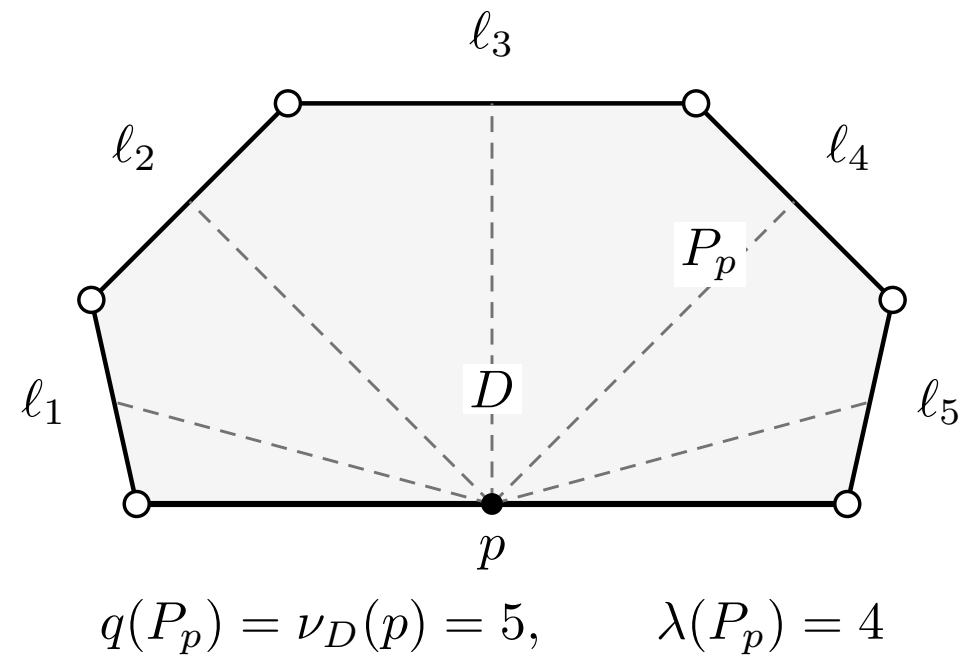
\begin{figure}[htbp]
\centering
\begin{tikzpicture}[x=1cm,y=1cm,font=\small]
  \coordinate (p0) at (-2.35,-1.10);
  \coordinate (p1) at (-2.65,.25);
  \coordinate (p2) at (-1.35,1.55);
  \coordinate (p3) at (1.35,1.55);
  \coordinate (p4) at (2.65,.25);
  \coordinate (p5) at (2.35,-1.10);
  \fill[black!4] (p0)--(p1)--(p2)--(p3)--(p4)--(p5)--cycle;
  \draw[arc] (p0)--(p1)--(p2)--(p3)--(p4)--(p5);
  \draw[boundary] (p5)--(p0);
  \path (p0)--(p1) node[midway,left=5pt] {$\ell_1$};
  \path (p1)--(p2) node[midway,above left=3pt] {$\ell_2$};
  \path (p2)--(p3) node[midway,above=5pt] {$\ell_3$};
  \path (p3)--(p4) node[midway,above right=3pt] {$\ell_4$};
  \path (p4)--(p5) node[midway,right=5pt] {$\ell_5$};
  \coordinate (rho) at (0,-1.10);
  \foreach \i/\j in {0/1,1/2,2/3,3/4,4/5}{
    \draw[auxiliary] (rho)--($(p\i)!.5!(p\j)$);
  }
  \node[figlabel] at (1.44,.55) {$P_p$};
  \node[figlabel] at (0,-.34) {$D$};
  \node[blackvertex,label={[font=\small]below=5pt:$p$}]
    at (0,-1.10) {};
  \foreach \i in {0,...,5}{\node[whitevertex] at (p\i) {};}
  \node at (0,-1.91) {$q(P_p)=\nu_D(p)=5,\qquad\lambda(P_p)=4$};
\end{tikzpicture}
\caption{The polygon $P_p$ containing the $\bullet$-point $p$
and the five endpoint occurrences of its dual fan.
The dissection-side occurrences are listed clockwise after
the boundary side. The dashed segments are the portions
of the dual arcs $\ell_1^*,\ldots,\ell_5^*$ from $p$ to
their crossings with $\ell_1,\ldots,\ell_5$; the dual arcs
continue beyond these crossings.}
\label{fig:elementary-polygon}
\end{figure}

By duality,
\begin{equation}
\label{eq:dual-valency}
 q(P_p)=\nu_D(p).
\end{equation}
Define the polygon weight by
\begin{equation}
\label{eq:polygon-weight}
 \lambda(P_p):=q(P_p)-1=\nu_D(p)-1.
\end{equation}
The exclusion of the zero surface implies $q(P_p)=\nu_D(p)\geq1$.
Figure~\ref{fig:elementary-polygon} illustrates the correspondence.

\subsection{Global dimension from the dual dissection}

The global-dimension formula can be read directly from the endpoint fans of
$D$; see, for example, \cite[Theorem~5.10]{LiuGaoHuang} and
\cite[Corollary~2.28]{Chang}.

\begin{proposition}
\label{prop:dual-global-dimension}
Let $\Delta$ be an admissible dissection of a marked surface as above, and
put $D:=\Delta^*$.  Then $A=A(\Delta)$ is a finite-dimensional gentle algebra
of finite global dimension, and
\begin{equation}
\label{eq:surface-global-dimension}
 \gldim A
 =\max_{p\in M^\bullet}\bigl(\nu_D(p)-1\bigr)
 =\max_{P\in\mathcal P(\Delta)}\lambda(P).
\end{equation}
\end{proposition}

Put $M_i^\bullet:=M^\bullet\cap\partial_iS$ and
$m_i:=|M_i^\bullet|$.  Write
$M_i^\bullet=\{p_{i,1},\ldots,p_{i,m_i}\}$, and let
$P_{i,j}:=P_{p_{i,j}}$ be the polygon in $\mathcal P(\Delta)$ containing
$p_{i,j}$.
The geometric interpretation of the AAG invariant gives
the multiset $\AG(A)=[(m_1,n_1),\ldots,(m_b,n_b)]$, where
\begin{equation}
\label{eq:surface-formulas}
\begin{aligned}
 n_i&=\sum_{p\in M_i^\bullet}\bigl(\nu_D(p)-1\bigr)
     =\sum_{j=1}^{m_i}\lambda(P_{i,j}).
\end{aligned}
\end{equation}
Indeed, the boundary $\partial_iS$ specifies one AAG orbit, whose forbidden
threads are the full fans based at the points of $M_i^\bullet$.
Their lengths are $\nu_D(p)-1$, with trivial threads contributing zero;
see \cite{AvellaAlaminosGeiss,OpperPlamondonSchroll}.
For the one-vertex algebra $k$, we use the convention
$\AG(k)=[(2,0)]$, as supplied by its surface model:
a single dissection arc in a disc with two points of each colour.
Both dual endpoint degrees are one.

\subsection{Euler identities}

We shall repeatedly use the following standard numerical consequences of the
admissible-dissection model.

\begin{proposition}
\label{prop:Euler-identities}
Let $A$ be a gentle algebra of finite global dimension represented by
an admissible dissection $\Delta$ of a surface of genus $g$ with
$b$ boundary components.  Put
\[
 m:=\sum_{i=1}^b m_i
\]
and let $E:=|\Delta|$.  Then
\begin{equation}
\label{eq:Euler-identities}
 E=m+2g+b-2,
 \qquad
 \sum_{i=1}^b n_i=m+4g+2b-4.
\end{equation}
In particular, the AAG list determines the genus.
\end{proposition}

\begin{proof}
The equality
\[
 |\Delta|=|M^{\bullet}|+2g+b-2
\]
is the unpunctured case of the admissible-dissection arc-count formula
\cite[Proposition~3.12]{AmiotPlamondonSchroll}; here
$|M^{\bullet}|=m$.
Every dual arc contributes two endpoint occurrences.  Therefore
\[
 \sum_{i=1}^b n_i
 =\sum_{p\in M^\bullet}\bigl(\nu_D(p)-1\bigr)
 =2|D|-m=m+4g+2b-4,\qedhere
\]
where $|D|=|\Delta|=m+2g+b-2.$
\end{proof}

\subsection{The complete line-field invariant}

We recall the complete derived invariant of a finite-dimensional
gentle algebra; see
\cite{LekiliPolishchuk,AmiotPlamondonSchroll}.
A grading on $\Delta$ determines a homotopy class of line fields
on $S$.  For the ordinary algebra $A=A(\Delta)$, we denote the
class associated with the zero grading on $\Delta$ by $\eta_A$,
or simply by $\eta$ when $A$ is fixed.  We use the construction
of this class from the dual dissection $D=\Delta^*$ given in
\cite{AmiotPlamondonSchroll}.

Let $(S,M,\eta)$ be a graded marked surface.
For an oriented immersed closed curve $\gamma$, write
$w_\eta(\gamma)$ for its winding number relative to $\eta$,
using the convention of
\cite[Definition~3.2]{AmiotPlamondonSchroll}.
This number is invariant under regular homotopy.
With the boundary orientation fixed above, the AAG interpretation
of boundary winding for $\eta=\eta_A$ gives
\begin{equation}
\label{eq:boundary-winding}
 w_i:=w_{\eta}(\partial_iS)=m_i-n_i.
\end{equation}
Equivalently, the endpoint degrees of $D$ give
\[
 w_i=2m_i-\sum_{p\in M_i^\bullet}\nu_D(p).
\]
Thus the AAG invariant determines all boundary winding numbers.

When $g=0$, no additional internal datum is required: the boundary winding
numbers determine the mapping-class-group orbit of the line field.

\begin{figure}[htbp]
\centering
\begin{tikzpicture}[x=1cm,y=1cm,font=\small]
  \fill[black!3] (-2.55,-1.75) rectangle (2.55,1.75);
  \draw[seam] (-2.55,-1.75) rectangle (2.55,1.75);
  \draw[identification] (-1.85,1.75)--(-1.05,1.75);
  \draw[identification] (-1.85,-1.75)--(-1.05,-1.75);
  \draw[identification] (-2.55,-1.18)--(-2.55,-.45);
  \draw[identification] (2.55,-1.18)--(2.55,-.45);
  \node[above=4pt] at (-1.45,1.75) {$s_1$};
  \node[below=4pt] at (-1.45,-1.75) {$s_1$};
  \node[left=5pt] at (-2.55,-.85) {$s_2$};
  \node[right=5pt] at (2.55,-.85) {$s_2$};
  \draw[arc] (-2.55,.52)--(2.55,.52);
  \draw[flow] (.42,.52)--(1.04,.52);
  \draw[arc] (-.58,-1.75)--(-.58,1.75);
  \draw[flow] (-.58,.88)--(-.58,1.38);
  \node[above=4pt] at (1.65,.52) {$\alpha$};
  \node[left=5pt] at (-.58,1.2) {$\beta$};
  \fill[white] (1.12,-.72) circle (.48);
  \draw[boundary] (1.12,-.72) circle (.48);
  \draw[flow] (1.12,-.72) ++(5:.48)
    arc[start angle=5,end angle=105,radius=.48];
  \node at (1.12,-.72) {$\partial_1 S$};
  \node at (0,-2.47) {$g=1,\quad b=1,\quad [\alpha]\cdot[\beta]=1$};
\end{tikzpicture}
\caption{A torus with one boundary component in a fundamental
square. Equally labelled dashed sides are identified in the
arrow directions. The boundary component $\partial_1 S$
is oriented so that $S$ lies to its right.
Marked points and the line field are omitted.}
\label{fig:line-field-data}
\end{figure}
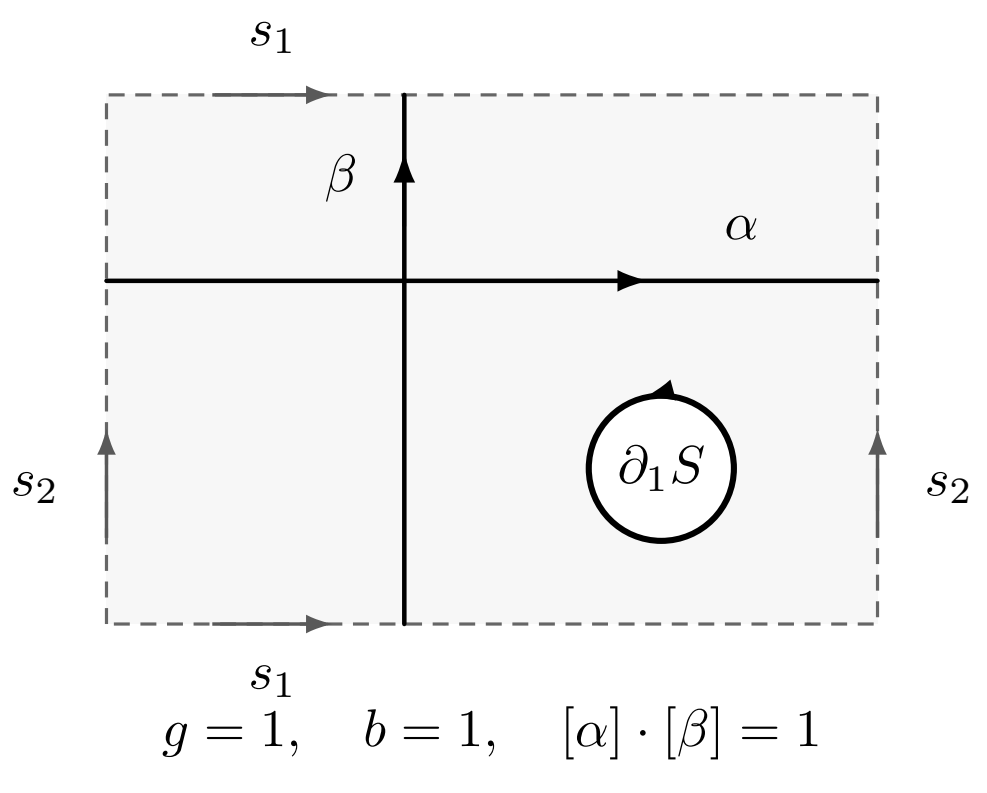

When $g=1$, choose oriented non-separating simple closed curves
$\alpha,\beta$ intersecting once positively, whose classes modulo the
boundary classes form a symplectic basis; see
\cref{fig:line-field-data} for the one-boundary case.  Define
\begin{equation}
\label{eq:genus-one-gcd}
 d(\eta):=
 \gcd\bigl(
   w_\eta(\alpha),w_\eta(\beta),
   w_1+2,\ldots,w_b+2
 \bigr),
\end{equation}
where $\gcd(0,\ldots,0)=0$.

When $g\geq2$, there are two successive pieces of information.
\begin{enumerate}[label=\textup{(\roman*)},leftmargin=2.6em]
\item
Set $\sigma(\eta)=0$ if every oriented closed curve has even winding number,
and set $\sigma(\eta)=1$ otherwise.

\item
If $\sigma(\eta)=0$ and
\[
 w_i\equiv2\pmod4
 \qquad(1\leq i\leq b),
\]
let $\overline S$ be the closed surface obtained by capping all boundary
components of $S$ with discs.  The line field induces a quadratic form
$q_\eta$ on $H_1(S;\mathbb Z/2)$, characterised on oriented simple closed
curves by
\[
 q_\eta([\gamma])=\frac12 w_\eta(\gamma)+1\pmod2.
\]
The boundary congruences give $q_\eta([\partial_iS])=0$, so $q_\eta$
descends to a quadratic form $\overline q_\eta$ on
$H_1(\overline S;\mathbb Z/2)$. We write
$\Arf(\eta):=\Arf(\overline q_\eta)$ for the classical Arf invariant
of the spin structure determined by $\eta$, and record this invariant as well;
see \cite[Lemma~7.5]{AmiotPlamondonSchroll}.  In all other higher-genus
cases no Arf datum is needed. 
\end{enumerate}

\begin{remark}
\label{rem:even-line-field}
The parity condition $\sigma(\eta)=0$ says precisely that $\eta$ is
orientable. Choose a representative of $\eta$ and a Riemannian metric
on $S$, and let $\pi:TS\to S$ be the bundle projection. The double cover
$$
 L=\{\,v\in TS\mid \|v\|=1,\ \operatorname{span}(v)=\eta(\pi(v))\,\}
   \longrightarrow S
$$
is classified by a class $a\in H^1(S;\ZZ/2)$, and
$w_\eta(\gamma)\equiv\langle a,[\gamma]\rangle\pmod2$.  In particular the
parity of $w_\eta$ is a homomorphism $H_1(S;\ZZ/2)\to\ZZ/2$, so
$\sigma(\eta)=0$ if and only if the curves of a basis of $H_1(S;\ZZ/2)$
have even winding number.  Equivalently, $\eta$ is even if and only if
it is induced by a nowhere-vanishing vector field, which is then unique
up to homotopy; compare the spin-structure interpretation of the Arf
datum in~\textup{(ii)}.
\end{remark}

We collect the internal data into the symbol
\begin{equation}
\label{eq:tau-line-field}
 \tau(\eta):=
 \begin{cases}
  \varnothing,
     & g=0,\\
  d(\eta),
     & g=1,\\
  \bigl(\sigma(\eta),\Arf(\eta)\bigr),
     & g\geq2,\ \sigma(\eta)=0,\ \text{and }
       w_i\equiv2\pmod4\text{ for all }i,\\
  \sigma(\eta),
     & g\geq2\text{ otherwise}.
 \end{cases}
\end{equation}

The following result is proved in \cite[Theorems~6.1 and~7.4]{AmiotPlamondonSchroll}.

\begin{proposition}
\label{prop:complete-invariant}
Let $A$ and $B$ be gentle algebras of finite global dimension.
Then
\[
 \Db(\modu A)\simeq\Db(\modu B)
\]
if and only if, after a permutation of the boundary components, the two
algebras have the same AAG invariant and
\begin{equation}
\label{eq:complete-invariant-tau}
 \tau(\eta_A)=\tau(\eta_B).
\end{equation}
\end{proposition}

\section{Prescribed-degree quadrangulations}
\label{sec:quadrangulations}

We construct auxiliary bipartite quadrangulations on closed oriented
surfaces to realise prescribed endpoint degrees and boundary assignments.
The resulting marked surfaces and admissible dissections will be obtained
from these quadrangulations through a ribbon-graph construction.

\begin{figure}[htbp]
\centering
\begin{tikzpicture}[
  x=.95cm,y=.95cm,font=\small,
  tr edge/.style={line width=.65pt,line cap=round,line join=round},
  tr chosen/.style={line width=1.25pt,line cap=round},
  tr diagonal/.style={draw=black!65,densely dashed,line width=.6pt},
  tr graph/.style={line width=.95pt,line cap=round},
  tr boundary/.style={line width=1.1pt},
  tr arc/.style={line width=.85pt,line cap=round},
  tr dual/.style={draw=black!65,densely dashed,line width=.65pt},
  tr white/.style={circle,draw=black,fill=white,line width=.65pt,
    inner sep=0pt,minimum size=4.7pt},
  tr black/.style={circle,fill=black,inner sep=0pt,minimum size=4.4pt},
  tr title/.style={font=\small,align=center},
  tr label/.style={font=\footnotesize,inner sep=1pt},
  tr flow/.style={line width=.7pt,-{Latex[length=2mm]}},
  tr pairing/.style={line width=.8pt,-{Latex[length=1.8mm]}},
  tr corner/.style={line width=.6pt}
]

\begin{scope}[xshift=-3.4cm]
  \node[tr title] at (0,2.65) {(a) $\mathcal Q$ on the sphere};
  \coordinate (qa) at (2.1,0);
  \coordinate (qb) at (.5,0);
  \coordinate (qc) at (-1.1,-1.05);
  \coordinate (qd) at (-1.65,0);
  \coordinate (qi) at (0,1.2);
  \coordinate (qo) at (0,-1.3);

  \draw[tr diagonal] (qa)--(qb)--(qd)--(qc);
  \draw[tr diagonal] (qd)
    .. controls (-2.6,2.45) and (3.1,2.45) .. (qa);

  \draw[tr chosen] (qi)--(qa);
  \draw[tr chosen] (qi)--(qb);
  \draw[tr edge]   (qi)--(qd);
  \draw[tr edge]   (qo)--(qa);
  \draw[tr edge]   (qo)--(qb);
  \draw[tr chosen] (qo)--(qc);
  \draw[tr edge]   (qd)--(qo);
  \draw[tr chosen] (qd)
    .. controls (-3.25,-1.65) and (-1.55,-2.65) .. (qo);

  \foreach \p in {qa,qb,qc,qd}{
    \node[tr black] at (\p) {};
  }
  \node[tr label,anchor=west] at (2.25,0) {$v_1$};
  \node[tr label] at (.76,-.20) {$v_2$};
  \node[tr label,anchor=east] at (-1.28,-1.02) {$v_3$};
  \node[tr label,anchor=east] at (-1.87,.10) {$v_4$};
  \node[tr white,
    label={[font=\footnotesize]above=4pt:$f_{\mathrm{in}}$}]
    at (qi) {};
  \node[tr white,
    label={[font=\footnotesize]below right=3pt:$f_{\mathrm{out}}$}]
    at (qo) {};

  \node[tr label] at (1.02,.34) {$Q_1$};
  \node[tr label] at (-.45,.40) {$Q_2$};
  \node[tr label] at (-1.04,-1.61) {$Q_3$};
  \node[tr label] at (-1.25,1.00) {$Q_4$};
\end{scope}

\draw[tr flow] (-.48,.25)--(.55,.25);
\node[tr label,align=center] at (.04,.80) {black\\diagonals};

\begin{scope}[xshift=3.4cm]
  \node[tr title] at (0,2.65)
    {(b) $G$ with distinguished corners};
  \coordinate (ga) at (.15,1.45);
  \coordinate (gb) at (1.5,.10);
  \coordinate (gc) at (-2.05,-.75);
  \coordinate (gd) at (-1.2,.10);
  \draw[tr graph] (ga)--(gb)--(gd)--(ga);
  \draw[tr graph] (gd)--(gc);

  \draw[tr corner] ($(ga)+(225:.28)$)
    arc[start angle=225,end angle=315,radius=.28];
  \draw[tr corner] ($(gb)+(135:.32)$)
    arc[start angle=135,end angle=180,radius=.32];

  \draw[tr corner] ($(gc)+(45:.25)$)
    arc[start angle=45,end angle=405,radius=.25];

  \draw[tr corner] ($(gd)+(45:.28)$)
    arc[start angle=45,end angle=225,radius=.28];

  \foreach \p in {ga,gb,gc,gd}{
    \node[tr black] at (\p) {};
  }

  \node[tr label,anchor=south] at (.15,1.63) {$v_1$};
  \node[tr label,anchor=west] at (1.68,.10) {$v_2$};
  \node[tr label,anchor=east] at (-2.42,-.75) {$v_3$};
  \node[tr label] at (-.90,-.18) {$v_4$};

  \node[tr label] at (.15,.91) {$c_1$};
  \node[tr label] at (.94,.31) {$c_2$};
  \node[tr label] at (-2.05,-1.25) {$c_3$};
  \node[tr label] at (-1.72,.60) {$c_4$};

  \node[tr label] at (.98,.98) {$e_1$};
  \node[tr label] at (.23,-.16) {$e_2$};
  \node[tr label] at (-1.32,-.65) {$e_3$};
  \node[tr label] at (-.83,1.02) {$e_4$};

  \node[tr label] at (.22,.49) {$R_{\mathrm{in}}$};
  \node[tr label] at (1.43,1.88) {$R_{\mathrm{out}}$};
\end{scope}

\draw[tr flow] (3.4,-1.67)--(3.4,-3.23);
\node[tr label,anchor=west] at (3.58,-2.14) {reconstruction};

\begin{scope}[yshift=-7.25cm,xshift=-4cm]
  \node[tr title] at (4,4.15)
    {(c) A rectangular model of the annulus};
  \coordinate (bl) at (0,0);
  \coordinate (bm) at (4,0);
  \coordinate (br) at (8,0);
  \coordinate (tl) at (0,3.2);
  \coordinate (tm) at (4,3.2);
  \coordinate (tr) at (8,3.2);

  \fill[black!8] (tl)--(tr)--(bm)--cycle;
  \fill[black!4] (bl)--(tl)--(bm)--cycle;
  \fill[black!4] (tr)--(br)--(bm)--cycle;
  \fill[white] (tl)--(tm)
    .. controls (3.3,2.1) and (1.7,2.1) .. (tl)--cycle;

  \draw[tr boundary] (bl)--(br);
  \draw[tr boundary] (tl)--(tr);

  \draw[tr arc] (bl)--(tl);
  \draw[tr arc] (br)--(tr);
  \draw[tr arc] (bm)--(tl);
  \draw[tr arc] (tl)
    .. controls (1.7,2.1) and (3.3,2.1) .. (tm);
  \draw[tr arc] (bm)--(tr);

  \draw[tr pairing] (0,.84)--(0,1.23);
  \draw[tr pairing] (8,.84)--(8,1.23);

  \draw[tr dual] (2,0)--(0,1.6);
  \draw[tr dual] (8,1.6)--(6,0);
  \draw[tr dual] (2,0)
    .. controls (1.8,1.8) and (4.8,1.6) .. (6,3.2);
  \draw[tr dual] (2,3.2)
    .. controls (2.2,1.3) and (4.5,2.6) .. (6,3.2);
  \draw[tr dual] (6,3.2)--(6,0);

  \node[tr label,anchor=east] at (-.17,1.94) {$\ell_1$};
  \node[tr label,anchor=west] at (8.17,1.94) {$\ell_1$};
  \node[tr label,fill=white] at (1.08,2.34) {$\ell_2$};
  \node[tr label,fill=white] at (3.10,2.55) {$\ell_3$};
  \node[tr label,fill=white] at (6.83,2.26) {$\ell_4$};

  \node[tr white,label={[font=\footnotesize]below=4pt:$y_1$}]
    at (bl) {};
  \node[tr white,label={[font=\footnotesize]below=4pt:$y_2$}]
    at (bm) {};
  \node[tr white,label={[font=\footnotesize]below=4pt:$y_1$}]
    at (br) {};
  \node[tr white,label={[font=\footnotesize]above=4pt:$y_4$}]
    at (tl) {};
  \node[tr white,label={[font=\footnotesize]above=4pt:$y_3$}]
    at (tm) {};
  \node[tr white,label={[font=\footnotesize]above=4pt:$y_4$}]
    at (tr) {};

  \node[tr black,label={[font=\footnotesize]below=4pt:$p_2$}]
    at (2,0) {};
  \node[tr black,label={[font=\footnotesize]below=4pt:$p_1$}]
    at (6,0) {};
  \node[tr black,label={[font=\footnotesize]above=4pt:$p_3$}]
    at (2,3.2) {};
  \node[tr black,label={[font=\footnotesize]above=4pt:$p_4$}]
    at (6,3.2) {};

  \node[tr label] at (1.02,1.35) {$P_{v_2}$};
  \node[tr label] at (1.35,2.77) {$P_{v_3}$};
  \node[tr label] at (4.03,1.10) {$P_{v_4}$};
  \node[tr label] at (6.98,1.35) {$P_{v_1}$};

  \draw[tr arc] (1.12,-.72)--(1.80,-.72);
  \node[tr label,anchor=west] at (1.95,-.72) {$\Delta$};
  \draw[tr dual] (4.05,-.72)--(4.73,-.72);
  \node[tr label,anchor=west] at (4.88,-.72) {$D=\Delta^*$};

\end{scope}
\end{tikzpicture}

\caption{An example of the passage from a quadrangulation
to a dissected annulus.} \label{fig:quadrangulation-transfer-example}
\end{figure}
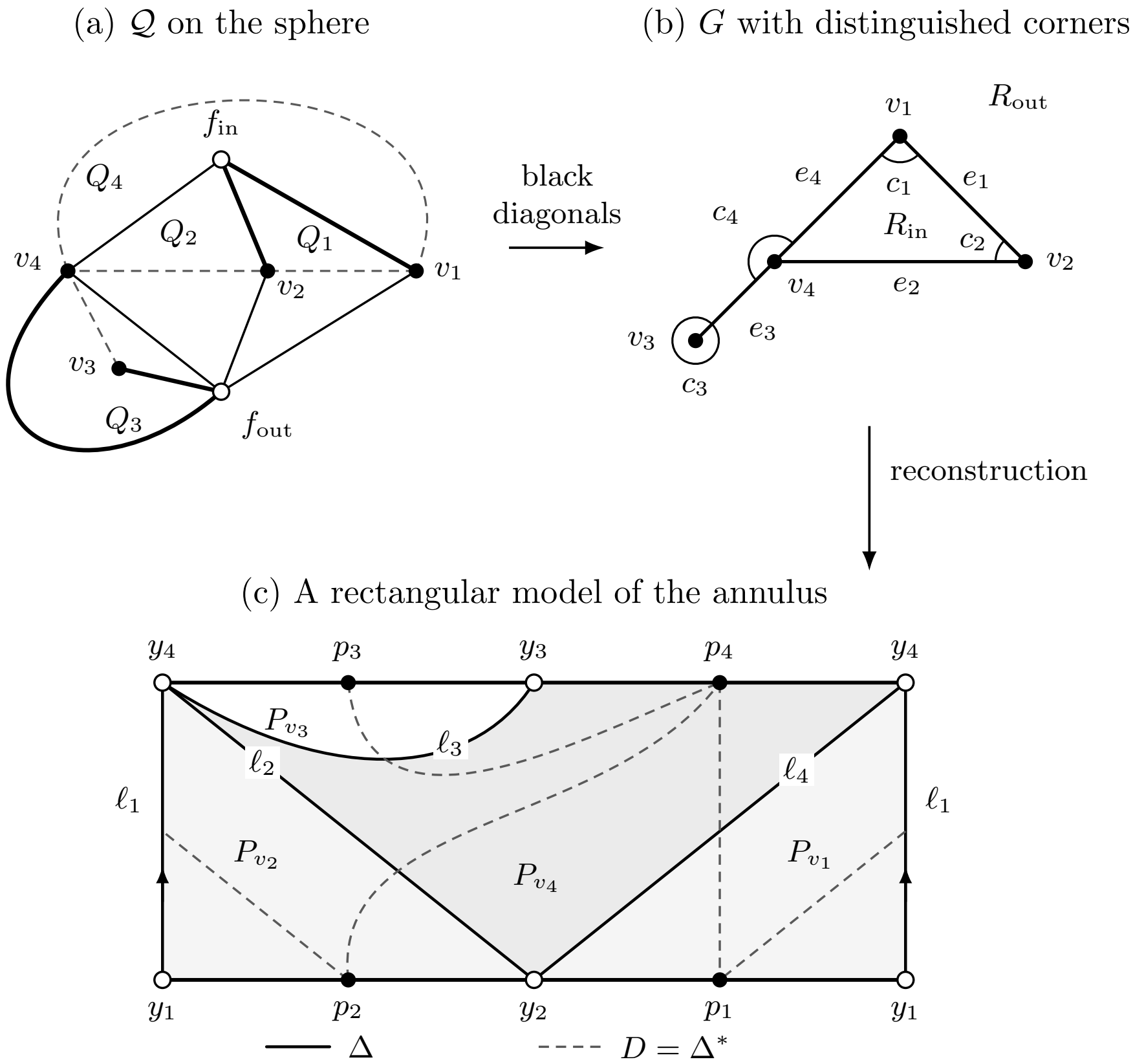

A \emph{cellular quadrangulation} is an embedding of a finite connected
graph in a closed oriented surface such that every complementary
component is an open disc and every face boundary walk has length four.
We use the convention that parallel edges are allowed, and
a face boundary walk may repeat vertices or edges.  Thus a quadrilateral
face need not have four distinct vertices or four distinct edges.
Degrees and face lengths are counted with incidence multiplicity.
All quadrangulations considered below are bipartite, with black vertices
$v$ and auxiliary white vertices $f_i$.  
See \cref{fig:quadrangulation-transfer-example} \textup{(a)} for an example.

The construction used in the following two subsections proceeds in the direction
\[
  \mathcal Q
  \xrightarrow{\text{black diagonals}}
  (G,(c_v)_{v\in V})
  \xrightarrow{\text{reconstruction}}
  (S,M,\Delta),
  \qquad D=\Delta^*.
\]
In the first step, each quadrilateral face of $\mathcal Q$
contributes a diagonal joining its two black corner occurrences.
These diagonals form an embedded ribbon graph $G$.
Each auxiliary white vertex $f_i$ lies in a complementary
disc of $G$, and these discs correspond to the boundary
components of its oriented thickening $S$.
For each black vertex $v$, a distinguished edge $vf_i$
selects the corner at $v$ facing the corresponding disc.
Under the reconstruction described in
\cref{subsec:ribbon-reconstruction}, $v$ becomes a marked
$\bullet$-point on the corresponding boundary component,
and $G$ yields the dual dissection $D=\Delta^*$. 

The general construction is developed in stages below.
Before diving into the abstract construction, let us first consider an example to build some intuition.

\begin{example}
Consider the spherical bipartite quadrangulation $\mathcal Q$ in
\cref{fig:quadrangulation-transfer-example} (a).
Its black vertices are $v_1,v_2,v_3,v_4$. The white vertex
$f_{\mathrm{in}}$ is joined once to each of $v_1,v_2,v_4$;
$f_{\mathrm{out}}$ is joined once to each of $v_1,v_2,v_3$
and twice to $v_4$. Its four face boundary walks are
\[
\begin{aligned}
 Q_1&=(f_{\mathrm{in}},v_1,f_{\mathrm{out}},v_2),&
 Q_2&=(f_{\mathrm{in}},v_2,f_{\mathrm{out}},v_4),\\
 Q_3&=(f_{\mathrm{out}},v_3,f_{\mathrm{out}},v_4),&
 Q_4&=(f_{\mathrm{in}},v_4,f_{\mathrm{out}},v_1).
\end{aligned}
\]

The diagonals connecting black vertices form the ribbon graph in
panel~\textup{(b)}, with edges
\[
 e_1=v_1v_2,\qquad e_2=v_2v_4,\qquad
 e_3=v_3v_4,\qquad e_4=v_4v_1.
\]
It consists of a triangle with a pendant edge at $v_4$.
Choose the distinguished incidences shown by thicker edges in
panel~\textup{(a)}: $v_1f_{\mathrm{in}}$, $v_2f_{\mathrm{in}}$,
$v_3f_{\mathrm{out}}$, and the outermost edge $v_4f_{\mathrm{out}}$.
They give corners $c_1,c_2$ facing $R_{\mathrm{in}}$ and
$c_3,c_4$ facing $R_{\mathrm{out}}$; specifically, $c_4$ lies
between $e_4$ and $e_3$. The corner at the leaf $v_3$ is unique.
The oriented thickening is an annulus. Panel~\textup{(c)}
represents the reconstructed dissection $\Delta$ in a rectangular
model of this annulus.
The dashed arcs show the
dual dissection $D=\Delta^*$. 
\end{example}

\subsection{Alternating polygons}

We begin with a planar degree-distribution lemma.

\begin{lemma}
\label{lem:alternating-polygon}
Let
$x_1,y_1,x_2,y_2,\ldots,x_p,y_p$ be the cyclically ordered vertex
occurrences of an alternating $2p$-gon, where $p\geq2$, with $x_j$ black
and $y_j$ white.  Let $r_1,\ldots,r_p\in\NN$ satisfy
\begin{equation}
\label{eq:alternating-demand}
 \sum_{j=1}^p r_j=p-2.
\end{equation}
Then the polygon has a noncrossing bipartite quadrangulation in which exactly
$r_j$ internal diagonals are incident with the black occurrence $x_j$.
\end{lemma}
\begin{figure}[htbp]
\centering
\begin{tikzpicture}[x=1cm,y=1cm,font=\small]
  \begin{scope}[xshift=-3.65cm]
    \foreach \i in {0,...,9}{\coordinate (z\i) at ({90-36*\i}:1.7);}
    \fill[black!12] (z9)--(z0)--(z1)--(z2)--cycle;
    \draw[arc] (z0)--(z1)--(z2)--(z3)--(z4)--(z5)--(z6)--(z7)--(z8)--(z9)--cycle;
    \draw[diagonal] (z9)--(z2);
    \foreach \i in {0,2,4,6,8}{\node[blackvertex] at (z\i) {};}
    \foreach \i in {1,3,5,7,9}{\node[whitevertex] at (z\i) {};}
    \node[above=4pt] at (z0) {$x_j$};
    \node[above right=4pt] at (z1) {$y_j$};
    \node[right=5pt] at (z2) {$x_{j+1}$};
    \node[above left=4pt] at (z9) {$y_{j-1}$};
    \node at (-.11,2.53) {$r_j=0$};
    \node[font=\scriptsize] at (.8,.04) {$r_{j+1}>0$};
    \node[figlabel,font=\scriptsize] at (.05,1.26) {ear};
    \node[panel] at (0,-2.12) {(a) Insert $y_{j-1}x_{j+1}$};
  \end{scope}
  \draw[flow] (-.67,0)--(.67,0);
  \begin{scope}[xshift=3.65cm]
    \foreach \i in {0,...,9}{\coordinate (z\i) at ({90-36*\i}:1.7);}
    \draw[arc] (z2)--(z3)--(z4)--(z5)--(z6)--(z7)--(z8)--(z9);
    \draw[diagonal] (z9)--(z2);
    \foreach \i in {2,4,6,8}{\node[blackvertex] at (z\i) {};}
    \foreach \i in {3,5,7,9}{\node[whitevertex] at (z\i) {};}
    \node[right=5pt] at (z2) {$x_{j+1}$};
    \node[above left=4pt] at (z9) {$y_{j-1}$};
    \node[panel] at (0,2.53) {$x_j,y_j$ removed};
    \node[panel] at (0,-.17) {$r'_{j+1}=r_{j+1}-1$};
    \node[panel] at (0,-2.12) {(b) The remaining $2(p-1)$-gon};
  \end{scope}
\end{tikzpicture}
\caption{The ear-removal step, illustrated for $p=5$.  The new diagonal
cuts off the quadrilateral $y_{j-1},x_j,y_j,x_{j+1}$.
Deleting the occurrences $x_j,y_j$ preserves the alternating colours;
only the residual demand at $x_{j+1}$ decreases.}
\label{fig:alternating-ear}
\end{figure}
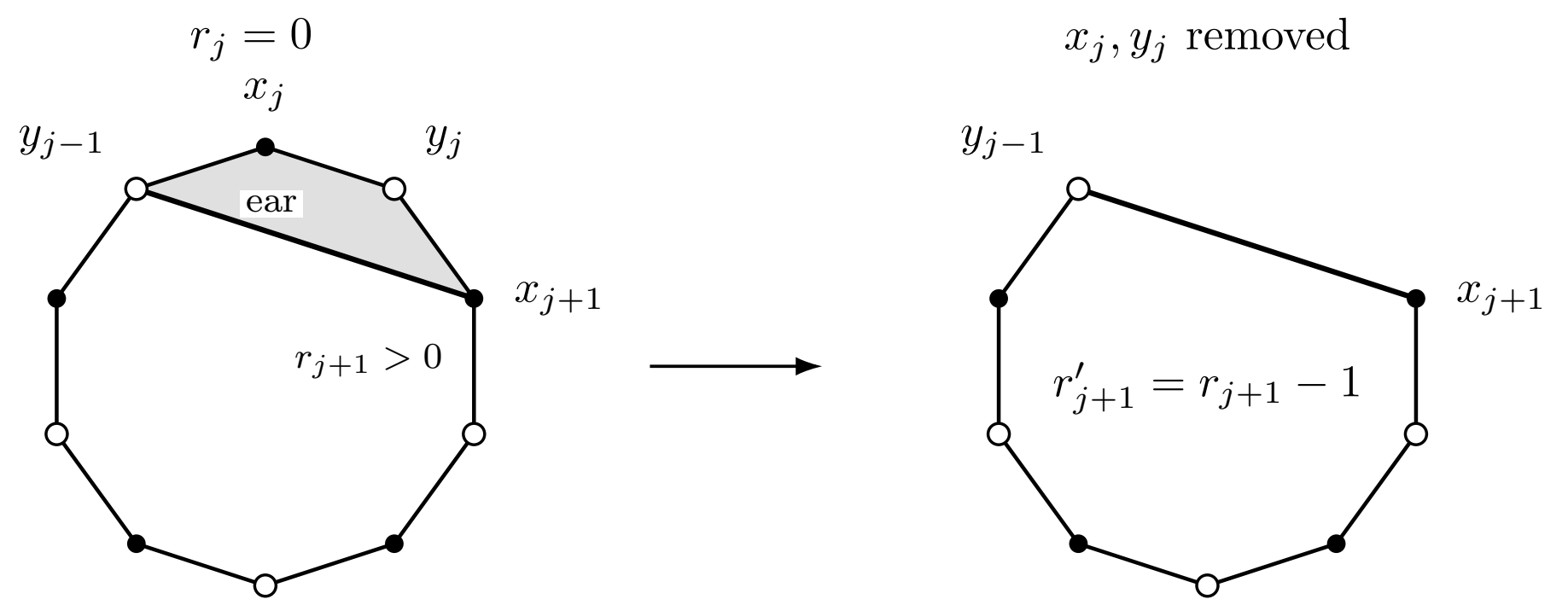

\begin{proof}
We argue by induction on $p$.  For $p=2$, the polygon is already a
quadrilateral and \eqref{eq:alternating-demand} forces all $r_j$ to vanish.
Suppose $p>2$. Since the sum is smaller than $p$, the cyclic list $(r_1,\ldots,r_p)$ contains at least one zero entry and at least one
positive entry. Therefore after cyclic
relabelling, assume that $r_j=0$ and $r_{j+1}>0$.  Draw the diagonal from
$y_{j-1}$ to $x_{j+1}$, with indices read cyclically.  It cuts off the
quadrilateral with boundary occurrences
$y_{j-1},x_j,y_j,x_{j+1}$, as in \cref{fig:alternating-ear}.

Delete the zero entry $r_j$ and replace $r_{j+1}$ by $r_{j+1}-1$.  The
resulting nonnegative cyclic list has sum $p-3$, so the induction hypothesis
quadrangulates the remaining alternating $2(p-1)$-gon.  Restoring the removed
quadrilateral adds precisely the new diagonal at $x_{j+1}$ and none at $x_j$,
which gives the required incidences.
\end{proof}

\subsection{A one-face bipartite core}

We recall the graph-theoretic terminology needed to construct a
cellular one-face embedding.  Let $G$ be a connected finite undirected
multigraph. Write $V(G)$ and $E(G)$ for the vertex and edge sets of $G$,
respectively.
Its \emph{maximum orientable genus} $\gamma_{\max}(G)$ is the largest
integer $h\geq0$ for which $G$ admits a cellular embedding in the closed
connected oriented surface of genus $h$.  Here cellular means that
every component of the complement is an open disc.

A \emph{spanning tree} $T$ of $G$ is a connected acyclic subgraph
containing every vertex of $G$.  The spanning subgraph $G-E(T)$,
obtained by deleting the edges of $T$ but retaining all vertices,
is called the \emph{cotree} relative to $T$. Note that a cotree need not be connected
or acyclic.  Write
\[
  \beta(G):=|E(G)|-|V(G)|+1
\]
for the cycle rank of $G$.  Since $|E(T)|=|V(G)|-1$, this is also the
number of edges in $G-E(T)$.  Let
$\operatorname{odd}(G-E(T))$ denote the number of connected components
of $G-E(T)$ having an odd number of edges.  Isolated vertices have
zero edges and do not contribute to this number.

Xuong's maximum-genus theorem \cite{Xuong} states that
\begin{equation}
\label{eq:Xuong}
  \gamma_{\max}(G)
  =
  \frac12\left(
    \beta(G)-\min_T\operatorname{odd}(G-E(T))
  \right),
\end{equation}
where the minimum is taken over all spanning trees of $G$.

For a cellular embedding of $G$ in a closed oriented surface of genus $h$ with $F$ faces, Euler's formula gives
\[
  F=1+\beta(G)-2h.
\]
Consequently, if $\beta(G)=2g$ and there is a spanning tree $T$ such
that every component of $G-E(T)$ has an even number of edges, then
\eqref{eq:Xuong} gives $\gamma_{\max}(G)=g$, and an embedding attaining
this genus has exactly one face.  This is the consequence used in
the next lemma.

\begin{lemma}
\label{lem:capacity-core}
Fix integers $b\geq1$ and $g\geq0$. Let
$V=V_1\amalg\cdots\amalg V_b$ be a disjoint union of finite sets, where $|V_i|=m_i\geq1$, and put
$m:=|V|$.  Let $q_v\geq1$ be integers satisfying
\begin{equation}
\label{eq:core-degree-sum}
 \sum_{v\in V}q_v=2(m+2g+b-2).
\end{equation}  Put $N:=m+b$.  There is a connected bipartite multigraph
$G$ with black vertex set $V$ and white vertex set
$\{f_1,\ldots,f_b\}$ such that:
\begin{enumerate}[label=\textup{(\roman*)}]
\item $G$ contains a distinguished edge $vf_i$ for every $v\in V_i$;
\item $\deg_G(v)\leq q_v$ for every black vertex $v$;
\item $|E(G)|=N-1+2g$;
\item $G$ has a cellular one-face embedding in the closed oriented surface
of genus $g$.
\end{enumerate}
\end{lemma}

\begin{proof}
We construct $G$ in the following three steps.  We first build a
tree $T$ on the prescribed vertex set, then add $2g$ edges to obtain
$G$, and finally apply Xuong's formula \eqref{eq:Xuong} using
the spanning tree $T\subseteq G$.

\textbf{Step 1: joining the distinguished stars.}
For every black vertex $v\in V_i$, include an edge $vf_i$ and declare it
distinguished.  These edges form $b$ disjoint stars.
If $b=1$, take this star as the tree $T$.
Suppose now that $b\geq2$.

We join the stars by repeatedly adding an edge from a black vertex
$v$ whose current degree is less than $q_v$ to a white vertex $f_j$
in a different connected component.  Each new edge merges two
components without creating a cycle and preserves all black-degree
bounds.
This operation is possible as long as more than one component
remains.  Indeed, every component contains a white vertex.
Moreover, each added edge increases the degree of exactly one
black vertex by one, so after $k<b-1$ additions the total remaining
black-degree capacity is
\[
  \sum_{v\in V}(q_v-1)-k
  =m+4g+2b-4-k
  \geq b-1-k>0,
\]
by \eqref{eq:core-degree-sum}, since $b\geq2$ and $m\geq b$.
Thus there is still a black vertex whose degree is less than $q_v$,
and it can be joined to a white vertex in another component.
After $b-1$ additions, we obtain a tree $T$ on all $N$ prescribed
vertices, containing every distinguished edge and satisfying
$\deg_T(v)\leq q_v$ for all $v\in V$.

\textbf{Step 2: adding the cotree edges.}
The tree $T$ has $N-1$ edges, each incident with exactly one black
vertex.  Hence its total remaining black-degree capacity is
\begin{equation}
\label{eq:residual-capacity}
  \sum_{v\in V}\bigl(q_v-\deg_T(v)\bigr)
  =2(N+2g-2)-(N-1)
  =N+4g-3.
\end{equation}
Since $\sum_{v\in V}q_v>0$, equation~\eqref{eq:core-degree-sum}
gives $N+2g\geq3$. Hence
\[
  N+4g-3=2g+(N+2g-3)\geq2g.
\]
We may therefore add $2g$ edges without exceeding any
black-degree bound, taking $f_1$ as their common white
endpoint. Parallel edges are allowed.
Let $G$ be the resulting graph.  Since we have added edges but
no vertices, $T$ is a spanning tree of $G$.  The graph $G$ is
connected and bipartite, contains every distinguished edge,
and satisfies
\[
  \deg_G(v)\leq q_v \quad(v\in V),
  \qquad
  |E(G)|=N-1+2g.
\]

\textbf{Step three: obtaining a one-face embedding.}
The edges of the cotree $G-E(T)$ are precisely the $2g$ edges
added in the preceding step.
If $g=0$, the cotree $G-E(T)$ has no edges.
If $g>0$, all its $2g$ edges lie in one connected component,
because they share the endpoint $f_1$; every other component
is an isolated vertex.  Thus
\[
  \operatorname{odd}(G-E(T))=0.
\]
Since $\beta(G)=2g$, Xuong's formula \eqref{eq:Xuong} gives
$\gamma_{\max}(G)=g$.  A cellular embedding attaining this genus
has
\[
  F=1+\beta(G)-2g=1
\]
face, as required.
\end{proof}

\subsection{Completion to a quadrangulation}

The remaining degree capacities in Lemma \ref{lem:capacity-core} (ii) can be filled inside the unique face.

\begin{theorem}
\label{thm:prescribed-quadrangulation}
Under the hypotheses of \cref{lem:capacity-core}, the closed oriented
surface of genus $g$ has a cellular bipartite quadrangulation $\cQ$ whose
black vertex set is $V$, whose white vertex set is
$\{f_1,\ldots,f_b\}$, and whose black degrees satisfy
\begin{equation}
\label{eq:quadrangulation-degrees}
 \deg_{\cQ}(v)=q_v\qquad(v\in V).
\end{equation}
Moreover, $\cQ$ contains a distinguished edge $vf_i$ for every $v\in V_i$.
\end{theorem}

\begin{proof}
Take the one-face embedding of the graph $G$ supplied by
\cref{lem:capacity-core}, and put $p:=|E(G)|=N-1+2g$.
We have $p\geq2$; otherwise $N=2$ and $g=0$, while
\eqref{eq:core-degree-sum} would give $\sum_vq_v=0$,
contrary to $q_v\geq1$.

Represent the unique face by a closed disc $B$ with $2p$
boundary sides following its boundary walk. Each edge of
$G$ labels two boundary sides. The side pairings determined
by the embedding present the ambient surface as a quotient
of $B$. The boundary vertices of $B$ alternate between black
and white; several may represent the same vertex of $G$.

For each black vertex $v$, put
$h_v:=q_v-\deg_G(v)\geq0$. Since the sum of the black
degrees of $G$ is $p$, we have
\begin{equation}
\label{eq:residual-demand}
  \sum_{v\in V}h_v=2(N+2g-2)-p=p-2.
\end{equation}

For each black vertex $v$, let
\[
  x_{v,1},\ldots,x_{v,\deg_G(v)}
\]
be its boundary occurrences on $\partial B$. Choose
nonnegative integers
\[
  r_{v,1},\ldots,r_{v,\deg_G(v)}
  \qquad\text{such that}\qquad
  \sum_s r_{v,s}=h_v.
\]
Since
\[
  \sum_{v,s}r_{v,s}=\sum_vh_v=p-2,
\]
\cref{lem:alternating-polygon} gives a noncrossing bipartite
quadrangulation of $B$ with exactly $r_{v,s}$ internal
diagonals incident with $x_{v,s}$.

Pass to the quotient determined by the original side
pairings. Only boundary points of $B$ are identified,
so each quadrilateral interior maps homeomorphically
onto an open disc, and distinct interiors have disjoint
images. Thus the images of the boundary sides and the
new diagonals form a cellular bipartite quadrangulation
$\cQ$ of the original surface.

At each black vertex $v$, the new diagonals contribute
altogether $h_v$ additional incidences, so
\[
  \deg_{\cQ}(v)=\deg_G(v)+h_v=q_v.
\]
The original core $G$ remains a subgraph of $\cQ$, so all
distinguished edges are retained.
\end{proof}

\begin{example}
The construction of \cref{thm:prescribed-quadrangulation}
is illustrated in \cref{fig:core-to-quadrangulation}
for $g=1$, $m=b=1$, and $q_v=4$. The unique face is
represented by a hexagonal disc, and one diagonal
subdivides this disc into two quadrilaterals.
\begin{figure}[htbp]
\centering
\begin{tikzpicture}[
  x=1cm,y=1cm,font=\small,
  arc/.style={line width=.85pt,line cap=round,line join=round},
  diagonal/.style={line width=1.2pt,line cap=round},
  whitevertex/.style={circle,draw=black,line width=.65pt,fill=white,
    inner sep=0pt,minimum size=4.7pt},
  blackvertex/.style={circle,fill=black,inner sep=0pt,minimum size=4.1pt},
  panel/.style={align=center},
  figlabel/.style={fill=white,inner sep=1.2pt},
  flow/.style={line width=.7pt,-{Latex[length=2.1mm,width=1.45mm]}},
  rotation/.style={draw=black!60,line width=.65pt,
    -{Latex[length=1.8mm,width=1.25mm]}},
  pairedge/.style={line width=.75pt,postaction={decorate},
    decoration={markings,mark=at position .54 with
      {\arrow{Latex[length=1.8mm,width=1.2mm]}}}}
]
  \begin{scope}[xshift=-3.65cm]
    \node[panel] at (0,2.35) {(a) The core rotation system};
    \foreach \c in {-1.25,1.25}{
      \begin{scope}[xshift=\c cm]
        \foreach \ang/\edge in {90/1,210/2,330/3}{
          \draw[arc] (0,.35)--++(\ang:.78);
          \node at ($(0,.35)+(\ang:1.06)$) {$e_{\edge}$};
        }
        \draw[rotation] ($(0,.35)+(65:.9)$)
          arc[start angle=65,end angle=-5,radius=.9];
      \end{scope}
    }
    \node[blackvertex] at (-1.25,.35) {};
    \node[anchor=south east,inner sep=0pt] at (-1.37,.49) {$v$};
    \node[whitevertex] at (1.25,.35) {};
    \node[anchor=south east,inner sep=0pt] at (1.13,.49) {$f_1$};
    \node[panel] at (0,-1.00)
      {$\operatorname{rot}(v)=\operatorname{rot}(f_1)=(e_1\ e_3\ e_2)$};
    \node[panel] at (0,-1.55)
      {$|V(G)|=2,\quad |E(G)|=3,\quad |F(G)|=1$};
    \node[panel] at (0,-2.03)
      {$\deg_G(v)=3,\qquad q_v=4$};
  \end{scope}

  \draw[flow] (-.65,.35)--(.65,.35);
  \node[panel,font=\scriptsize] at (0,.87) {unique face};

  \begin{scope}[xshift=3.65cm]
    \node[panel] at (0,2.35) {(b) Add the residual diagonal};
    \foreach \i in {0,...,5}{
      \coordinate (h\i) at ({90-60*\i}:1.65);
    }
    \fill[black!4] (h0)--(h1)--(h2)--(h3)--(h4)--(h5)--cycle;
    \foreach \i/\j in {0/1,2/1,2/3,4/3,4/5,0/5}{
      \draw[pairedge] (h\i)--(h\j);
    }
    \path (h0)--(h1) node[midway,above right=3pt] {$e_1$};
    \path (h1)--(h2) node[midway,right=5pt] {$e_2$};
    \path (h2)--(h3) node[midway,below right=3pt] {$e_3$};
    \path (h3)--(h4) node[midway,below left=3pt] {$e_1$};
    \path (h4)--(h5) node[midway,left=5pt] {$e_2$};
    \path (h5)--(h0) node[midway,above left=3pt] {$e_3$};
    \draw[diagonal] (h0)--(h3);
    \node[figlabel] at (0,.10) {$d$};
    \node at (.72,.10) {$Q_1$};
    \node at (-.72,.10) {$Q_2$};
    \foreach \i in {0,2,4}{\node[blackvertex] at (h\i) {};}
    \foreach \i in {1,3,5}{\node[whitevertex] at (h\i) {};}
    \node[panel] at (0,-2.03)
      {$h_v=1,\qquad \deg_{\mathcal Q}(v)=4$};
  \end{scope}
\end{tikzpicture}
\caption{A genus-one example with $m=b=1$ and $q_v=4$.
Equal half-edge labels in \textup{(a)} are paired, and the cyclic
orders are read clockwise, as indicated by the grey arrows.
The core has one face, with $e_1$ chosen as the distinguished edge.
The hexagonal disc in \textup{(b)} represents this face, with side
pairings indicated by arrows. The diagonal $d$ subdivides the disc
into two quadrilaterals. The resulting quadrangulation has two
vertices, four edges, and two faces, hence genus one.}
\label{fig:core-to-quadrangulation}
\end{figure}
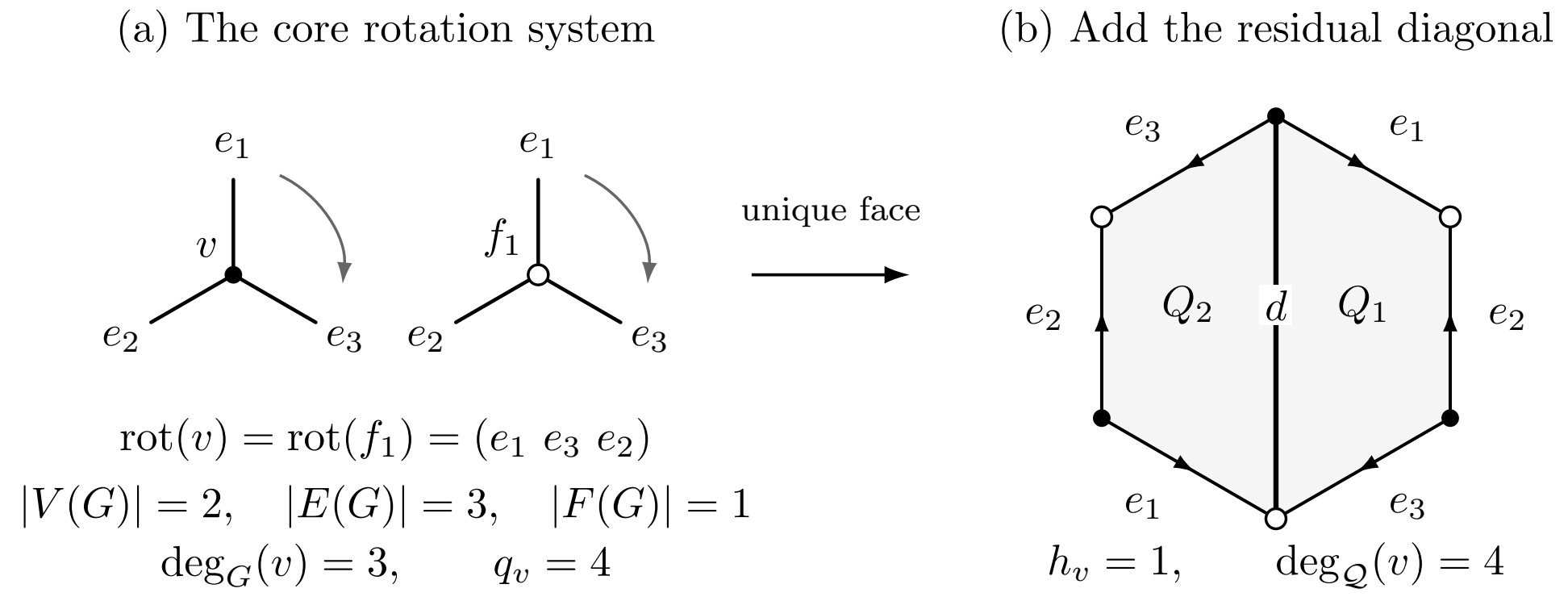
\end{example}

\section{All-genus realisation by admissible \texorpdfstring{$\circ$}{o}-dissections}
\label{sec:arc-realization}

We now pass from the quadrangulations of
\cref{sec:quadrangulations} to admissible dissections through
ribbon graphs.  We first recall the standard reconstruction
used in both the realisation theorem and the finite formula
of \cref{sec:finite-formula}.
\subsection{Ribbon graphs and dual dissections}
\label{subsec:ribbon-reconstruction}

We recall the standard correspondence between marked ribbon
graphs and dissections; see \cite[Section~3]{SchrollSurvey}
for an overview. 
We use the polygon-dual version, whose vertices correspond to
forbidden threads; see \cite[Remark~1.9]{OpperPlamondonSchroll}
and the construction preceding and in the proof of
\cite[Theorem~3.11]{LekiliPolishchuk}.
We describe the constructions and record the properties needed
below without repeating their standard proofs.

\begin{definition}
\label{def:distinguished-corner-ribbon-graph}
A \emph{ribbon graph} is a finite undirected multigraph $G$
equipped with a cyclic order on the half-edge occurrences at
each vertex.  Loops and parallel edges are allowed.  Throughout,
we assume that $G$ is connected and has at least one edge.

\begin{enumerate}[label=\textup{(\arabic*)},leftmargin=2.6em]
\item \emph{Half-edge encoding.}
Write $H$ for the half-edge set.  Each edge has two distinct
half-edges.  Let $\iota$ be the
fixed-point-free involution exchanging the half-edges of each
edge, and let $\pi$ be the successor permutation in the
prescribed cyclic orders, read clockwise in oriented vertex
charts. The length of the
cycle at $v$ is $\deg_G(v)$, so a loop contributes two.
Connectedness of $G$ is equivalent to transitivity of
$\langle\iota,\pi\rangle$ on $H$.

\item \emph{Thickening and boundary cycles.}
The \emph{oriented thickening} of $G$ is obtained by replacing
vertices with oriented discs and edges with bands, attaching
the half-edges in the prescribed clockwise orders and gluing
the bands so that the disc orientations extend over them.
The \emph{boundary cycles} are the cycles of $\pi\iota$.  They correspond
bijectively to the boundary components of the thickening.

\item \emph{Distinguished corners.}
A \emph{corner} is the sector between consecutive half-edge
occurrences in the clockwise cyclic order at a vertex. Index the corner from $\pi^{-1}(h)$ to $h$ by its succeeding
half-edge $h$. The corner indexed by $h$ belongs to the
boundary cycle containing $h$.
A degree-one vertex has a single corner.
A \emph{ribbon graph with distinguished corners} is a ribbon
graph with one chosen corner $c_v$ at each vertex $v$.
It is \emph{boundary-admissible} if every boundary cycle
contains at least one distinguished corner.
\end{enumerate}
\end{definition}

\smallskip
Let $G$ be boundary-admissible, and let $S$ be its oriented
thickening. At each vertex $v$, choose a point $p_v$ in the
boundary interval corresponding to its distinguished corner,
and declare these to be the $\bullet$-points. Every boundary
component contains at least one $\bullet$-point by boundary
admissibility. For each open boundary interval $J$ between
consecutive $\bullet$-points, choose a $\circ$-point $y_J\in J$.
Set
\[
  M:=\{p_v\mid v\in V(G)\}\amalg\{y_J\mid J\}.
\]
The marked points alternate in colour on every boundary component.

We now construct the dissection $\Delta_G$ on the marked surface. In each edge band,
choose a \emph{cocore}, namely a properly embedded transverse
segment joining the two sides of the band that lie in
$\partial S$. Every interval $J$ contains at least one cocore
endpoint, since the boundary walk between consecutive
distinguished corners traverses at least one edge-band side.
Keeping the cocores unchanged away from the boundary, deform
them near $\partial S$ so that all endpoints lying in each $J$
meet at $y_J$. Arrange the terminal segments as a noncrossing
fan at $y_J$, preserving the order of their original endpoints
along $J$. Denote the resulting family of arcs by $\Delta_G$.

The standard reconstruction makes $\Delta_G$ an admissible
dissection of $(S,M)$, as recorded below. Write $P_v$ for
the polygon of $\Delta_G$ containing $p_v$.

For the reconstructed dissection, put
\[
  D_G:=\Delta_G^*,\qquad
  A_G:=A(\Delta_G),\qquad
  \AG(G):=\AG(A_G).
\]
Geometrically, a representative of $D_G$ is obtained from the
embedded core graph $G\subset S$ by moving each vertex $v$
to $p_v$ through its distinguished corner, while keeping
the edge-band cores unchanged.

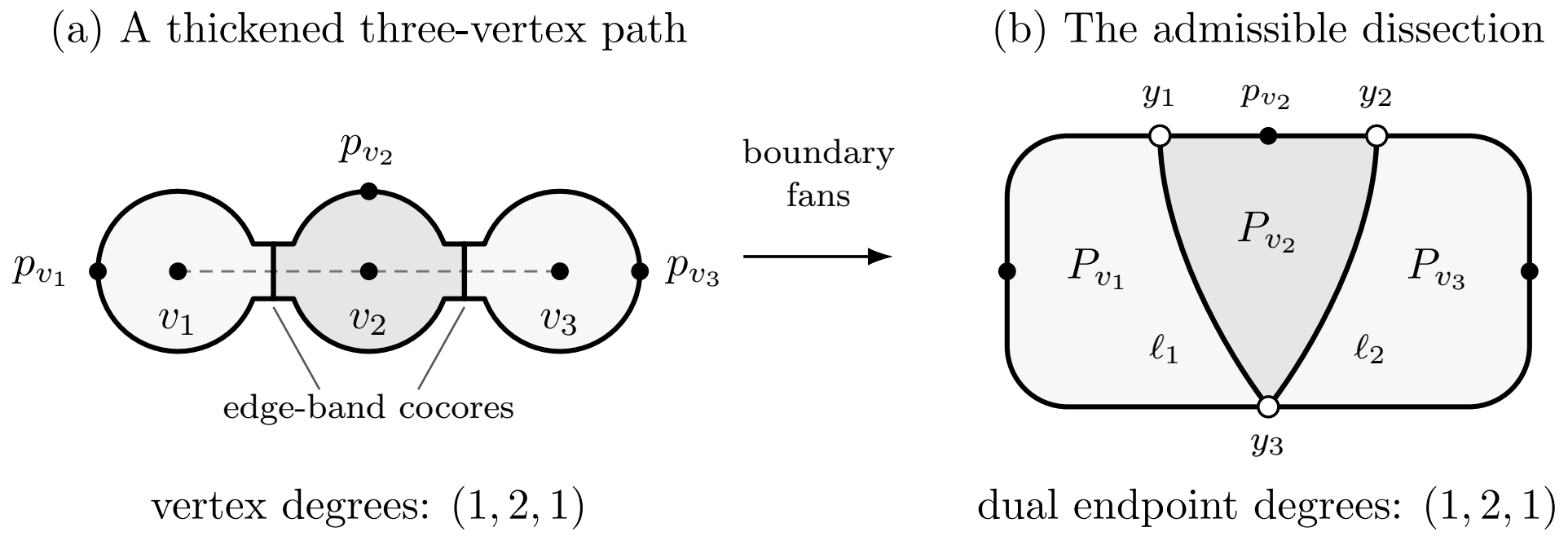
\begin{figure}[htbp]
\centering
\begin{tikzpicture}[x=1cm,y=1cm,font=\small]
  \begin{scope}[xshift=-3.65cm]
    \node[panel] at (0,1.95) {(a) A thickened three-vertex path};
    \pgfmathsetmacro{\bh}{.65*sin(20)}
    \pgfmathsetmacro{\dx}{.65*cos(20)}
    \def\ribbonoutline{
      ({-1.55+\dx},\bh)--(-\dx,\bh)
      arc[start angle=160,end angle=20,radius=.65]
      --({1.55-\dx},\bh)
      arc[start angle=160,end angle=-160,radius=.65]
      --(\dx,-\bh)
      arc[start angle=-20,end angle=-160,radius=.65]
      --({-1.55+\dx},-\bh)
      arc[start angle=-20,end angle=-340,radius=.65]--cycle}
    \fill[black!3] \ribbonoutline;
    \begin{scope}
      \clip \ribbonoutline;
      \fill[black!10] (-.775,-.8) rectangle (.775,.8);
    \end{scope}
    \draw[boundary] \ribbonoutline;
    \draw[auxiliary] (-1.55,0)--(1.55,0);
    \foreach \x/\i in {-1.55/1,0/2,1.55/3}{
      \node[spinevertex] at (\x,0) {};
      \node[below=5pt] at (\x,0) {$v_{\i}$};
    }
    \foreach \x in {-.775,.775}{
      \draw[diagonal] (\x,-\bh)--(\x,\bh);
      \draw[line width=.55pt] (\x-.045,-\bh)--(\x+.045,-\bh);
      \draw[line width=.55pt] (\x-.045,\bh)--(\x+.045,\bh);
    }
    \node[blackvertex,label={[font=\small]left=5pt:$p_{v_1}$}] at (-2.2,0) {};
    \node[blackvertex,label={[font=\small]above=5pt:$p_{v_2}$}] at (0,.65) {};
    \node[blackvertex,label={[font=\small]right=5pt:$p_{v_3}$}] at (2.2,0) {};
    \node[panel,font=\scriptsize] at (0,-1.12) {edge-band cocores};
    \draw[leader] (-.4,-.96)--(-.775,-.29);
    \draw[leader] (.4,-.96)--(.775,-.29);
    \node[panel] at (0,-1.92) {vertex degrees: $(1,2,1)$};
  \end{scope}
  \draw[flow] (-.61,.12)--(.61,.12);
  \node[panel,font=\scriptsize] at (0,.8) {boundary\\fans};
  \begin{scope}[xshift=3.65cm]
    \node[panel] at (0,1.95) {(b) The admissible dissection};
    \coordinate (y1) at (-.88,1.10);
    \coordinate (y2) at (.88,1.10);
    \coordinate (y3) at (0,-1.10);
    \fill[black!3,rounded corners=14pt] (-2.12,-1.10) rectangle (2.12,1.10);
    \fill[black!10] (y1)--(y2)
      .. controls (.87,.20) and (.27,-.76) .. (y3)
      .. controls (-.27,-.76) and (-.87,.20) .. (y1)--cycle;
    \draw[boundary,rounded corners=14pt] (-2.12,-1.10) rectangle (2.12,1.10);
    \draw[diagonal] (y1) .. controls (-.87,.20) and (-.27,-.76) .. (y3);
    \draw[diagonal] (y2) .. controls (.87,.20) and (.27,-.76) .. (y3);
    \node at (-1.38,.02) {$P_{v_1}$};
    \node at (0,.32) {$P_{v_2}$};
    \node at (1.38,.02) {$P_{v_3}$};
    \node[font=\scriptsize] at (-.83,-.63) {$\ell_1$};
    \node[font=\scriptsize] at (.83,-.63) {$\ell_2$};
    \node[blackvertex] at (-2.12,0) {};
    \node[blackvertex,label={[font=\scriptsize]above=5pt:$p_{v_2}$}] at (0,1.10) {};
    \node[blackvertex] at (2.12,0) {};
    \foreach \i in {1,2}{\node[whitevertex,label={[font=\scriptsize]above=5pt:$y_{\i}$}] at (y\i) {};}
    \node[whitevertex,label={[font=\scriptsize]below=5pt:$y_3$}] at (y3) {};
    \node[panel] at (0,-1.92) {dual endpoint degrees: $(1,2,1)$};
  \end{scope}
\end{tikzpicture}
\caption{Ribbon reconstruction for a path with three vertices.
The boundary black point $p_{v_i}$ records the distinguished corner at $v_i$;
the interior black points $v_i$ are vertices of the spine, not marked points.
The two lower cocore endpoints are gathered at $y_3$.  The resulting disc
has three polygons, each with one black point, and
$\nu_{D_G}(p_{v_i})=q(P_{v_i})=\deg_G(v_i)$.}
\label{fig:ribbon-reconstruction}
\end{figure}

Let $\eta_G:=\eta_{A_G}$ be the line field associated with the
zero grading on $\Delta_G$, constructed using $D_G$ as in
\cref{sec:surface-data}.  Identify the boundary cycles
$F_1,\ldots,F_b$ of $G$ with the boundary components of $S$,
and set
\[
  V_i(G):=\{v\in V(G)\mid p_v\in F_i\},
  \qquad
  \degmax(G):=\max_{v\in V(G)}\deg_G(v).
\]
The construction is illustrated in \cref{fig:ribbon-reconstruction}.
The following proposition records this standard reconstruction
in our notation.  Each half-edge at $v$ corresponds to a
dissection-side occurrence of $P_v$, so the numerical identities
follow from \eqref{eq:dual-valency},
\cref{prop:dual-global-dimension}, and
\eqref{eq:surface-formulas}.  The assertion about the line field
uses the construction in
\cite[Section~3.3]{AmiotPlamondonSchroll}.

\begin{proposition}
\label{prop:ribbon-reconstruction}
The family $\Delta_G$ is an admissible $\circ$-dissection of $(S,M)$, with
\begin{equation}
\label{eq:ribbon-dual-degrees}
  \nu_{D_G}(p_v)=q(P_v)=\deg_G(v)
  \qquad(v\in V(G)).
\end{equation}
It is independent of the choices up to homotopy under the corresponding
identification of marked surfaces; in particular, $A_G$ and $\eta_G$ are
well defined up to isomorphism and homotopy, respectively.
Moreover, $A_G$ is a connected finite-dimensional gentle algebra of
finite global dimension, and
\begin{equation}
\label{eq:ribbon-AAG-gldim}
  \AG(G)=
  \left[
    \left(|V_i(G)|,\sum_{v\in V_i(G)}(\deg_G(v)-1)\right)
  \right]_{i=1}^b,
  \qquad
  \gldim A_G=\degmax(G)-1.
\end{equation}
\end{proposition}

\smallskip
Conversely, let $\Delta$ be an admissible $\circ$-dissection of
a connected marked surface $(S,M)$ as in \cref{sec:surface-data}.
Put $D:=\Delta^*$ and $A_\Delta:=A(\Delta)$, and let
$\eta:=\eta_{A_\Delta}$ be the canonical line field associated
with the zero grading on $\Delta$.

To construct $G_\Delta$, for each $\bullet$-point $p$, let $P$ be the polygon of
$\Delta$ containing it. Move $p$ slightly into the interior
of $P$, deforming the incident segments of $D$ locally
while preserving their order. This gives an embedded ribbon
graph $G_\Delta$, whose vertices correspond to the polygons
of $\Delta$ and whose edges correspond to the arcs of $D$,
hence also to those of $\Delta$. Its cyclic orders are read
clockwise with respect to the orientation of $S$. At each
vertex, distinguish the corner corresponding to the original
boundary position of $p$.

This construction and the reconstruction above are inverse
up to the corresponding identification of marked surfaces;
see \cite[Propositions~1.16 and~1.25]{OpperPlamondonSchroll}
and \cite[Proposition~3.13]{AmiotPlamondonSchroll}.
Their compatibility with canonical line fields follows from
\cite[Section~3.3]{AmiotPlamondonSchroll}.
We record these facts for use in \cref{sec:finite-formula}.

\begin{lemma}
\label{lem:polygon-dual-reconstruction}
The ribbon graph $G_\Delta$ is boundary-admissible, and its
oriented thickening is homeomorphic to $S$. Under the
corresponding identification of marked surfaces, the
reconstruction in \cref{prop:ribbon-reconstruction} yields
a dissection homotopic to $\Delta$. Consequently,
\begin{equation}
\label{eq:dual-reconstruction-algebra}
  A_{G_\Delta}\cong A_\Delta,
  \qquad
  \eta_{G_\Delta}\simeq\eta.
\end{equation}
\end{lemma}

\subsection{Realisation of prescribed endpoint degrees}

We now apply the quadrangulation theorem to the dual ribbon graph.
The black diagonals record the prescribed endpoint degrees, while the
distinguished corners record their boundary assignments.

\begin{theorem}
\label{thm:all-genus-realization}
Let $g\geq0$ and $b\geq1$, and let
$V=V_1\amalg\cdots\amalg V_b$, where $|V_i|=m_i\geq1$, and put
$m:=|V|$.  Assign an integer $\lambda_v\geq0$ to every $v\in V$, and suppose
that
\begin{equation}
\label{eq:all-genus-Euler-condition}
 \sum_{v\in V}\lambda_v=m+4g+2b-4.
\end{equation}
Then there is an admissible $\circ$-dissection $\Delta$ on a surface of
genus $g$ with $b$ boundary components whose dual $D=\Delta^*$ has
$\bullet$-points $p_v$, $v\in V$, with $p_v\in\partial_iS$ for $v\in V_i$
and
\begin{equation}
\label{eq:realised-lambda}
 \nu_D(p_v)-1=\lambda(P_v)=\lambda_v\qquad(v\in V),
\end{equation}
where $P_v$ is the polygon containing $p_v$.  The algebra $A(\Delta)$
is connected and has finite global dimension
$\gldim A(\Delta)=\max_{v\in V}\lambda_v$.
\end{theorem}

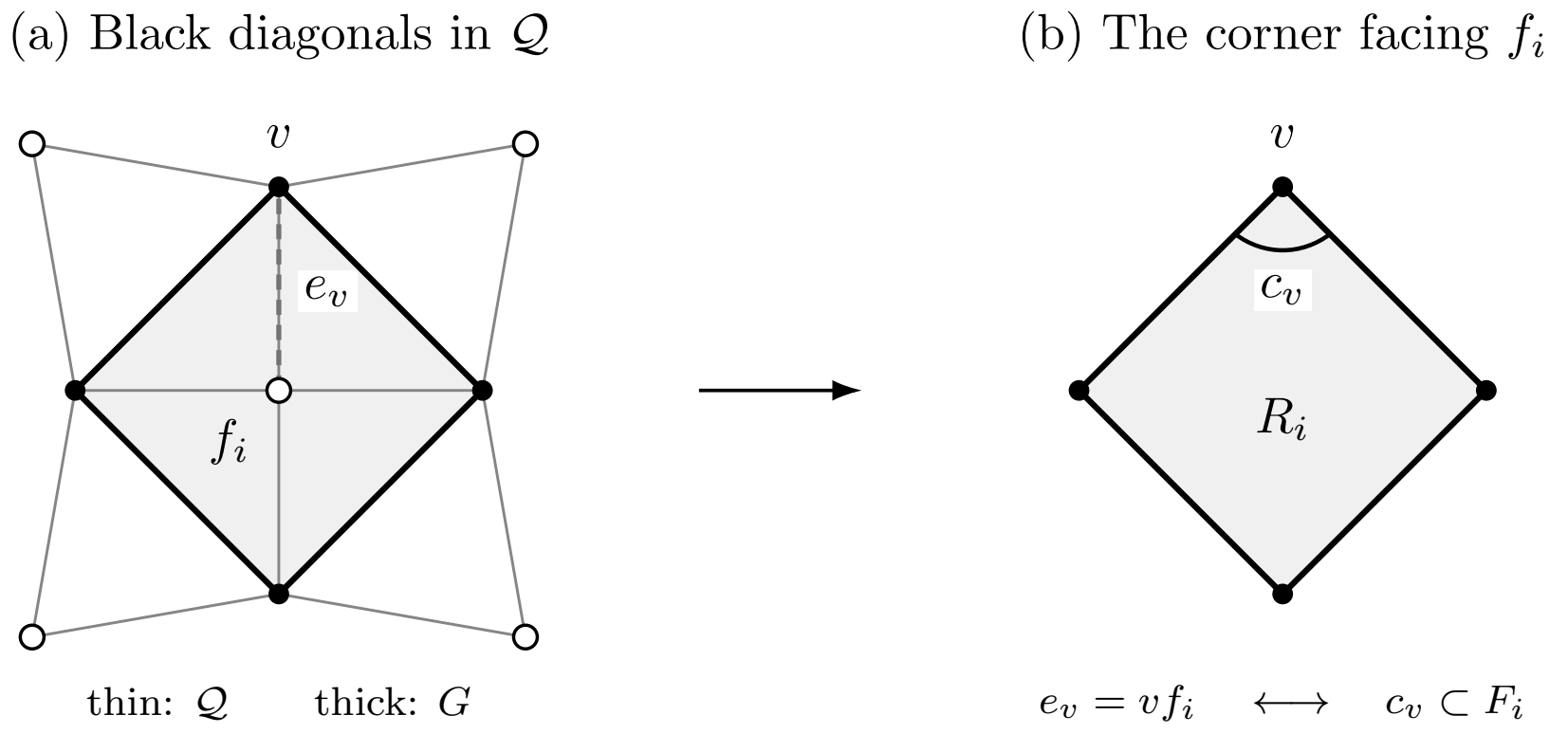
\begin{figure}[htbp]
\centering
\begin{tikzpicture}[x=1cm,y=1cm,font=\small]
  \begin{scope}[xshift=-3.5cm]
    \node[panel] at (0,2.47) {(a) Black diagonals in $\cQ$};

    \coordinate (f)  at (0,0);
    \coordinate (v1) at (0,1.42);
    \coordinate (v2) at (1.42,0);
    \coordinate (v3) at (0,-1.42);
    \coordinate (v4) at (-1.42,0);
    \coordinate (u1) at (1.72,1.72);
    \coordinate (u2) at (1.72,-1.72);
    \coordinate (u3) at (-1.72,-1.72);
    \coordinate (u4) at (-1.72,1.72);

    \fill[black!6] (v1)--(v2)--(v3)--(v4)--cycle;
    \draw[construction] (v1)--(u1)--(v2)--(u2)--(v3)--(u3)--(v4)--(u4)--cycle;
    \foreach \i in {1,...,4}{\draw[construction] (f)--(v\i);}
    \draw[diagonal] (v1)--(v2)--(v3)--(v4)--cycle;
    \draw[auxiliary,line width=1pt] (f)--(v1);

    \node[figlabel] at (.34,.69) {$e_v$};

    \foreach \i in {1,...,4}{
      \node[blackvertex] at (v\i) {};
      \node[whitevertex] at (u\i) {};
    }
    \node[whitevertex,label={[font=\small]below left=3pt:$f_i$}] at (f) {};
    \node[above=4pt] at (v1) {$v$};

    \node[panel,font=\scriptsize] at (0,-2.18) {thin: $\cQ$\qquad thick: $G$};
  \end{scope}

  \draw[flow] (-.57,0)--(.57,0);

  \begin{scope}[xshift=3.5cm]
    \node[panel] at (0,2.47) {(b) The corner facing $f_i$};

    \fill[black!6] (0,1.42)--(1.42,0)--(0,-1.42)--(-1.42,0)--cycle;
    \draw[diagonal] (0,1.42)--(1.42,0)--(0,-1.42)--(-1.42,0)--cycle;
    \foreach \ang in {90,0,270,180}{
      \node[blackvertex] at (\ang:1.42) {};
    }
    \node[above=4pt] at (0,1.42) {$v$};

    \draw[arc] (-.33,1.09) .. controls (-.13,.94) and (.13,.94) .. (.33,1.09);

    \node[figlabel] (cvlabel) at (0,.7) {$c_v$};

    \node at (0,-.20) {$R_i$};
    \node[panel,font=\scriptsize] at (0,-2.18)
      {$e_v=vf_i\quad\longleftrightarrow\quad c_v\subset F_i$};
  \end{scope}
\end{tikzpicture}
\caption{The black-diagonal construction, shown around an auxiliary white
vertex $f_i$.  Each quadrilateral contributes one edge of $G$.  Its
half-faces incident with $f_i$ form the complementary region $R_i$.
The distinguished incidence $e_v=vf_i$ selects the corner $c_v$ at $v$
facing that region, hence the corresponding boundary cycle $F_i$ of the
thickening.  In panel \textup{(b)}, the small curved arc marks the corner
$c_v$.  This is a local picture; the labels refer to incidences and
need not represent distinct vertices globally.}
\label{fig:black-diagonal-corners}
\end{figure}

\begin{proof}
Set $q_v:=\lambda_v+1$.  Condition~\eqref{eq:all-genus-Euler-condition} is equivalent to
\begin{equation}
\label{eq:q-degree-sum}
 \sum_{v\in V}q_v=2(m+2g+b-2).
\end{equation}
By \cref{thm:prescribed-quadrangulation} and \eqref{eq:q-degree-sum},
the closed genus-$g$ surface $\widehat S$ admits a cellular bipartite
quadrangulation $\cQ$ with black vertex set $V$, white vertices
$f_1,\ldots,f_b$, black degrees $q_v$, and a distinguished edge
$e_v$ joining $v$ to $f_i$ for each $v\in V_i$.

In each quadrilateral, draw the diagonal joining its two black corner
occurrences.  Let $G$ be the resulting embedded multigraph with vertex set
$V$. Let $F_{\cQ}$ be the number of quadrilateral faces.
Each quadrilateral corner at $v$ contributes one incident half-edge, so
\begin{equation}
\label{eq:black-diagonal-counts}
 \deg_G(v)=q_v,
 \qquad
 |E(G)|=F_{\cQ}=m+2g+b-2,
\end{equation}
where the edge count follows from $|E(\cQ)|=\sum_vq_v$ and
$4F_{\cQ}=2|E(\cQ)|$.

To see that $G$ is connected, fix a white vertex $f$ and list the
black endpoints of its incident edges in cyclic order, allowing
repetitions.  Consecutive entries are the two black corners of
a quadrilateral and are therefore joined by its new diagonal.
Thus all black neighbours of $f$ lie in one connected component
of $G$.  Since $\cQ$ is connected and bipartite, any two black
vertices can be joined by a path alternating between black and
white vertices.  Replacing each two-edge segment through a white
vertex by a walk along these diagonals gives a walk in $G$.
Therefore $G$ is connected.

Cut each quadrilateral along its black diagonal into two triangular
half-faces, each with exactly one white corner.  Across any original
edge of $\cQ$ incident with $f_i$, the two incident half-faces both
have white endpoint $f_i$.  Thus gluing along the original edges never joins half-faces
associated with different white vertices.  Conversely, the half-faces
incident with a fixed white vertex $f_i$ are joined in cyclic order
around $f_i$; the vertex $f_i$ and the interiors of the original edges
remain in the complement of $G$.  Since the black diagonals and black
vertices belong to $G$, the complement has exactly one connected
region $R_i$ for each white vertex $f_i$, and hence exactly $b$ regions.

To check that the regions are discs, choose a small closed regular
neighbourhood $N(G)$ and let $\widetilde R_i$ be the corresponding components
of $\widehat S\setminus\operatorname{int}N(G)$.  By
\eqref{eq:black-diagonal-counts},
\[
 \chi(N(G))=|V(G)|-|E(G)|=2-2g-b.
\]
Euler additivity along boundary circles gives
\[
 \sum_{i=1}^b\chi(\widetilde R_i)
 =\chi(\widehat S)-\chi(N(G))=b.
\]
Each $\widetilde R_i$ is a connected compact oriented surface with
nonempty boundary, so its Euler characteristic is at most one.  All $b$
terms must therefore equal one.  Hence every $\widetilde R_i$ is a disc,
and $G$ is cellular.  Its oriented thickening $S=N(G)$ has genus $g$ and
$b$ boundary components $F_i$, with $F_i$ facing $f_i$.  Equivalently,
one may identify $S$ with the closed surface with a small disc about each
white vertex removed.

For each $v\in V_i$, the half-edge of $e_v$ at $v$ lies in a
corner of $G$ facing the region $R_i$.  Denote this corner by
$c_v$ and declare it distinguished; see
\cref{fig:black-diagonal-corners}.  Then $c_v$ lies on $F_i$
exactly when $v\in V_i$.  Since $V_i\neq\varnothing$, every
boundary cycle contains a distinguished corner.  Thus $G$ is
boundary-admissible and $V_i(G)=V_i$.

Apply \cref{prop:ribbon-reconstruction}.  It gives a dissection $\Delta_G$
with dual $D_G$ and marked points $p_v\in F_i$
for $v\in V_i$, satisfying
\[
 \nu_{D_G}(p_v)=q(P_v)=\deg_G(v)=q_v.
\]
Taking $\Delta:=\Delta_G$ and $D:=D_G$ proves
\eqref{eq:realised-lambda}.  The same proposition gives connectedness and
finite global dimension of $A(\Delta)$, together with
$\gldim A(\Delta)=\max_v(q_v-1)=\max_v\lambda_v$.
\end{proof}

\subsection{The AAG-global dimension spectrum}

The only remaining numerical input is the following bounded-composition fact.

\begin{lemma}
\label{lem:bounded-distributions}
Let $m_i\geq1$ and $n_i\geq0$ for $1\leq i\leq b$, and put
$L:=\max_i\lceil n_i/m_i\rceil$ and $U:=\max_i n_i$.  For every integer
$d$ with $L\leq d\leq U$, there are integers $\lambda_{i,j}$,
$1\leq j\leq m_i$, such that
\begin{equation}
\label{eq:bounded-distributions}
 0\leq\lambda_{i,j}\leq d,
 \qquad
 \sum_{j=1}^{m_i}\lambda_{i,j}=n_i,
 \qquad
 \max_{i,j}\lambda_{i,j}=d.
\end{equation}
\end{lemma}

\begin{proof}
The inequality $d\geq L$ gives $n_i\leq m_id$ for every $i$, so $n_i$ can be
distributed among $m_i$ slots of capacity $d$.  Choose $i_0$ with
$n_{i_0}=U$ and set $\lambda_{i_0,1}:=d$.  Since
$0\leq n_{i_0}-d\leq(m_{i_0}-1)d$, the remainder fits into the other slots.
If $m_{i_0}=1$, these inequalities force $n_{i_0}=d$.  Distribute every other
$n_i$ among its slots with the same capacity.
\end{proof}

\begin{proof}[Proof of \cref{thm:intro-AAG-realization}]
Let $A'$ be a connected finite-dimensional gentle algebra of finite
global dimension with $\AG(A')=\AG(A)=[(m_i,n_i)]_{i=1}^b$.
Let $D$ be its dual dissection on the associated marked surface
$(S,M)$, with boundary components numbered according to these
AAG pairs.
Over $\partial_iS$, the $m_i$ nonnegative integers $\nu_D(p)-1$ have sum
$n_i$.  Their maximum is at least $\lceil n_i/m_i\rceil$, so
\eqref{eq:surface-global-dimension} gives $\gldim A'\geq L_A$.  Each is at most
$n_i\leq U_A$, which gives $\gldim A'\leq U_A$.

Conversely, fix $d\in[L_A,U_A]$.  Choose an array
$\lambda_{i,j}$ as in \cref{lem:bounded-distributions}.  By
\cref{prop:Euler-identities},
$\sum_{i,j}\lambda_{i,j}=\sum_i n_i=m+4g+2b-4$, where $g$ is the genus
determined by the AAG data.  Apply \cref{thm:all-genus-realization} with
$V_i:=\{v_{i,1},\ldots,v_{i,m_i}\}$ and
$\lambda_{v_{i,j}}:=\lambda_{i,j}$.  The resulting connected gentle algebra
has the same AAG list as $A$, and its global dimension is $d$ by
\eqref{eq:surface-global-dimension}.  This proves
\eqref{eq:intro-AAG-realization}.

Now let $B$ be any finite-dimensional algebra derived equivalent
to $A$. Derived equivalence preserves connectedness and
finiteness of global dimension \cite{KellerKrause,Rickard}.
Moreover, $B$ is gentle by \cite{SchroerZimmermann}, and
$\AG(B)=\AG(A)$ by \cite{AvellaAlaminosGeiss}.
The bounds just proved therefore give
\eqref{eq:intro-enclosure}.
\end{proof}

\begin{corollary}
\label{cor:gentle-boundedness}
Problem~\ref{prob:boundedness} has an affirmative answer for every connected
gentle algebra $A$ of finite global dimension.  One may take
\[
 b_A:=\max\{\gldim A-L_A,\,U_A-\gldim A\}.
\]
\end{corollary}

\begin{proof}
By \eqref{eq:intro-enclosure}, the global dimension of every algebra derived
equivalent to $A$ lies in $[L_A,U_A]$.
\end{proof}



\subsection{Derived classes determined by the AAG invariant}

Combining the realisation theorem with the sufficient conditions of
Lekili--Polishchuk \cite[Corollaries~3.14 and~3.16]{LekiliPolishchuk}, under which the
AAG data determine the line-field orbit, gives the following interval
formula.

\begin{corollary}
\label{thm:intro-derived-interval}
Let $A$ be a connected finite-dimensional gentle algebra of finite global
dimension, with AAG invariant
\[
 \AG(A)=[(m_i,n_i)]_{i=1}^b,
\] and let $g$ be the genus of its associated
surface. Suppose that one of the following conditions holds:
\begin{enumerate}[label=\textup{(\roman*)}]
\item $g=0$;
\item $g=1$ and $\gcd\{m_i-n_i+2\mid1\leq i\leq b\}=1$;
\item $g\geq2$ and $m_i-n_i$ is odd for at least one $i$.
\end{enumerate}
Then
\begin{equation}
\label{eq:intro-derived-interval}
 \GSpec(A)=[L_A,U_A].
\end{equation}
\end{corollary}

\begin{proof}
Let $B$ be a connected gentle algebra of finite global dimension with
$\AG(B)=\AG(A)$.  By \eqref{eq:boundary-winding}, $A$ and $B$ have the same
boundary winding numbers $w_i=m_i-n_i$ and the same numbers $m_i$ of
points of each colour on each boundary component.

If $g=0$, the boundary data determine the line-field orbit. If $g=1$, put $H:=\gcd_{1\leq i\leq b}(w_i+2)$.
By \eqref{eq:genus-one-gcd},
$d(\eta)=
 \gcd\bigl(w_\eta(\alpha),w_\eta(\beta),H\bigr).$
Thus the assumption $H=1$ implies
$d(\eta_A)=d(\eta_B)=1$. If $g\geq2$ and
some $w_i$ is odd, then the boundary curve itself has odd winding, so
$\sigma(\eta)=1$ for every such line field and no Arf invariant occurs.  The
mapping-class-group classification summarised in
\cref{prop:complete-invariant} therefore gives
$\Db(\modu B)\simeq\Db(\modu A)$ in each of the three cases.

Conversely, every finite-dimensional algebra derived equivalent
to $A$ is connected and gentle of finite global dimension,
with the same AAG invariant. Hence
$\AGSpec(A)=\GSpec(A)$ under the stated hypotheses, and
\eqref{eq:intro-derived-interval} follows from
\eqref{eq:intro-AAG-realization}.
\end{proof}

\begin{samepage}
\begin{example}
\label{ex:positive-genus-interval}
Consider the AAG data $[(3,4),(3,6)]$.  The Euler identities give $g=1$,
$b=2$, and $E=8$, while $L=2$ and $U=6$.  The boundary winding numbers are $w_1=-1$ and $w_2=-3$, so
$\gcd(w_1+2,w_2+2)=\gcd(1,-1)=1$.
The following polygon weights realise every integer in $[2,6]$:
\[
\begin{array}{c|c|c}
 d & (\lambda_{1,1},\lambda_{1,2},\lambda_{1,3})
   & (\lambda_{2,1},\lambda_{2,2},\lambda_{2,3})\\ \hline
 2&(2,1,1)&(2,2,2)\\
 3&(2,1,1)&(3,2,1)\\
 4&(2,1,1)&(4,1,1)\\
 5&(2,1,1)&(5,1,0)\\
 6&(2,1,1)&(6,0,0)
\end{array}
\]
Each row has boundary sums $4$ and $6$ and maximum entry $d$.
Theorem~\ref{thm:all-genus-realization} realises each row by an admissible
dissection.  The boundary gcd forces $d(\eta)=1$ for each resulting
algebra, so \cref{prop:complete-invariant} makes them all derived equivalent.
Their common ordinary spectrum is therefore $\{2,3,4,5,6\}$.
\end{example}
\end{samepage}
\section{A finite ribbon-graph formula for the spectrum}
\label{sec:finite-formula}

The preceding section determines the spectrum in an AAG fibre.
To restrict to a single derived-equivalence class, we must also fix
the internal line-field datum.  The reconstruction of
\cref{subsec:ribbon-reconstruction} associates an algebra and its
canonical line field with each boundary-admissible ribbon graph.
Formula~\eqref{eq:ribbon-AAG-gldim} expresses the AAG invariant
through vertex degrees and boundary assignments, and the global
dimension as $\degmax(G)-1$.  We now impose the remaining line-field
condition, prove the finite formula, and then explain how to compute it.

Let $A$ be a connected finite-dimensional gentle algebra of finite
global dimension, with AAG list \eqref{eq:intro-AAG-list}, genus $g$,
and internal line-field datum $\tau_A:=\tau(\eta_A)$.
The Euler identities fix the vertex and edge counts of the
corresponding ribbon graphs:
\begin{equation}
\label{eq:ME-fixed}
 m:=\sum_i m_i,
 \qquad
 E:=m+2g+b-2=\frac12\left(m+\sum_i n_i\right).
\end{equation}
Define $\mathfrak R(A)$ to be the set of isomorphism classes of
boundary-admissible ribbon graphs $G$ with distinguished corners
satisfying:
\begin{enumerate}[label=\textup{(\roman*)}]
\item $|V(G)|=m$ and $|E(G)|=E$;
\item the oriented thickening has genus $g$ and $b$ boundary components;
\item after numbering the boundary cycles,
\begin{equation}
\label{eq:R-boundary-data}
 |V_i(G)|=m_i,
 \qquad
 \sum_{v\in V_i(G)}(\deg_G(v)-1)=n_i
 \quad(1\leq i\leq b);
\end{equation}
\item $\tau(\eta_G)=\tau_A$.
\end{enumerate}
Here $V_i(G)$ consists of the vertices whose distinguished corners
lie on the $i$th boundary cycle, and $\eta_G$ is the canonical line
field supplied by \cref{prop:ribbon-reconstruction}.

\begin{proof}[Proof of \cref{thm:intro-finite-formula}]
Let $B$ be a finite-dimensional algebra derived equivalent to $A$.
Then $B$ is a connected gentle algebra of finite global dimension. It therefore admits a surface model with a dissection $\Delta_B$ such that $B\cong A(\Delta_B)$.  By
\cref{lem:polygon-dual-reconstruction}, the graph
$G_B:=G_{\Delta_B}$ is boundary-admissible and reconstructs both
$B$ and its canonical line field.  The Euler identities give the
counts in \eqref{eq:ME-fixed}, and \eqref{eq:ribbon-AAG-gldim}
gives the boundary conditions.  Derived equivalence preserves the
complete invariant, so $\tau(\eta_{G_B})=\tau_A$.  Hence
$G_B\in\mathfrak R(A)$ and
$\gldim B=\degmax(G_B)-1$.

Conversely, let $G\in\mathfrak R(A)$.
Proposition~\ref{prop:ribbon-reconstruction} gives a connected
finite-dimensional gentle algebra $A_G$ of finite global dimension.
Its AAG list equals that of $A$ by \eqref{eq:ribbon-AAG-gldim},
and its internal line-field datum is $\tau_A$ by condition~(iv).
Proposition~\ref{prop:complete-invariant} therefore gives
$\Db(\modu A_G)\simeq\Db(\modu A)$, and
$\gldim A_G=\degmax(G)-1$.  This proves both inclusions in
\eqref{eq:intro-finite-formula}.

Finally, $2E$ is fixed.  On a labelled $2E$-element half-edge set
there are only finitely many fixed-point-free involutions, vertex
permutations, and choices of distinguished corners.  Thus
$\mathfrak R(A)$ is finite.  The local winding rule and finite
procedure below give explicit tests for conditions~(i)--(iv),
and hence an explicit enumeration of $\mathfrak R(A)$.
\end{proof}

To evaluate the formula, we enumerate $\mathfrak R(A)$ and record
$\degmax(G)-1$ for each retained graph.  Conditions~(i)--(iii)
are read from the ribbon graph and its boundary cycles. 
The condition~(iv) can be checked by computing winding numbers directly from the
half-edge orders of $G$, which we will recall in the following paragraphs.

\begin{definition}
\label{def:ribbon-closed-walk}
An oriented edge of $G$ is an edge with a chosen direction of
traversal.  Write $\overline e$ for the reverse traversal
of $e$.  An \emph{oriented closed edge walk} is a sequence
\[
  C=(v_0,e_1,v_1,\ldots,e_s,v_s),
  \qquad s\geq1,\qquad v_s=v_0,
\]
where $e_r$ is traversed from $v_{r-1}$ to $v_r$. The walk is
\emph{cyclically reduced} if
\[
  e_{r+1}\neq\overline{e_r}
  \qquad(1\leq r\leq s),\qquad e_{s+1}:=e_1.
\]
A \emph{simple graph cycle} is a cyclically reduced closed edge
walk for which $v_0,\ldots,v_{s-1}$ are pairwise distinct.
Thus loops and pairs of distinct parallel edges give simple
cycles of lengths one and two, respectively.  Closed walks
without a chosen starting point are considered up to cyclic
rotation, but not reversal of orientation.
\end{definition}

\begin{figure}[htbp]
\centering
\resizebox{0.94\linewidth}{!}{%
\begin{tikzpicture}[
  x=1cm,y=1cm,
  font=\scriptsize,
  panel/.style={font=\scriptsize,align=center}
]
  \foreach \shift/\tag/\entry/\exit/\value in
    {-3.05/a/2/5/2,3.05/b/5/2/-2}{
    \begin{scope}[xshift=\shift cm]
      \node[panel] at (0,2.05)
        {(\tag) $i=\entry,\ j=\exit$};

      \fill[black!4] (0,0) circle[radius=.92];
      \foreach \k/\ang in
        {1/210,2/150,3/90,4/30,5/-30,6/-90}{
        \draw[arc] (0,0)--(\ang:1.28);
        \node at (\ang:1.57) {$h_{\k}$};
      }

      \draw[arc,line width=.55pt]
        (210:.34) arc[start angle=210,end angle=270,radius=.34];
      \node at (240:.76) {$c_v$};

      \ifnum\entry=2
        \draw[auxiliary,draw=black!65,line width=.85pt,
              preaction={draw=white,line width=2.2pt,solid},
              -{Latex[length=2mm,width=1.35mm]}]
          (150:1.16)
          .. controls (-.30,.35) and (.40,.30) ..
          (-30:1.16);
      \else
        \draw[auxiliary,draw=black!65,line width=.85pt,
              preaction={draw=white,line width=2.2pt,solid},
              -{Latex[length=2mm,width=1.35mm]}]
          (-30:1.16)
          .. controls (.40,.30) and (-.30,.35) ..
          (150:1.16);
      \fi

      \node[spinevertex] at (0,0) {};
      \node at (.24,-.38) {$v$};
      \node[panel] at (0,-2.08)
        {$\varepsilon(\entry,\exit)=\value$};
    \end{scope}
  }

  \node[panel] at (0,-2.75)
    {$W_v=(h_1,h_2,h_3,h_4,h_5,h_6)$};

\end{tikzpicture}%
}
\caption{Local turns on the ribbon graph in
\eqref{eq:ribbon-local-turn}. Solid spokes represent half-edge
occurrences at $v$, and the shaded disc is a vertex neighbourhood
in the oriented thickening. The corner $c_v$ determines the
displayed linear order, read clockwise from $h_1$ to $h_6$.
The dashed arrow represents the smoothed passage from the entry
half-edge $h_i$ to the exit half-edge $h_j$. The two intermediate occurrences $h_3,h_4$ give contribution
$2$ in \textup{(a)} and $-2$ for the reversed passage
in \textup{(b)}.}
\label{fig:ribbon-local-turn}
\end{figure}
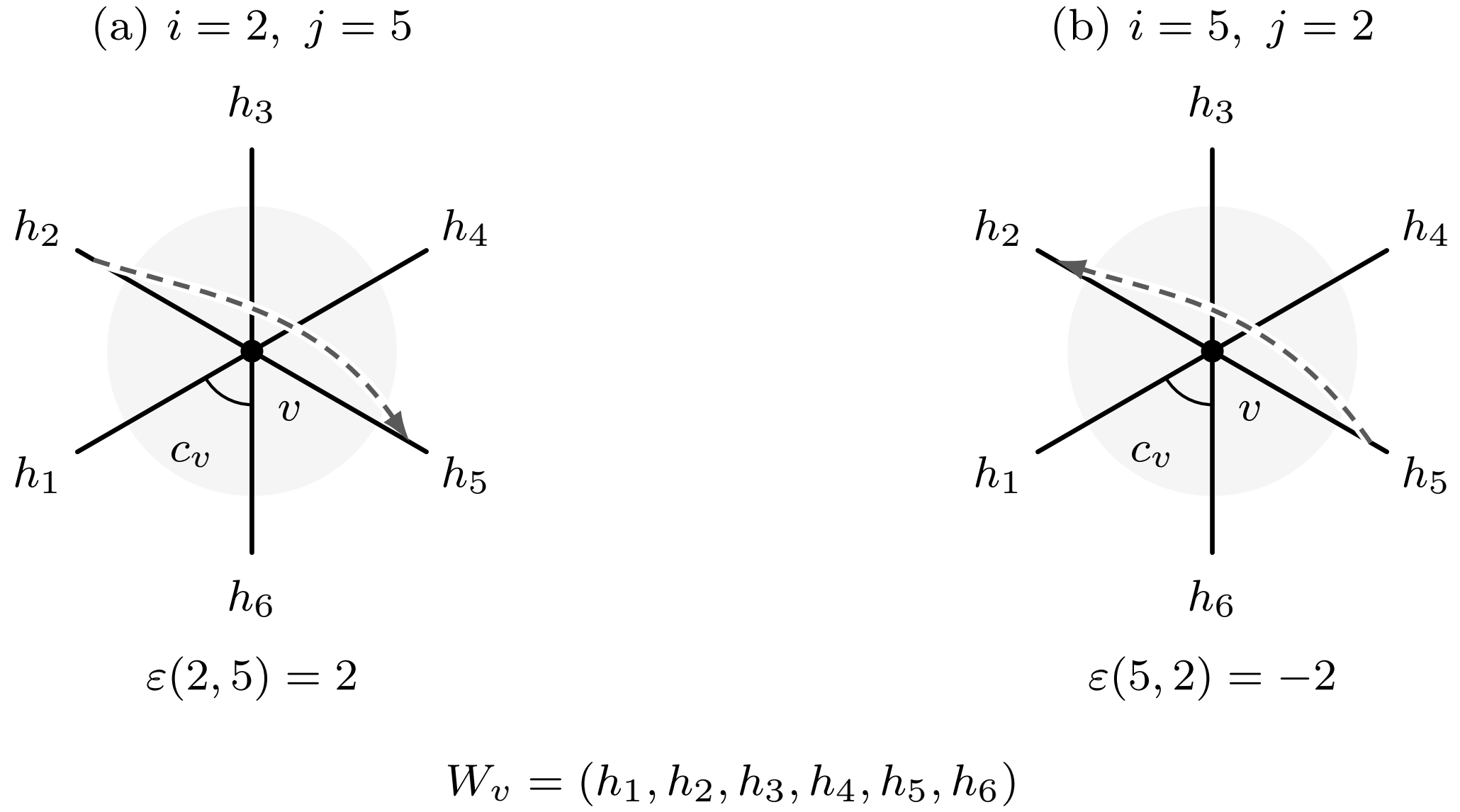

At a vertex $v$, read the half-edge occurrences clockwise,
starting immediately after the distinguished corner $c_v$.
The resulting linear list
\[
  W_v=(h_1,\ldots,h_d),\qquad d=\deg_G(v),
\]
is called the \emph{vertex word} at $v$.  Its position function
is defined by $p(h_k)=k$ for $1\leq k\leq d$; see \cref{fig:ribbon-local-turn}.
With our conventions, this is the order along the associated
forbidden thread of $A_G$.

We also write vertex words using edge labels, with each entry
representing the corresponding half-edge occurrence.  Repeated
labels are therefore treated as distinct occurrences; in
particular, the two ends of a loop occupy distinct positions.
When the vertex is understood, we simply say \emph{word}.

Let $C=(v_0,e_1,v_1,\ldots,e_s,v_s)$ be a cyclically reduced
closed edge walk.
Realise $C$ in the oriented thickening of $G$ by following
longitudinal strands in the edge bands and joining the incoming
and outgoing strands by simple connectors in the vertex discs,
without inserting curls.  Denote the resulting immersed closed
curve by $\gamma$; it may be chosen embedded when $C$ is a simple
graph cycle.  Winding numbers refer to this smoothed curve.

If the incoming and outgoing half-edges at a passage have
positions $i$ and $j$, respectively, then $i\neq j$, and the
local contribution is
\begin{equation}
\label{eq:ribbon-local-turn}
 \varepsilon(i,j)=
 \begin{cases}
  j-i-1,&i<j,\\
  1-(i-j),&i>j.
 \end{cases}
\end{equation}
Its absolute value is the number of half-edge occurrences strictly
between the two positions in the vertex word.  Adjacent positions
contribute zero, and reversing a turn changes its sign; see
\cref{fig:ribbon-local-turn}.

If the successive passages have entry and exit positions
$(i_1,j_1),\ldots,(i_s,j_s)$, then
\begin{equation}
\label{eq:ribbon-winding-sum}
  w_{\eta_G}(\gamma)=\sum_{r=1}^{s}\varepsilon(i_r,j_r).
\end{equation}
Each passage is counted separately, even when a vertex or an edge
is visited more than once, and the sum includes the closing turn.

We use the winding convention of
\cite[Definition~3.2]{AmiotPlamondonSchroll}.
The local rule follows from the construction of the line
field associated with $D_G$ and the transverse-curve rule
in \cite[Lemma~3.18]{AmiotPlamondonSchroll}.
For a related combinatorial formula, see
\cite[Remark~3.15]{LekiliPolishchuk}.
Geometrically, a vertex $v$ corresponds to the polygon $P_v$
of $\Delta_G$ containing $p_v$, and its half-edge occurrences
correspond to the dissection-side occurrences of $P_v$.
A turn from position $i$ to position $j$ is represented
by a curve segment through the interior of $P_v$ joining
the corresponding side occurrences. This explains the
local rule, but the calculation itself uses only the
words $W_v$ and the edge walk.
For the description using curves transverse to the dual dissection,
see \cite[Lemma~3.18]{AmiotPlamondonSchroll}; we use the ribbon-graph
rule throughout.

Now we are able to collect the following finite procedure computes $\GSpec(A)$.
\begin{enumerate}[label=\textup{(\arabic*)},leftmargin=2.6em]
\item
On $H=\{1,\ldots,2E\}$, enumerate fixed-point-free involutions
$\iota$, permutations $\pi$ with $m$ cycles, and one distinguished
corner in each $\pi$-cycle.  Retain the connected graphs,
equivalently those for which $\langle\iota,\pi\rangle$ acts
transitively on $H$.

\item
Compute the cycles of $\pi\iota$.  Retain graphs with $b$ boundary
cycles satisfying boundary admissibility and
\eqref{eq:R-boundary-data}, allowing all numberings of the boundary
cycles.  The genus is then $g$ by $m-E=2-2g-b$.

\item
Read the vertex words $W_v$ and their position functions $p$.
If $g>0$, choose simple closed curves
$\alpha_1,\beta_1,\ldots,\alpha_g,\beta_g$
in the oriented thickening of $G$, with $\alpha_j$ and $\beta_j$
meeting transversely once for each $j$ and all other pairs
disjoint.  Orient the curves so that their algebraic intersection
numbers satisfy
\[
  [\alpha_i]\cdot[\beta_j]=\delta_{ij},
  \qquad
  [\alpha_i]\cdot[\alpha_j]
  =
  [\beta_i]\cdot[\beta_j]=0.
\]
Here the intersection of $\alpha_j$ with $\beta_j$ is positive
when their ordered tangent directions agree with the orientation
of the surface.  Their homology classes form an integral
symplectic basis modulo the boundary classes; we call such an
oriented collection a \emph{geometric symplectic basis}.
One may orient each $\alpha_j$ arbitrarily and then orient
$\beta_j$ to make their intersection positive.

A triangulation and the usual polygonal-schema construction
provide such curves \cite{GrossTucker}.  Put these oriented curves into normal position with respect
to the edge bands and vertex discs. Represent them by
cyclically reduced edge walks whose smoothings are isotopic,
with their chosen orientations, to the original curves.
Compute their winding numbers from the entry and exit
positions using \eqref{eq:ribbon-winding-sum}.
Keep the geometric representatives: integer winding is not
obtained by linearly extending from formal homology combinations.
The boundary winding numbers are $w_i=m_i-n_i$.

\item
For $g=0$, there is no internal datum to check.  For $g=1$, compute
$d(\eta_G)$ using \eqref{eq:genus-one-gcd}.
For $g\geq2$, the winding parities of the displayed curves and
the boundary curves determine $\sigma(\eta_G)$, because winding
parity factors through mod-$2$ homology.  In this higher-genus
case, if $\sigma(\eta_G)=0$ and $w_i\equiv2\pmod4$ for every $i$,
also compute
\[
 \Arf(\eta_G)=\sum_{j=1}^g
 \left(\frac{w_{\eta_G}(\alpha_j)}2+1\right)
 \left(\frac{w_{\eta_G}(\beta_j)}2+1\right)\pmod2;
\]
see \cite[Lemma~7.5 and Theorem~7.4]{AmiotPlamondonSchroll}.
Reversing the orientation of a basis curve changes the sign
of its winding number but not its quadratic value, since in
this even case
\[
  -\frac{w_{\eta_G}(\gamma)}2+1
  \equiv
  \frac{w_{\eta_G}(\gamma)}2+1
  \pmod2.
\]
Thus the Arf formula is independent of the orientations chosen
for the individual basis curves. Retain precisely those graphs with $\tau(\eta_G)=\tau_A$.

\item
Record $\degmax(G)-1$ for each retained graph and remove duplicate
values.  By \cref{thm:intro-finite-formula}, the resulting set is
$\GSpec(A)$.
\end{enumerate}
The input datum $\tau_A$ can be computed by the same
winding-number calculations on the polygon-dual ribbon graph
of a dissection model of $A$.
Enumerating labelled graphs is sufficient: passing to
isomorphism classes only removes repetitions and does not
change the resulting set of global dimensions.

\begin{remark}
\label{rem:finite-not-efficient}
The enumeration in \cref{thm:intro-finite-formula} is exact but
need not be efficient.  Finiteness alone also follows by enumerating
gentle bound quivers with a fixed number of vertices.  The ribbon-graph
formula expresses the global-dimension objective and the complete
derived-invariant constraints in the same combinatorial data.
\end{remark}

In the following we compute the finite formula in a fibre with three simple modules, the
smallest number for which the AAG invariant can fail to distinguish gentle
derived classes \cite[Remark~4.4]{KalckDerivedUnique}.  The underlying
three-class classification is known: Kalck classified the
gentle algebras of finite global dimension with three vertices and four
arrows, and proved that two of the classes are derived unique
\cite[Corollary~4.7]{KalckDerivedUnique}. These three derived classes were also distinguished using spherelike
objects by Hochenegger, Kalck and Ploog
\cite[Section~1.1]{HKP2019}.  The two extreme representatives
also occur in \cite[Section~9]{AmiotPlamondonSchroll}, where their genus-one
gcd invariants are computed.  The following calculation recovers their
spectra directly from the ribbon-graph formula and displays the effect of
retaining the internal invariant.

\begin{proposition}
\label{prop:two-four-spectra}
Let $A$ be a connected gentle algebra of finite global dimension with
$\AG(A)=[(2,4)]$.  Its surface has genus one and one boundary component,
and $d(\eta_A)\in\{0,1,2\}$.  More precisely,
\begin{equation}
\label{eq:two-four-spectra}
 \GSpec(A)=
 \begin{cases}
  \{2\},&d(\eta_A)=0,\\
  \{2,3,4\},&d(\eta_A)=1,\\
  \{4\},&d(\eta_A)=2.
 \end{cases}
\end{equation}
All three cases occur.  In contrast, the common AAG fibre has spectrum
$\{2,3,4\}$.
\end{proposition}
\begin{proof}
The Euler identities give $m=2$, $E=3$, $g=1$, and $b=1$.
The boundary winding is $-2$, so
\[
  d(\eta)=\gcd\bigl(w_\eta(\alpha),w_\eta(\beta)\bigr)
\]
for oriented simple closed curves $\alpha,\beta$ intersecting
once; this gcd is unchanged if either curve is reversed.
By \cref{prop:complete-invariant}, it is the only internal datum
needed to distinguish derived classes in this AAG fibre.

Every graph in the finite formula has two vertices of positive
degree and three edges.  Up to interchanging the vertices, the
degree pair is therefore $(3,3)$, $(4,2)$, or $(5,1)$.
Since there is one boundary cycle, every choice of distinguished
corners is boundary-admissible.

For the case of degree pair $(3,3)$, the graph is either a theta graph or two loops joined by a bridge.
The latter has genus zero and is excluded.  Label the theta edges
$1,2,3$, all oriented from $v_1$ to $v_2$.  
After relabelling,
take $W_{v_1}=(1,2,3)$; the second word is one of its three
cyclic rotations.

Let $\alpha$ and $\beta$ be the simple closed smoothings of
$1\overline2$ and $1\overline3$, respectively.
Contracting the common tree edge $1$ leaves two interleaving
loops, so the curves can be chosen to intersect once.
We detail the case $W_{v_2}=(2,3,1)$, illustrated in
\cref{fig:two-four-theta-basis}. For the walk $1\overline2$, the incoming and outgoing half-edges
have positions $(3,1)$ at $v_2$ and $(2,1)$ at $v_1$.
For $1\overline3$, the corresponding position pairs are
$(3,2)$ at $v_2$ and $(3,1)$ at $v_1$.
The passages at $v_1$ are the closing turns; the numbers
refer to positions in the respective vertex words,
rather than to edge labels.
Thus \eqref{eq:ribbon-winding-sum} gives
\[
\begin{aligned}
  w_\eta(\alpha)
    &=\underbrace{\varepsilon(3,1)}_{v_2}
      +\underbrace{\varepsilon(2,1)}_{v_1}
      =-1+0=-1,\\
  w_\eta(\beta)
    &=\underbrace{\varepsilon(3,2)}_{v_2}
      +\underbrace{\varepsilon(3,1)}_{v_1}
      =0-1=-1.
\end{aligned}
\]
The other second words $(1,2,3)$ and $(3,1,2)$ give absolute
winding pairs $(0,0)$ and $(0,1)$, respectively.  Hence the
possible gcds are $0,1$, and the global dimension is $2$.
\begin{figure}[htbp]
\centering
\resizebox{0.98\linewidth}{!}{%
\begin{tikzpicture}[x=1cm,y=1cm,font=\small,
  theta edge/.style={
    draw=black!75,line width=.8pt,line cap=round
  },
  theta curve/.style={
    draw=black,line width=.85pt,densely dashed,
    preaction={draw=white,line width=2pt,solid},
    postaction={decorate},
    decoration={markings,
      mark=at position .20 with
        {\arrow{Latex[length=1.9mm]}}}
  },
  theta label/.style={
    font=\small,fill=white,inner sep=1pt
  },
  theta ccw/.style={
    draw=black!65,line width=.55pt,
    -{Latex[length=1.5mm]}
  }
]
  \foreach \shift/\which in {-3.65/2,3.65/3}{
\begin{scope}[xshift=\shift cm,yscale=-1]
    \coordinate (L) at (-2.4,0);
    \coordinate (R) at ( 2.4,0);

    \fill[black!5] (L) circle[radius=.60];
    \fill[black!5] (R) circle[radius=.60];

    \draw[theta edge] (L)
      .. controls (-1.2,1.8) and (1.2,1.8) .. (R);
    \draw[theta edge] (L)
      .. controls (-1.2,-2.25) and (1.1,0) .. (R);
    \draw[theta edge] (L)
      .. controls (-1.1,0) and (1.2,-2.25) .. (R);

    \ifnum\which=2
      \node[panel] at (0,-1.98)
        {(a) $\alpha\sim1\overline2$};

      \draw[theta edge,
            preaction={draw=white,line width=9pt,solid}]
        (-0.671,-0.986)
        .. controls (-0.255,-0.939) and (0.178,-0.781)
        .. (0.600,-0.604);

      \draw[theta curve] (-2.214,0.424)
        .. controls (-1.611,1.112) and (-0.806,1.460)
          .. (0,1.460)
        .. controls (0.806,1.460) and (1.611,1.112)
          .. (2.214,0.424)
        .. controls (2.333,0.289) and (2.740,0.130)
          .. (2.740,-0.080)
        .. controls (2.740,-0.480) and (2.311,-0.104)
          .. (2.134,-0.139)
        .. controls (1.553,-0.253) and (0.776,-0.701)
          .. (-0.004,-0.949)
        .. controls (-0.785,-1.196) and (-1.609,-1.227)
          .. (-2.218,-0.478)
        .. controls (-2.319,-0.354) and (-2.760,-0.210)
          .. (-2.760,-0.020)
        .. controls (-2.760,0.220) and (-2.333,0.289)
          .. (-2.214,0.424)
        -- cycle;

      \node[panel] at (0,1.83)
        {$w_\eta(\alpha)
          =\varepsilon(3,1)+\varepsilon(2,1)=-1$};

    \else
      \node[panel] at (0,-1.98)
        {(b) $\beta\sim1\overline3$};

      \draw[theta edge,
            preaction={draw=white,line width=9pt,solid}]
        (-0.600,-0.604)
        .. controls (-0.178,-0.781) and (0.255,-0.939)
        .. (0.671,-0.986);

      \draw[theta curve] (-2.214,0.424)
        .. controls (-1.611,1.112) and (-0.806,1.460)
          .. (0,1.460)
        .. controls (0.806,1.460) and (1.611,1.112)
          .. (2.214,0.424)
        .. controls (2.333,0.289) and (2.760,0.220)
          .. (2.760,-0.020)
        .. controls (2.760,-0.210) and (2.319,-0.354)
          .. (2.218,-0.478)
        .. controls (1.609,-1.227) and (0.785,-1.196)
          .. (0.004,-0.949)
        .. controls (-0.776,-0.701) and (-1.553,-0.253)
          .. (-2.134,-0.139)
        .. controls (-2.311,-0.104) and (-2.740,-0.480)
          .. (-2.740,-0.080)
        .. controls (-2.740,0.130) and (-2.333,0.289)
          .. (-2.214,0.424)
        -- cycle;

      \node[panel] at (0,1.83)
        {$w_\eta(\beta)
          =\varepsilon(3,2)+\varepsilon(3,1)=-1$};
    \fi

    \draw[arc,line width=.5pt] ($(L)+(0:.27)$)
      arc[start angle=0,end angle=56.31,radius=.27];
    \node at ($(L)+(28:.59)$) {$c_{v_1}$};

    \draw[arc,line width=.5pt] ($(R)+(123.69:.27)$)
      arc[start angle=123.69,end angle=180,radius=.27];
    \node at ($(R)+(152:.59)$) {$c_{v_2}$};

    \draw[theta ccw] ($(L)+(100:.82)$)
      arc[start angle=100,end angle=250,radius=.82];
    \draw[theta ccw] ($(R)+(-75:.82)$)
      arc[start angle=-75,end angle=75,radius=.82];

    \node[blackvertex] at (L) {};
    \node[blackvertex] at (R) {};
    \node at (-2.70,.50) {$v_1$};
    \node at ( 2.70,.50) {$v_2$};

    \node[theta label] at (-2.07,.78) {$1$};
    \node[theta label] at (-2.08,-.85) {$2$};
    \node[theta label] at (-1.52,.04) {$3$};

    \node[theta label] at ( 2.07,.78) {$1$};
    \node[theta label] at ( 1.52,.04) {$2$};
    \node[theta label] at ( 2.08,-.85) {$3$};
  \end{scope}
  }
\end{tikzpicture}%
}
\caption{The theta ribbon graph with vertex words
$W_{v_1}=(1,2,3)$ and $W_{v_2}=(2,3,1)$.
Both vertex rotations are clockwise $1\to2\to3\to1$;
$c_{v_1}$ lies between $3$ and $1$, whereas $c_{v_2}$ lies
between $1$ and $2$.  Solid lines represent $G$, and dashed arrows
show the closed smoothings of $1\overline2$ and $1\overline3$.
The vertex passages are drawn away from the distinguished corners:
$\alpha$ crosses the unused half-edge $3$ near $v_2$, and $\beta$
crosses the unused half-edge $2$ near $v_1$, without following
those edges.  The overpass is only a feature of the planar
projection, not a vertex; the band carrying the displayed curve
is drawn above the other band.}
\label{fig:two-four-theta-basis}
\end{figure}
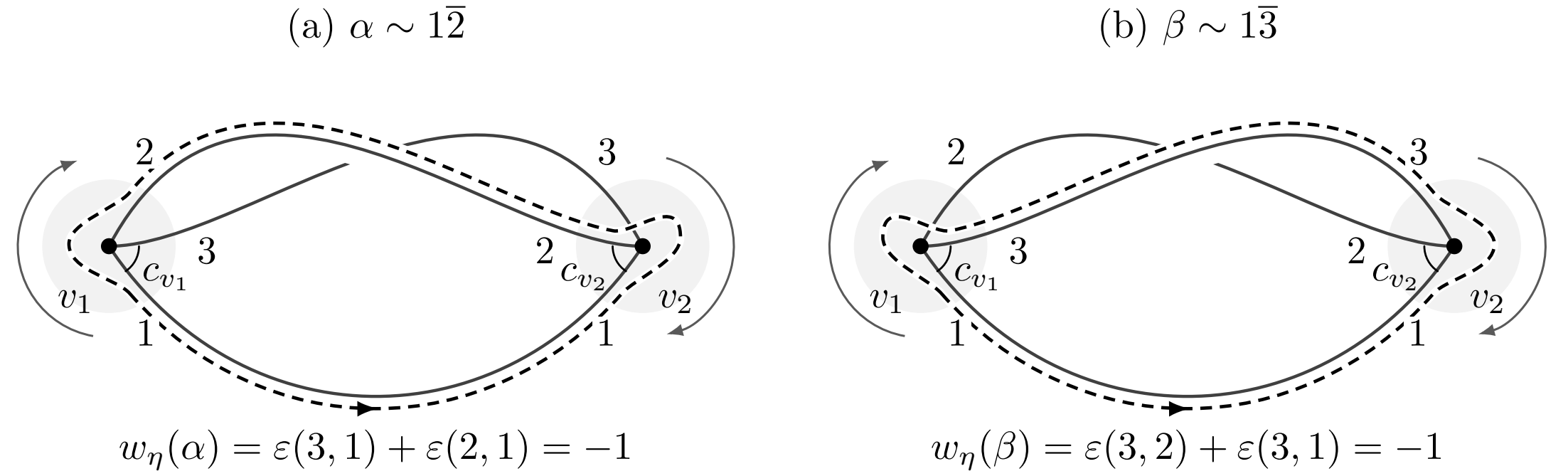

For the case of degree pair $(4,2)$, the graph has one loop $1$ at the degree-four vertex and two
joining edges $2,3$.  Genus one requires the cyclic word
$(1,2,1,3)$ at that vertex.  In every linear rotation, the two
occurrences of $1$ have positions differing by two, so the loop
has absolute winding one by \eqref{eq:ribbon-local-turn}.
Together with the smoothing of $2\overline3$, it gives a pair
intersecting once.  Thus $d(\eta)=1$ for every choice of corners,
and the global dimension is $3$.

For the case of degree pair $(5,1)$, the graph consists of two interleaving loops $1,2$ and a pendant
edge $3$.  Its vertex words are a cyclic rotation of
$(1,2,1,2,3)$ and the singleton $(3)$.
The two loops give a pair intersecting once.  Counting positions
with the pendant incidence included, the two occurrences of
each loop differ by either two or three, so its absolute winding
is either one or two.  Hence $d(\eta)\in\{1,2\}$.
Both values occur: the words $(1,2,1,2,3)$ and $(1,2,3,1,2)$
give absolute winding pairs $(1,1)$ and $(2,2)$, respectively.
The global dimension is $4$.

These cases give the complete list
\begin{equation}
\label{eq:two-four-degree-filter}
\begin{array}{c|c|c}
  \text{vertex degrees}
  &\text{possible }d(\eta)
  &\text{global dimension}\\ \hline
  (3,3)&0,1&2\\
  (4,2)&1&3\\
  (5,1)&1,2&4
\end{array}
\end{equation}
Every entry is realised by the ribbon graphs described above.
The reconstruction and the complete invariant therefore give
\eqref{eq:two-four-spectra} through
\cref{thm:intro-finite-formula}.
\end{proof}

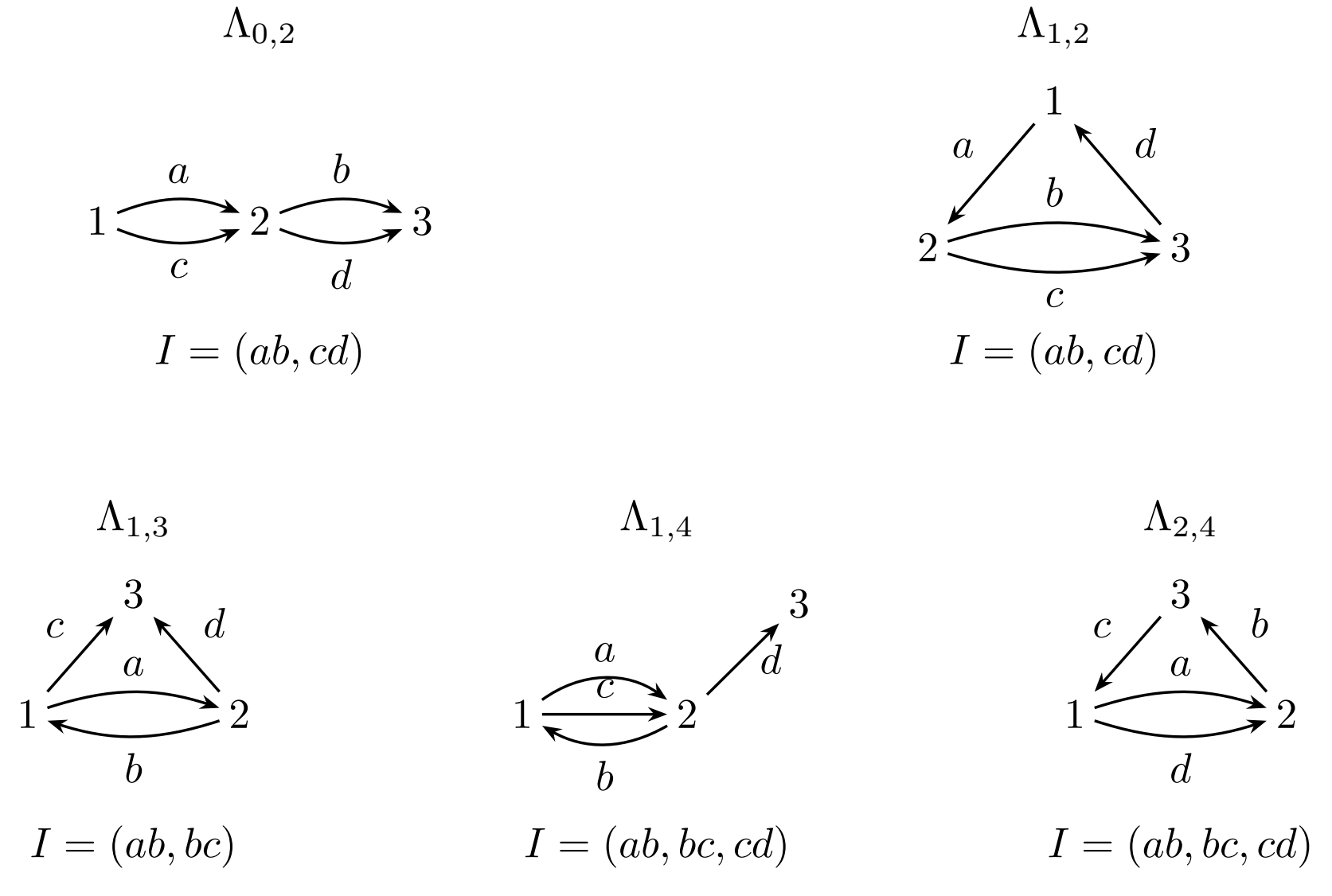
\begin{figure}[!htbp]
\centering
\begin{tikzpicture}[
 x=1cm,y=1cm,font=\small,
 qvertex/.style={inner sep=2pt,fill=white},
 qarrow/.style={-{Stealth[length=1.6mm,width=1.2mm]},line width=.65pt},
 qtitle/.style={font=\small,align=center}
]
 \begin{scope}[shift={(-3.3,2.35)}]
  \node[qtitle] at (0,1.60) {$\Lambda_{0,2}$};
  \node[qvertex] (a1) at (-1.35,0) {$1$};
  \node[qvertex] (a2) at (0,0) {$2$};
  \node[qvertex] (a3) at (1.35,0) {$3$};
  \draw[qarrow] (a1) to[bend left=22] node[above] {$a$} (a2);
  \draw[qarrow] (a1) to[bend right=22] node[below] {$c$} (a2);
  \draw[qarrow] (a2) to[bend left=22] node[above] {$b$} (a3);
  \draw[qarrow] (a2) to[bend right=22] node[below] {$d$} (a3);
  \node[qtitle] at (0,-1.10) {$I=(ab,cd)$};
 \end{scope}
 \begin{scope}[shift={(3.3,2.35)}]
  \node[qtitle] at (0,1.60) {$\Lambda_{1,2}$};
  \node[qvertex] (b1) at (0,1.00) {$1$};
  \node[qvertex] (b2) at (-1.05,-.22) {$2$};
  \node[qvertex] (b3) at (1.05,-.22) {$3$};
  \draw[qarrow] (b1) -- node[above left] {$a$} (b2);
  \draw[qarrow] (b2) to[bend left=17] node[above] {$b$} (b3);
  \draw[qarrow] (b2) to[bend right=17] node[below] {$c$} (b3);
  \draw[qarrow] (b3) -- node[above right] {$d$} (b1);
  \node[qtitle] at (0,-1.10) {$I=(ab,cd)$};
 \end{scope}
 \begin{scope}[shift={(-4.35,-1.75)}]
  \node[qtitle] at (0,1.60) {$\Lambda_{1,3}$};
  \node[qvertex] (c1) at (-.88,0) {$1$};
  \node[qvertex] (c2) at (.88,0) {$2$};
  \node[qvertex] (c3) at (0,1.00) {$3$};
  \draw[qarrow] (c1) to[bend left=18] node[above] {$a$} (c2);
  \draw[qarrow] (c2) to[bend left=18] node[below] {$b$} (c1);
  \draw[qarrow] (c1) -- node[above left] {$c$} (c3);
  \draw[qarrow] (c2) -- node[above right] {$d$} (c3);
  \node[qtitle] at (0,-1.10) {$I=(ab,bc)$};
 \end{scope}
 \begin{scope}[shift={(0,-1.75)}]
  \node[qtitle] at (0,1.60) {$\Lambda_{1,4}$};
  \node[qvertex] (d1) at (-1.12,0) {$1$};
  \node[qvertex] (d2) at (.25,0) {$2$};
  \node[qvertex] (d3) at (1.18,.92) {$3$};
  \draw[qarrow] (d1) to[bend left=36] node[above] {$a$} (d2);
  \draw[qarrow] (d1) -- node[above] {$c$} (d2);
  \draw[qarrow] (d2) to[bend left=30] node[below] {$b$} (d1);
  \draw[qarrow] (d2) -- node[right] {$d$} (d3);
  \node[qtitle] at (0,-1.10) {$I=(ab,bc,cd)$};
 \end{scope}
 \begin{scope}[shift={(4.35,-1.75)}]
  \node[qtitle] at (0,1.60) {$\Lambda_{2,4}$};
  \node[qvertex] (e1) at (-.88,0) {$1$};
  \node[qvertex] (e2) at (.88,0) {$2$};
  \node[qvertex] (e3) at (0,1.00) {$3$};
  \draw[qarrow] (e1) to[bend left=18] node[above] {$a$} (e2);
  \draw[qarrow] (e1) to[bend right=18] node[below] {$d$} (e2);
  \draw[qarrow] (e2) -- node[above right] {$b$} (e3);
  \draw[qarrow] (e3) -- node[above left] {$c$} (e1);
  \node[qtitle] at (0,-1.10) {$I=(ab,bc,cd)$};
 \end{scope}
\end{tikzpicture}
\caption{Complete quivers and generating zero relations for the five
representatives.  All arrows are displayed, and path products are read
from left to right.  Each algebra has AAG invariant $[(2,4)]$.}
\label{fig:two-four-quivers}
\end{figure}

For explicit algebra representatives, \cref{fig:two-four-quivers}
gives the five algebras
\[\Lambda_{0,2},\Lambda_{1,2},\Lambda_{1,3},
\Lambda_{1,4},\Lambda_{2,4}\] 
realising the pairs in
\eqref{eq:two-four-degree-filter}.
The first subscript records the gcd and the second the global
dimension.  
Their maximal forbidden threads and vertex words are
listed below.  Here $e_3$ denotes the trivial forbidden thread at
quiver vertex $3$, and repeated labels in a word denote distinct
occurrences.
\begin{center}
\small
\renewcommand{\arraystretch}{1.15}
\begin{tabular}{@{}cll@{}}
  Algebra & Forbidden threads & Vertex words\\ \hline
  $\Lambda_{0,2}$ & $ab,\,cd$
    & $(1,2,3);\ (1,2,3)$\\
  $\Lambda_{1,2}$ & $ab,\,cd$
    & $(1,2,3);\ (2,3,1)$\\
  $\Lambda_{1,3}$ & $abc,\,d$
    & $(1,2,1,3);\ (2,3)$\\
  $\Lambda_{1,4}$ & $abcd,\,e_3$
    & $(1,2,1,2,3);\ (3)$\\
  $\Lambda_{2,4}$ & $abcd,\,e_3$
    & $(1,2,3,1,2);\ (3)$
\end{tabular}
\end{center}
Pairing equal labels recovers the one-boundary ribbon graphs
used in the proof, with gcds $0,1,1,1,2$, respectively.
The longest forbidden threads have lengths $2,2,3,4,4$,
in agreement with the global dimensions in
\eqref{eq:two-four-degree-filter}.

The singleton spectra in \eqref{eq:two-four-spectra} also follow from the
derived uniqueness of Kalck's two representatives, which agree, up to
relabelling, with $\Lambda_{0,2}$ and $\Lambda_{2,4}$ above
\cite[Corollary~4.7]{KalckDerivedUnique}.  The ribbon-graph calculation makes the role of the internal gcd explicit:
it selects different degree pairs within one AAG fibre.  In particular,
$\Lambda_{0,2}$ and $\Lambda_{2,4}$ are counterexamples to the equality
$\GSpec(A)=\AGSpec(A)$.  Their different gcd invariants show that they
are not derived equivalent, and all three spectra in
\eqref{eq:two-four-spectra} are intervals. The next section determines exactly when internal gaps occur
and completes the classification of the spectra.

\section{Classification of noninterval spectra}
\label{sec:spectrum-classification}

Proposition~\ref{prop:two-four-spectra} shows that fixing the
internal line-field datum can shrink the AAG interval while
leaving an interval.  We now determine exactly when internal
gaps occur.  

Throughout this section, write
\[
 \AG(A)=[(m_i,n_i)]_{i=1}^b,\qquad
 m=\sum_i m_i,\qquad w_i=m_i-n_i,
 \qquad L=L_A,\quad U=U_A.
\]
A \emph{weight profile} consists of nonnegative integers
$(\lambda_v)_{v\in V}$, partitioned into sets $V_i$ with
$|V_i|=m_i$ and $\sum_{v\in V_i}\lambda_v=n_i$.
For a ribbon realisation, $\lambda_v=\deg(v)-1$ and the global dimension
is $\max_v\lambda_v$.  By \cref{thm:all-genus-realization}, every profile
with the Euler sum is realisable.  The remaining task is to control
$\tau(\eta)$ while realising the profile.  We use the vertex
words and half-edge positions fixed
in \cref{sec:finite-formula}, and compute winding numbers by
\eqref{eq:ribbon-winding-sum} throughout.

\subsection{Local constructions and reduced profiles}
\label{subsec:classification-constructions}

The following operations modify the vertex words to control
winding numbers.  In a decomposition $W_v=W_1W_2$, the factors
$W_1,W_2$ are consecutive blocks.
Juxtaposition denotes concatenation, and $|W_i|$ denotes the
number of half-edge occurrences in $W_i$.
We use edge labels in the displayed words; each new edge label
occurs exactly twice in the resulting collection, once for
each of its half-edges.

\begin{lemma}[Splitting a vertex]
\label{lem:classification-splitting}
Let $W_v=W_1W_2$ be a decomposition of the vertex word at $v$
into two consecutive blocks, either of which may be empty.
Split $v$ into two vertices $u,z$ joined by a new edge $e$,
with vertex words
\begin{equation}
\label{eq:classification-splitting}
  W_u=(W_1,e),\qquad W_z=(e,W_2).
\end{equation}
Here the two occurrences of $e$ represent its half-edges at
$u$ and $z$, and all old half-edge pairings are retained.
The two new distinguished corners lie on the boundary
containing the old one.  On that boundary, $(m_i,n_i)$
increases by $(1,1)$.  The genus, boundary windings, and
internal line-field invariant are unchanged. The operation is illustrated in
\cref{fig:classification-splitting}.
\end{lemma}
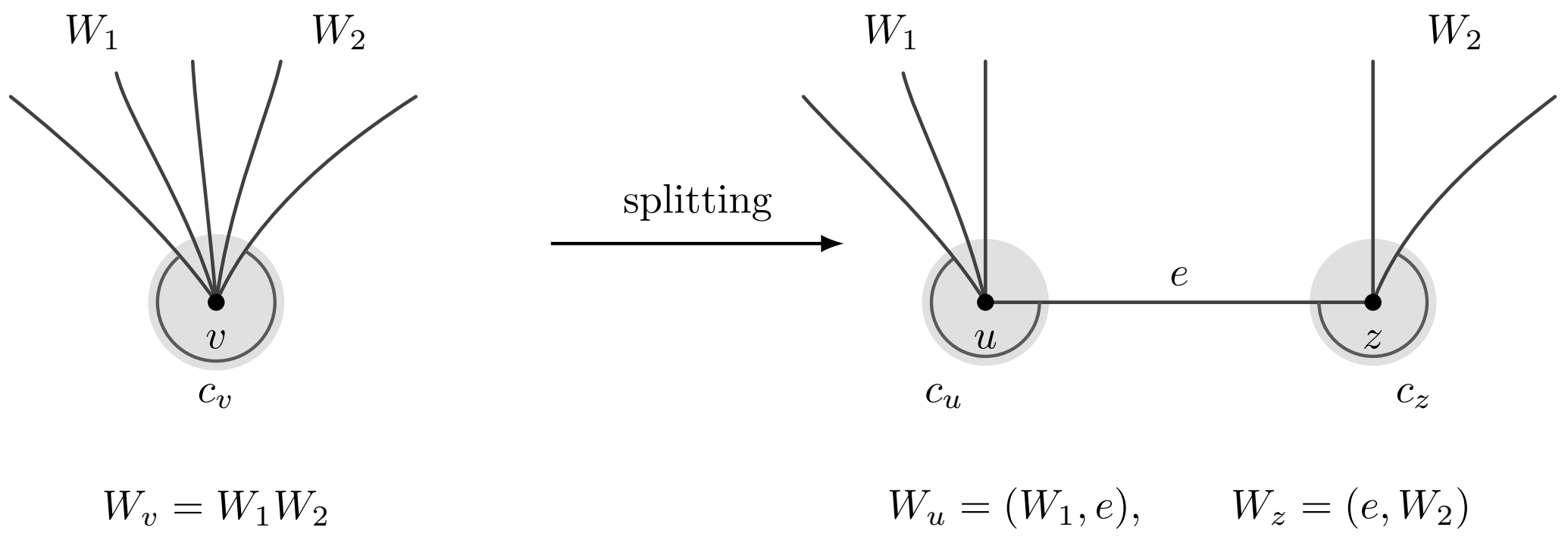
\begin{figure}[htbp]
\centering
\resizebox{0.92\linewidth}{!}{%
\begin{tikzpicture}[x=1cm,y=1cm,font=\small,
  split edge/.style={
    draw=black!75,line width=.85pt,line cap=round
  },
  split corner/.style={
    draw=black!65,line width=.8pt
  },
  split label/.style={
    font=\small,fill=white,inner sep=1.5pt
  },
  split arrow/.style={
    -{Latex[length=2mm]},line width=.8pt
  },
  vertex shade/.style={
    fill=black!12
  }
]

\begin{scope}[xshift=-4.1cm]
  \coordinate (v) at (0,0);

  \fill[vertex shade] (v) circle[radius=.58];

  \draw[split edge] (v) .. controls (-.45,.65) and (-1.25,1.35)
    .. (-1.75,1.75);
  \draw[split edge] (v) .. controls (-.25,.85) and (-.75,1.55)
    .. (-.85,1.95);
  \draw[split edge] (v) .. controls (-.08,.90) and (-.18,1.65)
    .. (-.20,2.05);

  \draw[split edge] (v) .. controls (.12,.90) and (.45,1.60)
    .. (.55,2.05);
  \draw[split edge] (v) .. controls (.35,.80) and (1.15,1.40)
    .. (1.70,1.75);

  \node[blackvertex] at (v) {};
  \node[below=3pt] at (v) {$v$};

  \node[split label] at (-1.05,2.30) {$W_1$};
  \node[split label] at (1.05,2.30) {$W_2$};

  \draw[split corner]
    ($(v)+(60:.50)$)
    arc[start angle=60,end angle=-230,radius=.50];
  \node[split label] at (0,-.8) {$c_v$};

  \node[panel] at (0,-1.75) {$W_v=W_1W_2$};
\end{scope}

\draw[split arrow] (-1.25,.50)--(1.25,.50);
\node[panel] at (0,.85) {splitting};

\begin{scope}[xshift=4.1cm]
  \coordinate (u) at (-1.65,0);
  \coordinate (z) at (1.65,0);

  \fill[vertex shade] (u) circle[radius=.54];
  \fill[vertex shade] (z) circle[radius=.54];

  \draw[split edge] (u) .. controls (-2.05,.65) and (-2.80,1.30)
    .. (-3.20,1.75);
  \draw[split edge] (u) .. controls (-1.85,.85) and (-2.25,1.55)
    .. (-2.35,1.95);
  \draw[split edge] (u) .. controls (-1.65,.90) and (-1.65,1.65)
    .. (-1.65,2.05);

  \draw[split edge] (z) .. controls (1.65,.90) and (1.65,1.65)
    .. (1.65,2.05);
  \draw[split edge] (z) .. controls (1.95,.80) and (2.75,1.40)
    .. (3.20,1.75);

  \draw[split edge] (u)--(z);
  \node[split label] at (0,.22) {$e$};

  \node[blackvertex] at (u) {};
  \node[blackvertex] at (z) {};
  \node[below=3pt] at (u) {$u$};
  \node[below=3pt] at (z) {$z$};

  \node[split label] at (-2.45,2.30) {$W_1$};
  \node[split label] at (2.35,2.30) {$W_2$};

  \draw[split corner]
    ($(u)+(0:.46)$)
    arc[start angle=0,end angle=-233,radius=.46];
  \node[split label] at (-2.00,-0.80) {$c_u$};

  \draw[split corner]
    ($(z)+(62:.46)$)
    arc[start angle=62,end angle=-180,radius=.46];
  \node[split label] at (2.00,-.80) {$c_z$};

  \node[panel] at (0,-1.75)
    {$W_u=(W_1,e),\qquad W_z=(e,W_2)$};
\end{scope}

\end{tikzpicture}%
}
\caption{Splitting a vertex.  The vertex word
$W_v=W_1W_2$ is replaced by
$W_u=(W_1,e)$ and $W_z=(e,W_2)$, where $e$ is a new edge.
The clockwise half-edge orders are preserved within each block.
The distinguished corners $c_u$ and $c_z$ lie on the same
boundary cycle as the original corner $c_v$.}
\label{fig:classification-splitting}
\end{figure}
\begin{proof}
Let $G'$ be the graph after splitting, and let $S,S'$ be the
oriented thickenings of $G,G'$.  Put
$a=|W_1|$ and $d=\deg_G(v)$.

We first examine the topology and boundary data.
In $S'$, the two vertex discs at $u,z$, together with the band
of $e$, form a disc.  Its attachment intervals for the old
half-edges have the same clockwise cyclic order as those of
the original vertex disc at $v$.  This gives an
orientation-preserving identification of the underlying
unmarked surfaces $S$ and $S'$, fixed outside the modified
neighbourhood.  In particular, the genus and the number of
boundary components are unchanged.

Let $F_i$ be the boundary containing $c_v$.  The boundary
interval corresponding to $c_v$ is replaced by the two
intervals corresponding to $c_z,c_u$, connected by one side
of the new edge band.  Thus both new distinguished corners
lie on the boundary corresponding to $F_i$.
If a block is empty, the corresponding new vertex is a leaf,
and the same boundary interval simply runs around its vertex
disc.  All other distinguished corners retain their boundary
assignments.

The new degrees are $a+1$ and $d-a+1$, so the new weights
are $a$ and $d-a$, replacing the old weight $d-1$.
Consequently,
\[
  (m_i',n_i')=(m_i+1,n_i+1),
  \qquad
  (m_k',n_k')=(m_k,n_k)\quad(k\neq i).
\]
It follows that $w_k'=m_k'-n_k'=w_k$ for every boundary
component.

Next we compare winding numbers.
Compare a cyclically reduced closed edge walk in $G$ with its
corresponding walk in $G'$. Denote by $\gamma$ and $\gamma'$ the curves in $S$ and $S'$ corresponding to these walks in $G$ and $G'$, respectively. A passage within one block stays
at the corresponding new vertex, whereas a passage between
the two blocks traverses $e$.  The connectors can be chosen
to correspond by isotopy inside the identified discs, without
inserting curls.

For a turn within either block, position differences are
unchanged, so the local contribution is unchanged.
For a turn from position $i$ in $W_1$ to position $j$ in
$W_2$, the old contribution is $\varepsilon(i,a+j)$.
The corresponding turns at $u$ and $z$ contribute
\[
\begin{aligned}
  \varepsilon(i,a+1)+\varepsilon(1,j+1)
    &=(a-i)+(j-1)\\
    &=a+j-i-1\\
    &=\varepsilon(i,a+j).
\end{aligned}
\]
The reverse passage gives the negative identity.
If a block is empty, the new edge is pendant and is not
traversed by a cyclically reduced closed walk; the positions
of the old half-edges are unchanged or shifted uniformly.
Thus \eqref{eq:ribbon-winding-sum} gives
\[
  w_{\eta_{G'}}(\gamma')=w_{\eta_G}(\gamma)
\]
for corresponding smoothed closed walks.

To consider the internal line-field invariant, we choose a geometric symplectic basis on $S$ and the corresponding
oriented curves on $S'$.  Their intersection numbers agree,
and the preceding calculation shows that their winding numbers
agree.  Together with the equality of the boundary windings,
this preserves the genus-one gcd in \eqref{eq:genus-one-gcd}.
For $g\geq2$, winding parity agrees on the basis curves and the
boundary curves, which generate homology, so $\sigma$ is
unchanged.  The conditions for the Arf invariant to be defined
are therefore the same on both surfaces.  When it is defined,
the values
\[
  \frac{w_{\eta_{G'}}(\gamma')}{2}+1
  \equiv
  \frac{w_{\eta_G}(\gamma)}{2}+1
  \pmod2
\]
agree on the corresponding symplectic basis curves, so the
Arf invariants agree as well.  In genus zero there is no
internal datum.  Hence
$\tau(\eta_{G'})=\tau(\eta_G)$ in every case.
\end{proof}
On the surface, the operation of splitting a vertex adds one marked point of each
colour to the boundary side of the polygon corresponding to
$v$, together with a new $\circ$-arc dividing the polygon
into two, each containing one $\bullet$-point.

We also use \emph{suppression}, the inverse of edge subdivision:
a bivalent vertex with word $(e,f)$ is removed by joining
its two incident edge segments. This preserves winding numbers,
since $\varepsilon(1,2)=\varepsilon(2,1)=0$.
\begin{lemma}[Adding a handle]
\label{lem:classification-handle}
Let $W_v=W_1PW_2$ be a vertex word, where $P$ is a consecutive
block of $r$ half-edge occurrences.  Replace $W_v$ by
\begin{equation}
\label{eq:classification-handle}
  (W_1,x,y,P,x,y,W_2),
\end{equation}
where $x,y$ are new loop edges, see \cref{fig:classification-handle}. This increases the genus by one and the
selected weight by four.  Boundary components and their distinguished
corner assignments are preserved; only the corresponding $n_i$
increases, by four.  Moreover:
\begin{enumerate}[label=\textup{(\roman*)}]
\item if $r$ is even, the new line field has a closed curve of odd winding;
\item if $r$ is odd, an even line field remains even;
\item if the Arf invariant is defined before the insertion
and $r$ is odd, it remains defined afterwards; its value
is unchanged for $r\equiv1\pmod4$ and changes by one
for $r\equiv3\pmod4$;
\item on a genus-zero starting surface, the resulting genus-one invariant is
\begin{equation}
\label{eq:classification-handle-gcd}
 \gcd(r+1,w_1'+2,\ldots,w_b'+2),
\end{equation}
where the $w_i'$ are the resulting boundary windings.
\end{enumerate}
\end{lemma}
\begin{figure}[htbp]
\centering
\resizebox{0.94\linewidth}{!}{%
\begin{tikzpicture}[x=1cm,y=1cm,font=\small,
  handle edge/.style={
    draw=black!75,line width=.85pt,line cap=round
  },
  handle corner/.style={
    draw=black!65,line width=.8pt
  },
  handle label/.style={
    font=\small,fill=white,inner sep=1.5pt
  },
  handle arrow/.style={
    -{Latex[length=2mm]},line width=.8pt
  },
  vertex shade/.style={
    fill=black!12
  }
]

\begin{scope}[xshift=-4.6cm]
  \coordinate (vL) at (0,0);

  \fill[vertex shade] (vL) circle[radius=.60];
  \draw[black!35,line width=.6pt] (vL) circle[radius=.60];

  \draw[handle edge] (vL) .. controls (-.52,.62) and (-1.28,1.30)
    .. (-1.92,1.80);
  \draw[handle edge] (vL) .. controls (-.22,.86) and (-.62,1.56)
    .. (-1.02,2.00);

  \draw[handle edge] (vL) .. controls (-.05,.94) and (-.06,1.65)
    .. (-.10,2.12);
  \draw[handle edge] (vL) .. controls (.12,.92) and (.22,1.60)
    .. (.28,2.04);

  \draw[handle edge] (vL) .. controls (.32,.82) and (.98,1.40)
    .. (1.10,1.92);
  \draw[handle edge] (vL) .. controls (.48,.74) and (1.15,1.26)
    .. (1.92,1.76);

  \node[blackvertex] at (vL) {};
  \node[below=3pt] at (vL) {$v$};

  \node[handle label] at (-1.18,2.38) {$W_1$};
  \node[handle label] at (.1,2.48) {$P$};
  \node[handle label] at (1.18,2.38) {$W_2$};

  \draw[handle corner]
    ($(vL)+(60:.52)$)
    arc[start angle=60,end angle=-230,radius=.52];
  \node[handle label] at (0,-.84) {$c_v$};

  \node[panel] at (0,-1.84) {$W_v=W_1PW_2$};
\end{scope}

\draw[handle arrow] (-1.48,.50)--(1.48,.50);
\node[panel] at (0,.86) {adding a handle};

\begin{scope}[xshift=4.6cm]
  \coordinate (vR) at (0,0);

  \fill[vertex shade] (vR) circle[radius=.60];
  \draw[black!35,line width=.6pt] (vR) circle[radius=.60];


  \draw[handle edge] (vR) .. controls (-.52,.62) and (-1.28,1.30)
    .. (-1.92,1.80);
  \draw[handle edge] (vR) .. controls (-.24,.86) and (-.64,1.56)
    .. (-1.02,2.00);

  \coordinate (xL)    at (-.46,2.08);
  \coordinate (yL)    at (-.18,2.18);
  \coordinate (PmidA) at (.04,2.18);
  \coordinate (PmidB) at (.24,2.16);
  \coordinate (xR)    at (.50,2.06);
  \coordinate (yR)    at (.82,1.94);
  \coordinate (WtwoA) at (1.36,1.92);
  \coordinate (WtwoB) at (1.92,1.76);

  \draw[handle edge] (vR) .. controls (-.10,.92) and (-.24,1.62) .. (xL);
  \draw[handle edge] (vR) .. controls (-.02,.96) and (-.04,1.72) .. (yL);
  \draw[handle edge] (vR) .. controls (.06,.98) and (.06,1.74) .. (PmidA);
  \draw[handle edge] (vR) .. controls (.14,.98) and (.18,1.72) .. (PmidB);
  \draw[handle edge] (vR) .. controls (.24,.94) and (.34,1.64) .. (xR);
  \draw[handle edge] (vR) .. controls (.36,.88) and (.54,1.48) .. (yR);
  \draw[handle edge] (vR) .. controls (.52,.82) and (.94,1.36) .. (WtwoA);
  \draw[handle edge] (vR) .. controls (.66,.74) and (1.22,1.26) .. (WtwoB);

  \draw[handle edge]
    (xL)
      .. controls (-.78,3.05) and (-.08,3.42)
      .. (.40,3.22)
      .. controls (.82,3.00) and (.86,2.48)
      .. (xR);
  \node[handle label] at (.08,3.45) {$x$};

  \draw[handle edge]
    (yL)
      .. controls (.06,2.78) and (.54,2.98)
      .. (1.06,2.82)
      .. controls (1.22,2.68) and (1.12,2.28)
      .. (yR);
  \node[handle label] at (1.18,3) {$y$};

  \node[blackvertex] at (vR) {};
  \node[below=3pt] at (vR) {$v$};

  \node[handle label] at (-1.18,2.38) {$W_1$};
  \node[handle label] at (.37,2.52) {$P$};
  \node[handle label] at (1.58,2.42) {$W_2$};

  \draw[handle corner]
    ($(vR)+(60:.52)$)
    arc[start angle=60,end angle=-230,radius=.52];
  \node[handle label] at (0,-.84) {$c_v$};

  \node[panel] at (0,-1.84)
    {$W_v=(W_1,x,y,P,x,y,W_2)$};
\end{scope}

\end{tikzpicture}%
}
\caption{Adding a handle at a vertex.  The vertex word
$W_v=W_1PW_2$ is replaced by
$W_v=(W_1,x,y,P,x,y,W_2)$, where $x$ and $y$ are new loop
edges.  The two loops are interleaved in the clockwise order,
producing one additional handle.  The distinguished corner
remains between $W_2$ and $W_1$.}
\label{fig:classification-handle}
\end{figure}
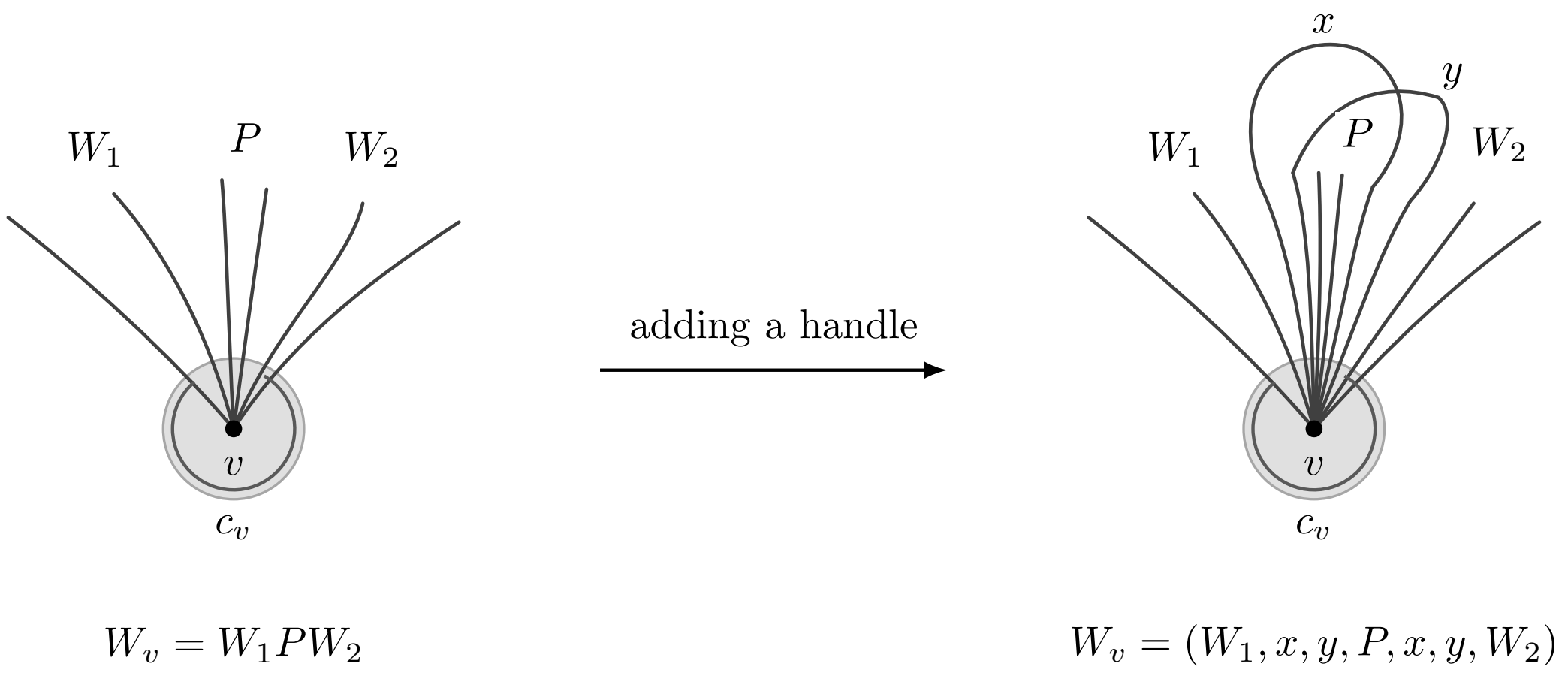
\begin{proof}
Let $G'$ be the graph after the insertion, and let $S,S'$
be the oriented thickenings of $G,G'$, respectively.
Distinguish the two occurrences of
each new loop by $x_1,x_2$ and $y_1,y_2$, so that the new
vertex word is
\[
  (W_1,x_1,y_1,P,x_2,y_2,W_2).
\]

We first consider the change of topology and boundary data.
The boundary cycles are the cycles of $\pi\iota$ and
$\pi'\iota'$, respectively. Suppose first that the two
insertion sites, immediately before and after $P$, are
distinct. Write
\[
  a\longmapsto b,
  \qquad
  c\longmapsto d
\]
for the old boundary transitions across these two sites,
respectively. Under $\pi'\iota'$, they are replaced by
\[
  a\longmapsto x_1\longmapsto y_2\longmapsto b,
  \qquad
  c\longmapsto x_2\longmapsto y_1\longmapsto d.
\]
For example, $\iota'(x_1)=x_2$ and $\pi'(x_2)=y_2$,
which explains the step $x_1\longmapsto y_2$.

If the two insertion sites coincide, the two detours
concatenate into a single detour containing all four new
half-edges. In either case, every old boundary cycle is
modified only by inserting these detours, and no new
boundary cycle is created. Thus the boundary components
and their distinguished corner assignments are preserved.

The insertion adds two edges without changing the vertex
set.  Hence
\[
  \chi(S')=\chi(S)-2.
\]
Since the number of boundary components is unchanged,
the genus increases by one.

At the modified vertex $v$, the degree and weight increase
by four.  Thus, on the boundary containing its
distinguished corner,
\[
  (m_i',n_i')=(m_i,n_i+4),
  \qquad w_i'=w_i-4,
\]
while all other AAG pairs and boundary windings remain
unchanged.  In particular, every boundary winding is
unchanged modulo four.

Now we consider winding numbers and parity.
The half-edges of each new loop occur at positions
differing by $r+2$.  Therefore, by
\eqref{eq:ribbon-local-turn}, the simple closed
smoothings $\gamma_x,\gamma_y$ of the loops satisfy
\[
  |w_{\eta_{G'}}(\gamma_x)|
  =|w_{\eta_{G'}}(\gamma_y)|=r+1.
\]
Their half-edges alternate in the vertex word, so
the two curves intersect once and form a handle pair.

We next compare the winding numbers of the old graph
cycles with those of their natural continuations in $G'$.
Every old half-edge retains its relative order, but its
position is shifted by
\[
  s(h)=
  \begin{cases}
    0,&h\in W_1,\\
    2,&h\in P,\\
    4,&h\in W_2.
  \end{cases}
\]
For a turn with old entry and exit positions $i,j$,
the new positions are $i+s_i,j+s_j$, where
$s_i,s_j\in\{0,2,4\}$ are the corresponding shifts.
Since their order is preserved,
\[
  \varepsilon(i+s_i,j+s_j)
  =\varepsilon(i,j)+s_j-s_i.
\]
Thus every old local contribution changes by an even
integer.  Summing over the turns shows that the
winding parity of every old graph cycle is preserved.

If $r$ is even, the new loop curves have odd winding,
so the new line field has $\sigma(\eta_{G'})=1$.
If $r$ is odd, both new loops have even winding.
Together with the unchanged parities on the old cycles
and boundary curves, this shows that an even line
field remains even.  This proves \textup{(i)} and
\textup{(ii)}.

For the Arf invariant, we suppose that the Arf invariant is defined and $r$ is
odd.  Let $q,q'$ be the quadratic forms on the capped
homology before and after the insertion.

Choose a geometric symplectic basis of the old surface,
represented by smoothed graph cycles $z_1,\ldots,z_{2g}$,
and regard these cycles as curves in the new surface.
Let $f_j\in\mathbb F_2$ be the parity of the passages
of $z_j$ between $P$ and its complement in the old
vertex word.

The position-shift calculation above gives
\[
  q'(z_j)=q(z_j)+f_j.
\]
Indeed, modulo four, an old turn changes its winding
by two exactly when one of its half-edges belongs to
$P$ and the other does not.

The alternating half-edges of the new loops give
\[
  z_j\cdot\gamma_x
  =z_j\cdot\gamma_y=f_j,
  \qquad
  \gamma_x\cdot\gamma_y=1
  \quad\text{in }\mathbb F_2.
\]
Moreover, since $r$ is odd,
\[
  q'(\gamma_x)=q'(\gamma_y)
  =\frac{r+3}{2}\pmod2.
\]

Put
\[
  z_j'=z_j+f_j(\gamma_x+\gamma_y).
\]
As
\[
  q'(\gamma_x+\gamma_y)=1,
\]
the classes $z_j'$ are orthogonal to
$\gamma_x,\gamma_y$, retain the mutual intersection
pairings of the old basis, and satisfy
\[
  q'(z_j')=q(z_j).
\]
Consequently, the new quadratic form is the orthogonal
sum of the old form and the plane spanned by
$\gamma_x,\gamma_y$.

The Arf contribution of this plane is
\[
  q'(\gamma_x)q'(\gamma_y)
  =\frac{r+3}{2}\pmod2.
\]
It is zero for $r\equiv1\pmod4$ and one for
$r\equiv3\pmod4$.
This proves \textup{(iii)}.

Finally, suppose that the starting surface has genus zero.
The insertion creates a single handle, represented by the
curves $\gamma_x,\gamma_y$, which form a geometric
symplectic basis modulo the boundary classes.
Their winding numbers have absolute value $r+1$.
Therefore \eqref{eq:genus-one-gcd} gives
\[
  d(\eta_{G'})
  =\gcd(r+1,w_1'+2,\ldots,w_b'+2),
\]
which proves \textup{(iv)}.
\end{proof}

\begin{lemma}[Adding a boundary]
\label{lem:classification-boundary}
Let $G$ be a boundary-admissible ribbon graph.
Subdivide an edge by a new vertex $v$, replacing it with
two edges $e,f$, and attach a loop $x$ at $v$.
Choose the distinguished corner at $v$ so that its vertex
word is
\[
  W_v=(x,e,f,x),
\]
and retain all old distinguished corners, see
\cref{fig:classification-boundary}.
The resulting graph $G'$ is boundary-admissible and has
one additional boundary component with AAG pair $(1,3)$.
All old AAG pairs and the genus are preserved.
The new boundary winding is $-2$, while the old boundary
windings and the internal line-field invariant are unchanged.
\end{lemma}
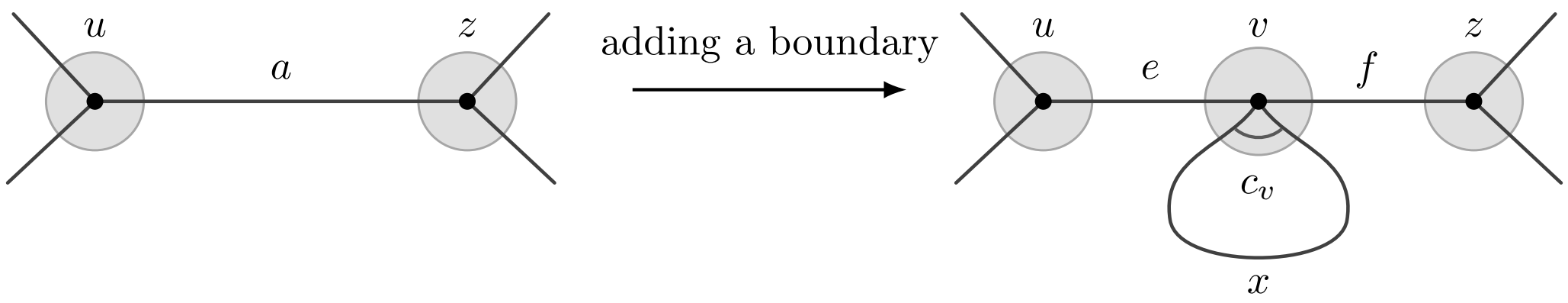
\begin{figure}[htbp]
\centering
\resizebox{0.92\linewidth}{!}{%
\begin{tikzpicture}[x=1cm,y=1cm,font=\small,
  bd edge/.style={
    draw=black!75,line width=.85pt,line cap=round
  },
  bd shade/.style={
    fill=black!12,draw=black!35,line width=.55pt
  },
  bd corner/.style={
    draw=black!65,line width=.8pt
  },
  bd label/.style={
    font=\small,fill=white,inner sep=1.5pt
  },
  bd arrow/.style={
    -{Latex[length=2mm]},line width=.8pt
  }
]

\begin{scope}[xshift=-4.2cm]
  \coordinate (u) at (-1.60,0);
  \coordinate (z) at ( 1.60,0);

  \filldraw[bd shade] (u) circle[radius=.42];
  \filldraw[bd shade] (z) circle[radius=.42];

  \draw[bd edge] (u)--(z);
  \node[bd label] at (0,.27) {$a$};

  \draw[bd edge] (u)--++(-.70,.75);
  \draw[bd edge] (u)--++(-.75,-.70);
  \draw[bd edge] (z)--++(.70,.75);
  \draw[bd edge] (z)--++(.75,-.70);

  \node[blackvertex] at (u) {};
  \node[blackvertex] at (z) {};
  \node[above=12pt] at (u) {$u$};
  \node[above=12pt] at (z) {$z$};

  \node[panel] at (0,-2.08)
    {(a) Before insertion};
\end{scope}

\draw[bd arrow] (-1.18,.10)--(1.18,.10);
\node[panel] at (0,.49) {adding a boundary};

\begin{scope}[xshift=4.2cm]
  \coordinate (u) at (-1.85,0);
  \coordinate (v) at (0,0);
  \coordinate (z) at (1.85,0);

  \filldraw[bd shade] (u) circle[radius=.42];
  \filldraw[bd shade] (v) circle[radius=.46];
  \filldraw[bd shade] (z) circle[radius=.42];

  \draw[bd edge] (u)--++(-.70,.75);
  \draw[bd edge] (u)--++(-.75,-.70);
  \draw[bd edge] (z)--++(.70,.75);
  \draw[bd edge] (z)--++(.75,-.70);

  \draw[bd edge] (u)--(v)--(z);
  \node[bd label] at (-.93,.27) {$e$};
  \node[bd label] at ( .93,.27) {$f$};

  \draw[bd edge]
    (v)
      .. controls (-.23,-.38) and (-.85,-.46)
      .. (-.76,-1.02)
      .. controls (-.69,-1.45) and (.69,-1.45)
      .. (.76,-1.02)
      .. controls (.85,-.46) and (.23,-.38)
      .. (v);

  \draw[bd corner]
    ($(v)+(-50:.31)$)
    arc[start angle=-50,end angle=-130,radius=.31];

  \node[blackvertex] at (u) {};
  \node[blackvertex] at (v) {};
  \node[blackvertex] at (z) {};

  \node[above=12pt] at (u) {$u$};
  \node[above=12pt] at (v) {$v$};
  \node[above=12pt] at (z) {$z$};

  \node[bd label] at (0,-.74) {$c_v$};
  \node[bd label] at (0,-1.57) {$x$};

  \node[panel] at (0,-2.08)
    {(b) $W_v=(x,e,f,x)$};
\end{scope}

\end{tikzpicture}%
}
\caption{Adding a boundary.  An edge $a$ is subdivided
into $e,f$, and a new loop $x$ is attached at $v$.
The clockwise vertex word is $W_v=(x,e,f,x)$,
starting immediately after the distinguished corner $c_v$.
The loop bounds the region corresponding to the new boundary
component.  The original boundary cycles are preserved.
The picture shows an edge with distinct endpoints; the
construction also applies when the original edge is a loop.}
\label{fig:classification-boundary}
\end{figure}
\begin{proof}
Let $G'$ be the graph after the insertion, and let $S,S'$
be the oriented thickenings of $G,G'$, respectively.

We examine the topology and boundary data.
Let $h,k$ be the two half-edges of the original edge,
with $\iota(h)=k$.  Denote the new half-edges at $v$
by $e_v,f_v,x_1,x_2$, ordered clockwise as
\[
  W_v=(x_1,e_v,f_v,x_2).
\]
The old half-edges $h,k$ become the other ends of $e,f$,
respectively, and
\[
  \iota'(h)=e_v,\qquad
  \iota'(k)=f_v,\qquad
  \iota'(x_1)=x_2.
\]
The successor permutation at $v$ is
\[
  \pi'_v=(x_1,e_v,f_v,x_2),
\]
while the cyclic orders at all old vertices are unchanged.

Write $\phi=\pi\iota$ and $\phi'=\pi'\iota'$ for the
boundary permutations.  Their transitions differ only
at the subdivided edge:
\[
\begin{aligned}
  h&\longmapsto f_v\longmapsto\phi(h),\\
  k&\longmapsto x_2\longmapsto e_v
       \longmapsto\phi(k),\\
  x_1&\longmapsto x_1.
\end{aligned}
\]
Thus each old boundary cycle is preserved, with additional
half-edge occurrences inserted into it, while $(x_1)$
is a new boundary cycle.  The distinguished corner at
$v$, between $x_2$ and $x_1$, lies on this new cycle.
All old distinguished corners retain their boundary
assignments, so $G'$ is boundary-admissible.

The operation adds one vertex and two edges.  Since it
also adds one boundary component, Euler's formula
\[
  |V(G)|-|E(G)|=2-2g-b
\]
shows that the genus is unchanged.

The new vertex has degree four and weight three.
Its distinguished corner is the only one on the new
boundary, so this boundary contributes the AAG pair
$(1,3)$.  All old vertex degrees and distinguished-corner
assignments are unchanged, hence so are the old AAG
pairs.  In particular, the new boundary winding is
\[
  w_{\mathrm{new}}=1-3=-2,
\]
and all old boundary windings are preserved.

We next compare winding numbers. Any closed edge walk in
$G$ has a natural continuation in $G'$, obtained by replacing
each traversal of the subdivided edge with a traversal of
$e$ followed by $f$, or of the same path in the reverse
direction. Let $\gamma$ and $\gamma'$ denote the closed
curves obtained by smoothing these walks in the respective
thickenings.

At the new vertex, such a walk turns between positions
$2$ and $3$, or between positions $3$ and $2$.  By
\eqref{eq:ribbon-local-turn},
\[
  \varepsilon(2,3)=\varepsilon(3,2)=0.
\]
All turns at old vertices have the same entry and exit
positions as before.  Consequently,
\eqref{eq:ribbon-winding-sum} gives
\[
  w_{\eta_{G'}}(\gamma')
  =w_{\eta_G}(\gamma)
\]
for corresponding smoothed closed walks.

We finally consider the internal line-field invariant.
The insertion creates no new handle, and the old graph
cycles still supply a geometric symplectic basis of the
capped surface.  Their intersection pairings and winding
numbers are unchanged.

For $g=1$, the genus-one gcd is preserved, since the
winding numbers of a handle basis and all old boundary
windings are unchanged, while the new boundary contributes
\[
  w_{\mathrm{new}}+2=0
\]
to \eqref{eq:genus-one-gcd}.

For $g\geq2$, winding parity is unchanged on the old
homology classes, and the new boundary has even winding.
Thus $\sigma(\eta_{G'})=\sigma(\eta_G)$.

Moreover, the Arf invariant is defined on $S'$ precisely
when it is defined on $S$, since
\[
  w_{\mathrm{new}}=-2\equiv2\pmod4.
\]
When defined, the corresponding symplectic basis curves
have the same winding numbers and hence the same quadratic
values
\[
  q([\gamma])
  =\frac{w_\eta(\gamma)}2+1\pmod2.
\]
Their intersection pairings also agree, so the Arf
invariants coincide.

Note that there is no internal datum in genus zero.  Therefore
the internal line-field invariant is preserved in
every genus.
\end{proof}

We next use vertex splitting to construct prescribed weight
profiles while keeping the boundary windings fixed.
A vertex splitting increases both the number of weights
and their sum by one, preserving their difference.
We therefore begin with weight profiles having the smallest possible
number of entries compatible with the prescribed boundary winding numbers.

\begin{definition}[Reduced weight profiles]
\label{def:classification-reduced-profiles}
Fix the boundary windings $w_1,\ldots,w_b$.
On boundary $\partial_i S$, a weight profile is a nonempty list
$(\lambda_1,\ldots,\lambda_r)$ of nonnegative integers
satisfying
\[
 r-\sum_{j=1}^r\lambda_j=w_i.
\]
The \emph{reduced weight profile} on boundary $\partial_iS$ is
the unique such list with the smallest possible number
of entries. Its length is $\max\{1,w_i\}$, and it is
given explicitly by
\begin{equation}
\label{eq:classification-reduced-profile}
 \begin{cases}
  (1-w_i),&w_i\leq1,\\[2pt]
  \underbrace{(0,\ldots,0)}_{w_i\text{ entries}},&w_i>1.
 \end{cases}
\end{equation}
Thus the first case consists of a single weight $1-w_i$,
whereas the second consists of $w_i$ zero weights.
\end{definition}

\begin{lemma}[Realisation from reduced profiles]
\label{lem:classification-profiles}
Let $G$ be a boundary-admissible ribbon graph whose weight
profile on each boundary is reduced.
Any prescribed weight profiles with the same boundary
windings can be realised by a ribbon graph obtained
from $G$ by a sequence of vertex splittings as in
\cref{lem:classification-splitting}.
These splittings preserve the genus, boundary windings,
and internal line-field invariant.
\end{lemma}

\begin{proof}
We work independently on each boundary $\partial_i S$, starting
with its prescribed weight profile.
As long as at least two weights remain and one is
positive, choose two weights $a,b$ with $a+b\geq1$
and replace them by $a+b-1$.
This numerical operation decreases both the number
of weights and their sum by one, so it preserves $w_i$.

The process terminates because the number of weights
strictly decreases. If only one weight remains, its
value is necessarily $1-w_i$. If at least two remain,
they are all zero, and their number is necessarily $w_i$.
Thus the terminal profile is the reduced profile
\eqref{eq:classification-reduced-profile}.

Now start from $G$ and reverse the recorded numerical
operations. Each reverse step replaces a weight
$a+b-1$ by two weights $a,b$.
The vertex carrying the weight $a+b-1$ has a vertex
word of length $a+b$. Split this word into consecutive
blocks of lengths $a$ and $b$, allowing either block
to be empty. By \cref{lem:classification-splitting},
the resulting vertices have weights $a$ and $b$.

Performing these splittings on each boundary produces
the prescribed profiles. The same lemma ensures that
the genus, boundary windings, and internal line-field
invariant are preserved throughout.
\end{proof}

\subsection{Higher-genus spectra}
\label{subsec:classification-higher-genus}

We now apply the preceding constructions to determine
the global-dimension spectra for $g\geq2$. Since global
dimension is the largest polygon weight, we seek weight
profiles with the prescribed AAG data that can be realised
with the required internal line-field invariant.

We first consider the cases $\sigma(\eta_A)=1$ or $b\geq2$.
In these cases, we can realise the reduced profiles with
the required line-field data and then use
\cref{lem:classification-profiles} to obtain every prescribed
weight profile. The remaining case of an even line field
on a surface with one boundary component will be treated
separately.

\begin{proposition}
\label{prop:classification-higher-general}
Suppose $g\geq2$.  If $\sigma(\eta_A)=1$ or $b\geq2$, then
\[
 \GSpec(A)=[L,U].
\]
In fact, every weight profile with these AAG data is realisable in the
prescribed derived class.
\end{proposition}

\begin{proof}
We prove that every weight profile with the prescribed AAG data
can be realised in the derived class of $A$.
The spectrum equality will then follow by considering the
possible maximum weights.

Suppose first that some boundary winding $w_i$ is odd, and
let $(\lambda_v)_{v\in V}$ be any prescribed weight profile.
By \cref{prop:Euler-identities},
\[
 \sum_{v\in V}\lambda_v
 =\sum_{i=1}^b n_i
 =m+4g+2b-4.
\]
Thus \cref{thm:all-genus-realization} gives a connected gentle
algebra $B$ of finite global dimension realising this profile
and satisfying $\AG(B)=\AG(A)$.
Its boundary windings are again the $w_i$, so
$\sigma(\eta_B)=\sigma(\eta_A)=1$.
No Arf datum is required in this case, and
\cref{prop:complete-invariant} therefore gives
$\Db(\modu B)\simeq\Db(\modu A)$.
This proves the required realisation statement whenever an
odd boundary winding occurs.

\smallskip
\noindent\emph{Reduction to smaller genus.}
Assume henceforth that all boundary windings are even.
By \cref{lem:classification-profiles}, it is enough to
construct a ribbon graph of genus $g$ with the reduced
profiles \eqref{eq:classification-reduced-profile} and
the internal line-field invariant of $A$.
Indeed, vertex splittings will then recover any prescribed
weight profile while preserving the genus, boundary windings,
and internal invariant.
For this construction, write $(\lambda_v)$ for the reduced
weights.

We shall first realise smaller weights in lower genus and
then restore the reduced weights by adding handles.
Since each handle insertion increases the selected weight
by four, we begin by determining how many multiples of four
can be subtracted while keeping all weights nonnegative.
Put
\[
 C=\sum_v\left\lfloor\frac{\lambda_v}{4}\right\rfloor,
 \qquad
 W_+=\sum_{w_i>0}w_i,
 \qquad
 t=\#\{i:w_i<0,\ w_i\equiv2\pmod4\}.
\]
By \eqref{eq:classification-reduced-profile}, a boundary
with $w_i<0$ carries the single reduced weight $1-w_i$.
Since $w_i$ is even, its contribution to $C$ is
\[
 \left\lfloor\frac{1-w_i}{4}\right\rfloor
 =
 \begin{cases}
  -w_i/4,&w_i\equiv0\pmod4,\\[2pt]
  (-w_i-2)/4,&w_i\equiv2\pmod4.
 \end{cases}
\]
For $w_i=0$, the reduced profile is $(1)$, whereas for
$w_i>0$ it consists entirely of zeros; both contribute
zero to $C$.
Summing these contributions gives
\[
 4C=-\sum_{w_i<0}w_i-2t
    =W_+-\sum_{i=1}^b w_i-2t.
\]
Using the Euler identity $\sum_iw_i=4-4g-2b$, we obtain
\begin{equation}
\label{eq:classification-capacity}
 C=g+\frac{2b-4+W_+-2t}{4}.
\end{equation}
If $t=b$, every boundary winding is negative and congruent
to two modulo four, so $W_+=0$ and $C=g-1$.
Otherwise $t\leq b-1$, and
\eqref{eq:classification-capacity} gives
$C\geq g-\frac12$.
Since $C$ is an integer, we then have $C\geq g$.

These bounds determine how far we can lower the genus.
More precisely, subtracting a total of $4h$ from the
weights, in multiples of four at individual vertices,
changes the Euler condition from genus $g$ to genus $g-h$.
Each subtraction of four increases the corresponding
boundary winding by four, so all boundary windings
retain their residues modulo four.
Adding handles at the same vertices subsequently restores
the original weights and boundary windings.

\smallskip
\noindent\emph{Realising the reduced profiles.}
We first treat $\sigma(\eta_A)=1$.
Since $C\geq g-1\geq1$, some reduced weight is at least four.
Subtract four from this weight.
The resulting nonnegative weights satisfy the Euler
condition for genus $g-1$, so
\cref{thm:all-genus-realization} gives a realisation.
At the modified vertex, apply
\cref{lem:classification-handle} with $r=0$.
This restores the weight and the genus and produces a
closed curve of odd winding.
The resulting graph therefore has the required reduced
profiles and satisfies $\sigma=1$, as desired.

It remains to consider $\sigma(\eta_A)=0$.
The hypothesis of the proposition then gives $b\geq2$.
We first construct an even line field in genus zero or
one, according to the preceding estimate for $C$, and
then restore the removed weights by handle insertions
that preserve evenness.

When $C\geq g$, subtract a total of $4g$ from the reduced
weights, in multiples of four, keeping all weights
nonnegative.
The resulting data satisfy the Euler condition for
genus zero and are realisable by
\cref{thm:all-genus-realization}.
Their boundary windings are still even.
Since winding parity factors through homology and
boundary classes generate first homology in genus zero,
the resulting line field is even.

The remaining possibility is $C=g-1$.
In this case every boundary winding is negative and
congruent to two modulo four, so each reduced weight
is congruent to three modulo four.
Subtracting all available multiples of four leaves
one weight three on each boundary.
The resulting data have genus one and boundary
windings all equal to $-2$.

To realise these data with an even line field, start
with $b=2$ and take two vertices with words
\[
 (1,2,3,4),\qquad(1,2,3,4),
\]
where equally labelled half-edges are paired.
The boundary permutation has two cycles, each containing
one distinguished corner.
With two vertices and four edges, the thickening has
genus one, and each boundary carries one weight three.
Every two-edge cycle has winding zero, because its
local contributions at the two vertices cancel.
These cycles generate first homology, so the line
field is even.
Repeated applications of
\cref{lem:classification-boundary} add the remaining
$b-2$ boundary components, each carrying one weight
three, while preserving evenness.

We have thus obtained an even starting line field in
both cases.
Restore the removed weights by adding handles at the
corresponding vertices.
Taking $r=1$ at each insertion preserves evenness by
\cref{lem:classification-handle} and gives a genus-$g$
realisation of the original reduced profiles.
If no Arf datum is required, this already gives the
desired internal invariant.

When an Arf datum is required, we use the choice of
the last handle insertion to obtain its prescribed value.
In this case $w_i\equiv2\pmod4$ for every $i$.
A vertex at which a handle is inserted must lie on a
boundary with negative winding, since the reduced
weights on positive boundaries are all zero.
Its weight remains congruent to three modulo four
throughout the construction and is therefore at least
three.
Its vertex word consequently has length at least four,
so both $r=1$ and $r=3$ are available.

Keep $r=1$ at all earlier insertions and choose either
$r=1$ or $r=3$ at the last one.
Before this last insertion, the surface has genus
$g-1\geq1$.
The even line field and the boundary congruences define
a quadratic form on the homology of its capped surface.
The quadratic-form calculation in the proof of
\cref{lem:classification-handle} therefore applies:
the choice $r=1$ preserves its Arf invariant, whereas
$r=3$ changes it by one.
This also applies when $g=2$, since the same calculation
is valid for the capped genus-one surface.
We can thus obtain either Arf value, and in particular
the one prescribed by $\eta_A$.

\smallskip
\noindent\emph{Recovering the prescribed profiles.}
We have now realised the reduced profiles with genus
$g$, boundary windings $w_i$, and the required internal
line-field invariant.
Applying \cref{lem:classification-profiles} produces
every target weight profile with the prescribed AAG
data while preserving this invariant.
By \cref{prop:complete-invariant}, the resulting
algebras belong to the derived class of $A$.

Together with the odd-boundary case, this proves the
assertion for every weight profile.
The global dimension of each resulting algebra is
the maximum of its weights, and
\cref{lem:bounded-distributions} shows that these
maxima attain every integer in $[L,U]$.
The reverse inclusion follows from
\eqref{eq:intro-enclosure}.
Hence $\GSpec(A)=[L,U]$.
\end{proof}
It remains to consider $b=1$ and $\sigma(\eta_A)=0$.
The following lemma gives a necessary condition on the
weight profile.

\begin{lemma}
\label{lem:one-boundary-parity}
Let $G$ be a boundary-admissible ribbon graph with distinguished corners
whose thickening has one boundary component.  Suppose that every closed
curve has even winding number relative to $\eta_G$.  Then $G$ has a vertex
of odd degree.  If $G$ has exactly two vertices, both have odd degree, and
$\gldim A_G$ is even.
\end{lemma}

\begin{proof}
Suppose, for a contradiction, that every vertex has
even degree. Colour the half-edge occurrences at each
vertex alternately by $0$ and $1$ in their cyclic order.
Since the degree is even, the last and first half-edges
also have opposite colours.

We may interchange the two colours independently at
each vertex. Choose a spanning tree of $G$ and fix
the colouring at a root vertex. Proceeding outwards
along the tree, choose the colouring at each new
vertex so that the two half-edges of every tree edge
have opposite colours.

We claim that the two half-edges of every remaining
edge also have opposite colours. Write $c(h)$ for
the colour of a half-edge $h$.
By \eqref{eq:ribbon-local-turn}, a turn between
half-edges $a,b$ at the same vertex contributes
\[
  1+c(a)+c(b)\pmod2
\]
to the winding number. This expression is unchanged
if both colours at that vertex are interchanged.
For a closed edge walk, the number of turns equals
the number of edge traversals, and each traversed
edge contributes the colours of its two half-edges.
Thus its winding parity is the sum of
$1+c(h)+c(\iota(h))$ over the traversed edges.

Now let $e=\{h,\iota(h)\}$ be an edge outside the tree,
and let $\gamma_e$ be the closed curve obtained by
smoothing the cycle formed by $e$ and the tree path
between its endpoints. Every edge in this cycle
other than $e$ is a tree edge, so its contribution
to the parity sum is zero. Since $\gamma_e$ has
even winding, we obtain
\[
  0\equiv w_{\eta_G}(\gamma_e)
   \equiv1+c(h)+c(\iota(h))\pmod2.
\]
The half-edges of $e$ therefore have opposite colours,
proving the claim. This also covers loops, for which
the tree path is empty.

Consequently, both the edge involution $\iota$ and
the vertex successor permutation $\pi$ reverse
colours. The boundary permutation $\pi\iota$ therefore
preserves colours, so each of its cycles contains
only one colour. Both colours occur, and hence there
are at least two boundary cycles. This contradicts
the assumption that the thickening has one boundary
component. Thus some vertex has odd degree.

If there are exactly two vertices, their degree sum
is even, so both degrees must be odd. Consequently,
$\gldim A_G=\degmax(G)-1$ is even.
\end{proof}

We therefore use two-vertex profiles as starting points.
Their weights sum to $4g$, and the preceding lemma
requires both weights to be even. The next lemma
realises every such pair with either Arf value.
These realisations will then be extended to more
general profiles by vertex splitting.

\begin{lemma}[Two even weights]
\label{lem:classification-two-even}
For every $g\geq2$ and every pair of nonnegative even integers
$\lambda_1,\lambda_2$ satisfying $\lambda_1+\lambda_2=4g$, there are
one-boundary ribbon realisations with these weights, $\sigma=0$,
and either Arf value.
\end{lemma}
\begin{proof}
At genus two the following seeds suffice, together with their swaps.
Equal labels are paired to form edges.
\begin{center}
\small
\renewcommand{\arraystretch}{1.12}
\begin{tabular}{@{}ccll@{}}
$\Arf$ & Weights & First word & Second word\\ \hline
$0$ & $(4,4)$ & $(3,2,1,3,2)$ & $(4,5,1,4,5)$\\
$0$ & $(2,6)$ & $(1,2,3)$ & $(4,5,1,2,3,4,5)$\\
$0$ & $(0,8)$ & $(1)$ & $(2,3,1,2,4,5,3,4,5)$\\
$1$ & $(4,4)$ & $(1,2,3,4,5)$ & $(1,2,3,4,5)$\\
$1$ & $(2,6)$ & $(1,2,3)$ & $(1,2,4,5,3,4,5)$\\
$1$ & $(0,8)$ & $(1)$ & $(2,3,4,5,1,2,3,4,5)$
\end{tabular}
\end{center}
All six boundary permutations have one cycle.  Use edge $1$ as a
spanning tree and smooth the fundamental walks associated with
edges $2,3,4,5$.  Applying \eqref{eq:ribbon-winding-sum} on each
original graph, with suitable orientations, gives the winding
vectors below, listed in the order of the six rows:
\[
\begin{gathered}
 (2,2,2,2),\quad(0,0,4,4),\quad(2,4,2,2),\\
 (0,0,0,0),\quad(0,2,2,2),\quad(4,4,4,4).
\end{gathered}
\]
Only after these winding computations do we contract edge $1$
to read the intersection pairing from the alternating half-edges.  For a quadratic form on the four-dimensional capped homology,
\[
 \sum_{x\in H_1(\overline S;\mathbb F_2)}(-1)^{q(x)}
 =(-1)^{\Arf(q)}4.
\]
The six sums are respectively $4,4,4,-4,-4,-4$, verifying the
displayed Arf values.
An $r=1$ handle adds four to either weight and preserves Arf.
Every even pair of sum $4g$ reduces, by subtracting four from suitable
entries, to a genus-two pair in the table or its swap.
\end{proof}

Combining the parity restriction with these two-vertex
constructions and vertex splitting gives the following
criterion for one-boundary profiles and the resulting
global-dimension spectra.

\begin{proposition}
\label{prop:classification-one-boundary}
Suppose $g\geq2$, $b=1$, and $\sigma(\eta_A)=0$.  A weight profile is
realisable in the derived class of $A$ if and only if it has an even
entry.  In particular, $m\geq2$.  If $m=2$, then
\begin{equation}
\label{eq:classification-gap-spectrum}
 \GSpec(A)=\{2g,2g+2,\ldots,4g\}.
\end{equation}
If $m\geq3$, put $N=m+4g-2$ and $L=\lceil N/m\rceil$.  Then
\begin{equation}
\label{eq:classification-one-boundary-interval}
 \GSpec(A)=[L+\varepsilon,N],\qquad
 \varepsilon=
 \begin{cases}
  1,&N=mL\text{ and }L\text{ is odd},\\
  0,&\text{otherwise}.
 \end{cases}
\end{equation}
\end{proposition}
\begin{proof}
By \cref{prop:Euler-identities}, the unique AAG pair is
$(m,N)$, where
\[
 N=m+4g-2,
 \qquad
 w_1=m-N=2-4g.
\]
In particular, $w_1\equiv2\pmod4$, so the Arf invariant
of $\eta_A$ is defined.

\smallskip
\noindent\emph{The realisation criterion.}
If a profile is realised in the derived class of $A$,
then \cref{lem:one-boundary-parity} gives a vertex of
odd degree in its ribbon graph.
Since the corresponding weight is one less than this
degree, the profile has an even entry.
For $m=1$, the only weight would be $4g-1$, which is
odd. Thus $m\geq2$.

Conversely, let $(\lambda_1,\ldots,\lambda_m)$ be a
profile with the prescribed AAG data and at least
one even entry.
We first reduce it to a two-entry profile to which
\cref{lem:classification-two-even} applies.
Choose an even entry $\lambda\leq4g$.
Such an entry exists: if a zero occurs, take
$\lambda=0$; otherwise every entry is at most
\[
 N-(m-1)=4g-1.
\]

Keep $\lambda$ fixed and merge the other entries
using the operation
\[
 (a,b)\longmapsto a+b-1.
\]
Each merge decreases both the number of remaining
entries and their sum by one.
Consequently, when $r$ of these entries remain,
their sum is
\[
 4g-\lambda+r-1.
\]
If $r\geq2$, this sum is positive, so at least one
entry is positive and can be merged with another
without producing a negative weight.
The process therefore continues until a single
entry remains, whose value is $4g-\lambda$.
We have thus reduced the target profile to the
nonnegative even pair
\[
 (\lambda,4g-\lambda).
\]

By \cref{lem:classification-two-even}, this pair has
a one-boundary ribbon realisation of genus $g$ with
$\sigma=0$ and Arf invariant $\Arf(\eta_A)$.
Reverse the recorded merges by applying
\cref{lem:classification-splitting}.
For the inverse of a merge $(a,b)\mapsto a+b-1$,
split the vertex word of length $a+b$ into consecutive
blocks of lengths $a$ and $b$; either block is allowed
to be empty.
These splittings recover the target profile while
preserving the genus, boundary winding, and internal
line-field invariant.
The resulting algebra has AAG invariant $[(m,N)]$
and the same internal invariant as $A$.
Hence \cref{prop:complete-invariant} places it in the
derived class of $A$, proving the criterion.

\smallskip
\noindent\emph{The spectrum.}
We now determine the possible maximum weights among
the profiles satisfying this criterion.

If $m=2$, the weights sum to $4g$.
The presence of an even entry therefore forces both
entries to be even.
Their maximum is an even integer between $2g$ and
$4g$, and every such value $d$ occurs in the profile
$(4g-d,d)$.
This proves \eqref{eq:classification-gap-spectrum}.

Now suppose that $m\geq3$, and put $L=\lceil N/m\rceil$.
Every profile has maximum in $[L,N]$.
Fix an integer $d$ in this interval.
By \cref{lem:bounded-distributions}, there is a profile
$(\lambda_1,\ldots,\lambda_m)$ satisfying
\[
 0\leq\lambda_j\leq d,
 \qquad
 \sum_{j=1}^m\lambda_j=N,
 \qquad
 \lambda_1=d.
\]
If this profile has an even entry, it is realisable
by the criterion just proved.

Otherwise all entries are positive and odd.
If $N<md$, at least one of the remaining entries is
less than $d$; after relabelling, assume that
$\lambda_2<d$.
Since $m\geq3$, we may replace
\[
 (\lambda_2,\lambda_3)
 \quad\text{by}\quad
 (\lambda_2+1,\lambda_3-1).
\]
Both new entries are even and lie in $[0,d]$.
The total weight is unchanged, and the maximum remains
$d$ because $\lambda_1=d$.
Thus the modified profile is realisable.

The only excluded possibility is therefore
$N=md$ with $d$ odd.
In this case the unique profile with sum $N$ and
maximum $d$ is $(d,\ldots,d)$, which has no even entry.
Moreover, $N=md$ forces $d=N/m=L$.
Hence the only possible exclusion is the lower
endpoint $L$, and it is excluded precisely when
$N=mL$ and $L$ is odd.
This is exactly the correction $\varepsilon$ in
\eqref{eq:classification-one-boundary-interval},
and proves that formula.
\end{proof}

\subsection{A winding bound in genus one}
\label{subsec:classification-winding-bound}

Parity and the Arf invariant settle the higher-genus cases.  In genus
one, the remaining obstruction is divisibility of winding numbers.
To control it, we need a bound that also detects the weight profile.  The relevant quantity is the number of
leaves of the ribbon graph, or equivalently the number of zero weights.

\begin{lemma}[The leaf bound]
\label{lem:classification-leaf-bound}
Let $G$ have genus one and $b$ boundary components, all with winding
$-2$.  If $\ell$ is its number of degree-one vertices and
$D=d(\eta_G)>0$, then
\begin{equation}
\label{eq:classification-leaf-bound}
 D\leq\ell+b.
\end{equation}
\end{lemma}
\begin{proof}
We first reduce to a graph with one boundary component
and only degree-one and degree-three vertices, while
preserving $D$ and the quantity $\ell+b$.
We then bound the windings of a geometric handle basis
in the resulting graph.

\smallskip
\noindent\emph{Reducing the vertex degrees.}
Suppress all bivalent vertices.
Their local turns contribute zero, so this preserves
winding numbers and creates no new leaf.
Boundary admissibility is also preserved.
Indeed, $w_i=-2$ gives $n_i=m_i+2$, so each boundary
has a distinguished corner at a vertex of degree at
least three; these corners survive the suppression.

Next, whenever a vertex word has length $q\geq4$,
split it after its first two occurrences.
The resulting vertices have degrees $3$ and $q-1$,
so no leaf is introduced.
Repeating \cref{lem:classification-splitting} produces
a graph whose vertices have degrees one or three,
without changing the genus, boundary windings, $D$,
or the number of leaves.
We continue to denote this graph by $G$.

\smallskip
\noindent\emph{Reducing the number of boundary components.}
Suppose that $b>1$.
We first find an edge with its two sides on distinct
boundary cycles and with one half-edge occurring first
or last in a trivalent vertex word.

If no such half-edge existed, then at every trivalent
vertex the first and last edges would each have the
same boundary cycle on both sides.
The two sectors adjacent to the distinguished corner
would therefore belong to the same boundary cycle as
that corner.
Thus all three sectors at the vertex would belong
to one boundary cycle.
A leaf has only one sector, and connectedness would
then force all sectors of the graph to belong to the
same boundary cycle, contrary to $b>1$.

Choose such an edge $e$ and such a half-edge $h$ at
a trivalent vertex $v$.
Detach $h$ from $v$ and attach it to a new leaf.
At $v$, retain the other two occurrences in their
original linear order; the new leaf has word $(h)$.
All half-edge pairings remain unchanged.

The edge $e$ is not a bridge, since the two sides of
a bridge belong to the same boundary cycle.
Hence the modified graph is connected.
If $\pi'$ denotes its vertex permutation, then
\[
 \pi'=(h\ \pi(h))\pi,
 \qquad
 \pi'\iota=(h\ \pi(h))(\pi\iota).
\]
The half-edges $h$ and $\pi(h)$ lie on the two distinct
boundary cycles adjacent to $e$, so this transposition
merges those cycles.
The edge count is unchanged and the vertex count
increases by one; since the boundary count decreases
by one, the genus remains one.

We next check the boundary windings.
Suppose the merged boundaries have data
$(m_i,n_i)$ and $(m_j,n_j)$.
Because $h$ was first or last, the distinguished corner
at $v$ belongs to one of these boundaries.
The weight at $v$ decreases from two to one, and the
new leaf contributes one distinguished corner of
weight zero.
Thus, before suppressing $v$, the merged boundary has
marked count and weight sum
\[
 m_i+m_j+1,
 \qquad
 n_i+n_j-1,
\]
respectively.
Its winding is therefore
\[
 w_i+w_j+2=-2.
\]
All other boundary windings are unchanged.

To compare the genus-one invariants, let $H=G-e$.
In the cellular embedding of $G$ in the capped torus,
deleting $e$ merges its two disc faces into a disc.
Thus $H$ is still connected and cellularly embedded
in a torus, and it carries a geometric handle basis.

View the curves of this basis in the original graph
and in the graph obtained by detachment, before
suppressing $v$.
The only changed vertex word is the one at $v$.
Since the removed occurrence was first or last,
all position differences between the remaining
occurrences are unchanged.
Hence \eqref{eq:ribbon-local-turn} gives the same
local contributions to these curves in both graphs,
and their winding numbers agree.
All shifted boundary windings are zero, so
\eqref{eq:genus-one-gcd} computes the same $D$ from
these two basis windings; see also the proof of
\cite[Theorem~7.4]{AmiotPlamondonSchroll}.

Finally, suppress the now-bivalent vertex $v$.
This preserves winding numbers and subtracts one
from both the marked count and the weight sum of
the merged boundary.
The resulting graph is again boundary-admissible,
has only degree-one and degree-three vertices,
and satisfies
\[
 (b,\ell,D)\longmapsto(b-1,\ell+1,D).
\]
In particular, $\ell+b$ is unchanged.

Repeating this operation $b-1$ times gives a graph
$G_0$ of genus one with one boundary component,
boundary winding $-2$, and
\[
 \ell'=\ell+b-1
\]
leaves.
Its genus-one invariant is still $D$.

\smallskip
\noindent\emph{Bounding the winding.}
Write $\eta_0:=\eta_{G_0}$, and let $T$ be the number
of trivalent vertices of $G_0$.
The degree sum and Euler identities give
\[
 T+\ell'-\frac{3T+\ell'}2=-1,
 \qquad\text{hence}\qquad
 T=\ell'+2.
\]
At a trivalent vertex, the local winding contribution
is zero between adjacent positions, while
\[
 \varepsilon(1,3)=1,
 \qquad
 \varepsilon(3,1)=-1.
\]
A simple graph cycle visits each vertex at most once
and cannot visit a leaf.
Its smoothing $\gamma$ therefore satisfies
\[
 \bigl|w_{\eta_0}(\gamma)\bigr|\leq T.
\]

If equality holds, then $\gamma$ passes through every
trivalent vertex, and all local contributions are
nonzero and have the same sign.
The unused middle half-edge at each vertex consequently
lies on the same side of the oriented curve.
Since all other vertices are leaves, every edge outside
$\gamma$ is attached on that side.
On the opposite side, a ribbon neighbourhood of
$\gamma$ extends to an annular collar of a boundary
component.
Thus $\gamma$ is boundary-parallel and hence separating.
It follows that every nonseparating simple graph cycle
has absolute winding at most $T-1$.

Choose a spanning tree $R$ of $G_0$.
Since
\[
 |E(G_0)|-|V(G_0)|+1=2,
\]
there are exactly two edges outside $R$, and their
fundamental cycles are simple graph cycles.
A regular neighbourhood of $R$ is a disc.
The two remaining edge bands attach to this disc
with alternating ends, because the thickening has
genus one and one boundary component.
Consequently, the smoothings of the two fundamental
cycles can be isotoped to meet once.
After choosing orientations, they form a geometric
handle basis $\alpha,\beta$.
In particular, both are nonseparating, and their
absolute windings are at most $T-1$.

Since the boundary winding is $-2$,
\eqref{eq:genus-one-gcd} gives
\[
 D=\gcd\bigl(
     w_{\eta_0}(\alpha),w_{\eta_0}(\beta)
    \bigr).
\]
As $D>0$, at least one of these windings is nonzero.
It is divisible by $D$, and hence
\[
 D\leq T-1=\ell'+1=\ell+b,
\]
as required.
\end{proof}
\subsection{Genus-one spectra and completion of the classification}
\label{subsec:classification-genus-one}

We now compute the spectrum for each value of the genus-one gcd.

\begin{proposition}
\label{prop:classification-torus}
Suppose $g=1$, and put $D=d(\eta_A)$.
\begin{enumerate}[label=\textup{(\roman*)}]
\item If some $w_i\neq-2$, then $\GSpec(A)=[L,U]$.
\item If all $w_i=-2$ and $D=0$, then
\begin{equation}
\label{eq:classification-torus-zero}
 \GSpec(A)=
 \begin{cases}
  [2,m],&b=1,\\
  [L,U],&b\geq2.
 \end{cases}
\end{equation}
The first case requires $m\geq2$.
\item If all $w_i=-2$ and $D>0$, then $1\leq D\leq m$ and
\begin{equation}
\label{eq:classification-torus-positive}
 \GSpec(A)=
 \left\{d\in[L,U]:
    \sum_{i=1}^b\left\lceil\frac{m_i+2}{d}\right\rceil
       \leq m-D+b\right\}.
\end{equation}
This set is a nonempty terminal interval of $[L,U]$.  For $b=1$ it is
\begin{equation}
\label{eq:classification-torus-one-boundary}
 \left[\left\lceil\frac{m+2}{m-D+1}\right\rceil,m+2\right].
\end{equation}
\end{enumerate}
In particular, every genus-one spectrum is an interval.
\end{proposition}

\begin{proof}
We use \cref{prop:complete-invariant} throughout:
once the prescribed AAG data and the genus-one invariant
$D$ have been realised, the resulting algebra belongs
to the derived class of $A$.
Its global dimension is the maximum weight.

We shall also use the following consequence of
\cref{lem:classification-splitting}.
A weight $\lambda$ can be replaced by any prescribed
$r$ nonnegative weights whose sum is $\lambda+r-1$.
Indeed, merge these weights numerically by
$(a,b)\mapsto a+b-1$.
When $k\geq2$ entries remain, their sum is
$\lambda+k-1>0$, so a positive entry can be merged
with another without producing a negative weight.
The final entry is $\lambda$, and reversing the
merges gives the required vertex splittings.

\smallskip
\noindent\emph{(i) Some boundary winding differs from $-2$.}
Put
\[
 H=\gcd_i(w_i+2)>0.
\]
By \eqref{eq:genus-one-gcd}, $D$ is a positive divisor
of $H$.
The Euler identity gives
\[
 \sum_i(w_i+2)=0.
\]
Since these integers are not all zero, some $w_j+2$
is a negative multiple of $H$, and hence
\[
 w_j+2\leq-H\leq-D.
\]
The reduced profile on this boundary consists of
the single weight
\[
 \lambda_j=1-w_j\geq D+3.
\]

Decrease $\lambda_j$ by four.
The resulting nonnegative weights satisfy the Euler
condition for genus zero, so they admit a realisation
by \cref{thm:all-genus-realization}.
At the modified vertex, the degree is
$\lambda_j-3\geq D$, and a consecutive block of length
$r=D-1$ can therefore be chosen.
Apply \cref{lem:classification-handle} with this value
of $r$.
The insertion restores the original reduced profiles
and all the original boundary windings, while
\eqref{eq:classification-handle-gcd} gives
\[
 d(\eta)
 =\gcd(D,w_1+2,\ldots,w_b+2)
 =D.
\]
By \cref{lem:classification-profiles}, vertex splittings
now produce every target weight profile with the
prescribed AAG data and the same invariant $D$.
Together with \cref{lem:bounded-distributions} and
\eqref{eq:intro-enclosure}, this proves
$\GSpec(A)=[L,U]$.

For the rest of the proof, assume that every boundary
winding is $-2$.
Then
\[
 n_i=m_i+2
 \qquad(1\leq i\leq b).
\]

\smallskip
\noindent\emph{(ii) The case $D=0$.}
Suppose first that $b\geq2$.
The two vertex words
\[
 (1,2,3,4),\qquad(1,2,3,4)
\]
give a genus-one graph with two boundary components,
each carrying a single weight three.
A handle basis can be chosen among its two-edge
cycles, whose local contributions cancel, so its
genus-one invariant is zero.
Apply \cref{lem:classification-boundary} repeatedly
to add the remaining $b-2$ boundaries.
This gives the reduced profile consisting of one
weight three on each boundary, still with $D=0$.
By \cref{lem:classification-profiles}, every target
profile can therefore be realised.
The spectrum is $[L,U]$ by
\cref{lem:bounded-distributions} and
\eqref{eq:intro-enclosure}.

Now let $b=1$, so the total weight is $m+2$.
We first show that no weight can exceed $m$.
Otherwise all the other weights would have sum
at most one.
After suppressing the bivalent vertices, there would
be a single vertex of degree at least three, with
all remaining vertices being leaves.
The first Betti number is two, so the graph consists
of two loops at this vertex together with pendant
edges.
Genus one forces the endpoints of the two loops to
alternate in the cyclic order, and the loops form
a handle basis.
By \eqref{eq:ribbon-local-turn}, a loop has winding
zero only when its two occurrences are adjacent
in the linear vertex word.
The alternating endpoints prevent both loop windings
from vanishing, contradicting $D=0$.

In particular, $m=1$ is impossible, since its only
weight would be three.
For $m\geq2$, the total weight $m+2$ also forces the
maximum to be at least two.
Thus the spectrum is contained in $[2,m]$.

To realise every value in this interval, start from
the two words
\[
 (1,2,3),\qquad(1,2,3),
\]
which give one boundary, weights $(2,2)$, and $D=0$.
Keep one weight two and replace the other by $m-1$
nonnegative weights of total $m$.
By \cref{lem:bounded-distributions}, their maximum
can be any prescribed integer $d\in[2,m]$.
The splitting observation at the beginning of the
proof realises each such choice while preserving $D$.
The resulting profile has $m$ entries, total weight
$m+2$, and maximum $d$.
This proves \eqref{eq:classification-torus-zero}.

\smallskip
\noindent\emph{(iii) The case $D>0$.}
We first establish a more precise realisation criterion.
For a target profile with the prescribed AAG data,
let $\ell_i$ be the number of zero entries on the
$i$th boundary and put $\ell=\sum_i\ell_i$.
We claim that this profile is realisable in the
derived class of $A$ if and only if
\[
 D\leq\ell+b.
\]
Necessity is \cref{lem:classification-leaf-bound}.
Moreover, each boundary has positive total weight
$m_i+2$, so $\ell\leq m-b$.
Applying the leaf bound to a profile of $A$ gives
\[
 1\leq D\leq m.
\]

To prove sufficiency, first suppose that $D\geq b$
and that $D\leq\ell+b$.
Choose integers $p_i$ such that
\begin{equation}
\label{eq:classification-necklace-choice}
 1\leq p_i\leq\ell_i+1,
 \qquad
 \sum_i p_i=D.
\end{equation}
Such a choice exists because the lower bounds sum
to $b$ and the upper bounds sum to $\ell+b$.

Construct a cyclic necklace of $b$ vertices, with
two edges $a_i,b_i$ between vertex $i$ and vertex
$i+1$, using cyclic indices.
At vertex $i$, take the word
\begin{equation}
\label{eq:classification-necklace}
 (a_{i-1},b_{i-1},
  c_{i,1},\ldots,c_{i,p_i-1},
  a_i,b_i),
\end{equation}
where each $c_{i,j}$ is a distinct edge ending at
a new leaf with singleton word $(c_{i,j})$.
For $b=1$, the $a$- and $b$-edges are loops.

Tracing the boundary permutation shows that the
$i$th boundary contains the distinguished corner
at vertex $i$ and those at its $p_i-1$ attached leaves.
These are the $b$ boundary cycles of the graph.
Since vertex $i$ has degree $p_i+3$, this boundary
has one weight $p_i+2$ and $p_i-1$ zero weights.
Its AAG pair is therefore $(p_i,p_i+2)$.
The graph has $D$ vertices and $D+b$ edges, so
Euler's formula gives genus one.

We next check its genus-one invariant.
For $b\geq2$, let $\alpha$ be the longitudinal cycle
following the $a$-edges.
At vertex $i$, it enters through $a_{i-1}$ and leaves
through $a_i$, at positions $1$ and $p_i+2$.
Thus
\[
 w_{\eta_G}(\alpha)
 =\sum_i\varepsilon(1,p_i+2)
 =\sum_i p_i
 =D.
\]
The two-edge cycle $\beta=a_i b_i^{-1}$ has winding
zero, since its turns are between adjacent positions
at both vertices.
The union of these two cycles consists of the three
paths $a_i$, $b_i$, and the remaining $a$-path.
Their cyclic orders at the two endpoints agree,
so their regular neighbourhood is a torus with one
boundary component.
Consequently, the smoothings of $\alpha$ and $\beta$
can be arranged to meet once and form a geometric
handle basis.
Since all boundary windings are $-2$, the gcd is
$\gcd(D,0)=D$.

For $b=1$, the two loops have alternating endpoints
and form a handle basis.
Their position pairs are $(1,D+2)$ and $(2,D+3)$,
so each can be oriented to have winding $D$.
The gcd is again $D$.

We now recover the prescribed profile.
By \eqref{eq:classification-necklace-choice}, the
$i$th target boundary has at least $p_i-1$ zero
entries.
Identify these with the existing zero weights and
keep them fixed.
The remaining target entries number
\[
 r_i=m_i-p_i+1
\]
and have total weight
\[
 m_i+2=(p_i+2)+(r_i-1).
\]
The splitting observation at the beginning of the
proof therefore replaces the weight $p_i+2$ by
exactly these remaining entries.
Doing this on every boundary gives the target AAG
data and target profile while preserving $D$.
The resulting algebra is derived equivalent to $A$.

If $D<b$, apply the same necklace construction with
$D$ boundary components and all $p_i=1$.
It has one weight three on each boundary and gcd $D$.
Add $b-D$ boundary components using
\cref{lem:classification-boundary}.
The resulting graph has the reduced profile
consisting of one weight three on each of the $b$
boundaries, still with gcd $D$.
By \cref{lem:classification-profiles}, every target
profile is realisable.
Since $D<b\leq\ell+b$, this completes the proof of
the claimed profile criterion.

We now translate this criterion into a condition
on the maximum weight.
Fix an integer $d\in[L,U]$ and put
\[
 q_i=\left\lceil\frac{m_i+2}{d}\right\rceil.
\]
Because $d\geq L$, we have $1\leq q_i\leq m_i$.
Any profile with maximum $d$ needs at least $q_i$
positive entries on the $i$th boundary.
Thus, if it is realisable and has $\ell$ zero entries,
\begin{equation}
\label{eq:classification-positive-slots}
 \sum_i\left\lceil\frac{m_i+2}{d}\right\rceil
 \leq m-\ell
 \leq m-D+b.
\end{equation}

Conversely, the largest possible number of zero
entries is attained by taking, on the $i$th boundary,
the profile
\[
 \left(
  \underbrace{d,\ldots,d}_{q_i-1},
  m_i+2-(q_i-1)d,
  \underbrace{0,\ldots,0}_{m_i-q_i}
 \right).
\]
The middle entry lies in $[1,d]$, so this profile
has exactly $q_i$ positive entries.
On a boundary with $m_i+2=U\geq d$, one of the
entries is $d$.
Hence the global maximum is exactly $d$, and the
total number of zero entries is
\[
 \ell=m-\sum_iq_i.
\]
The profile criterion just proved therefore gives
\[
 d\in\GSpec(A)
 \quad\Longleftrightarrow\quad
 \sum_i\left\lceil\frac{m_i+2}{d}\right\rceil
 \leq m-D+b.
\]
Together with \eqref{eq:intro-enclosure}, this proves
\eqref{eq:classification-torus-positive}.

The sum in this condition is nonincreasing in $d$.
At $d=U$, it equals $b$, and
$b\leq m-D+b$ because $D\leq m$.
Thus the allowed set is a nonempty terminal interval
of $[L,U]$.

Finally, when $b=1$, the condition becomes
\[
 \left\lceil\frac{m+2}{d}\right\rceil
 \leq m-D+1
 \quad\Longleftrightarrow\quad
 d\geq
 \left\lceil\frac{m+2}{m-D+1}\right\rceil.
\]
Since $1\leq m-D+1\leq m$, this lower bound is at
least $L$.
With $U=m+2$, we obtain
\eqref{eq:classification-torus-one-boundary}.
\end{proof}
\begin{example}
\label{ex:classification-interval-endpoints}
The interval property does not imply equality with the AAG interval.
For $\AG(A)=[(3,5)]$ and $D=2$, the genus-one formula gives
$\GSpec(A)=[3,5]$, while $[L,U]=[2,5]$.
In genus two, $\AG(A)=[(3,9)]$ with $\sigma=0$ also gives
$\GSpec(A)=[4,9]$ instead of $[3,9]$: the only profile with maximum
three is $(3,3,3)$, which violates the even-weight condition.
These endpoint exclusions differ from the internal gaps occurring
when there are exactly two weights.
\end{example}

\begin{proof}[Proof of \cref{thm:intro-complete-spectrum}]
Part~\textup{(i)} is \cref{thm:intro-derived-interval}.

For part~\textup{(ii)}, apply \cref{prop:classification-torus}.
In its positive-$D$ case, the left-hand side of the inequality in
\eqref{eq:intro-spectrum-LD} is nonincreasing in $r$.
Moreover, $L\geq2$, and the inequality holds at $r=U$ because its
left-hand side is $b$ and $D\leq m$.
Consequently, the membership condition in
\eqref{eq:classification-torus-positive} defines exactly $[L_D,U]$.
The one-boundary formula is
\eqref{eq:classification-torus-one-boundary}.

Part~\textup{(iii)} follows from
\cref{prop:classification-higher-general,prop:classification-one-boundary}.
For $b=1$, the Euler identity gives $U=n_1=m+4g-2$, so the condition
defining $\varepsilon$ agrees with
\eqref{eq:classification-one-boundary-interval}.
The same higher-genus propositions show that the formulas are independent
of the Arf invariant whenever it is defined.
\end{proof}

\begin{proof}[Proof of \cref{thm:noninterval-classification}]
In genus zero, \cref{thm:intro-derived-interval} gives $[L,U]$.
Every genus-one spectrum is an interval by
\cref{prop:classification-torus}.  For $g\geq2$,
\cref{prop:classification-higher-general,prop:classification-one-boundary}
leave exactly the case $b=1$, $m=2$, $\sigma=0$.
The Euler identity then gives $\AG(A)=[(2,4g)]$, and
\eqref{eq:classification-gap-spectrum} computes its spectrum.
Both Arf values occur by \cref{lem:classification-two-even} and
distinguish exactly two derived classes by
\cref{prop:complete-invariant}.
\end{proof}

The classification gives the following explicit counterexample
with five simple modules.

\begin{figure}[htbp]
\centering
\begin{tikzpicture}[
  x=2cm,
  >=Latex,
  line width=.7pt,
  font=\small,
  vertex/.style={circle,inner sep=2pt}
]
  \foreach \i in {1,...,5}{
    \node[vertex] (v\i) at ({\i-1},0) {$\i$};
  }

  \foreach \i in {1,...,4}{
    \pgfmathtruncatemacro{\j}{\i+1}
    \draw[->]
      (v\i) to[bend left=22]
      node[midway,above] {$a_{\i}$} (v\j);
    \draw[->]
      (v\i) to[bend right=22]
      node[midway,below] {$b_{\i}$} (v\j);
  }
\end{tikzpicture}
\caption{The quiver $Q$ defining $A_5$.
The relations are $a_i a_{i+1}=0$ and
$b_i b_{i+1}=0$ for $1\leq i\leq3$:
two consecutive upper arrows or two consecutive lower
arrows compose to zero.}
\label{fig:counterexample-quiver}
\end{figure}
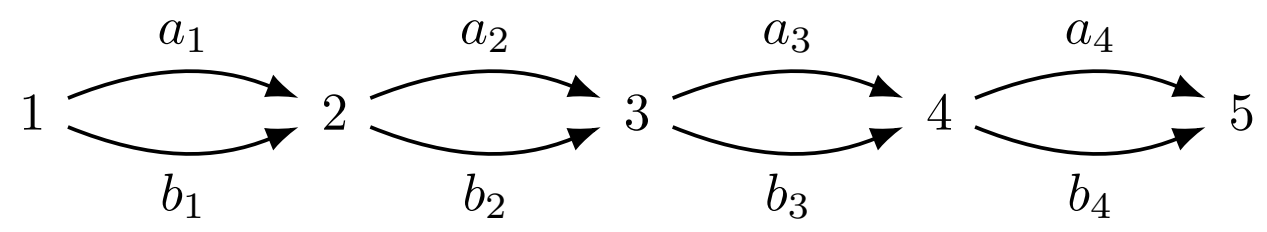

\begin{example}
\label{thm:intro-counterexample}
Let $Q$ be the quiver in
\cref{fig:counterexample-quiver}, and put
\[
  A_5
  :=kQ/(a_i a_{i+1},b_i b_{i+1}\mid1\leq i\leq3).
\]
This is a connected gentle algebra with five simple modules.
Its polygon-dual ribbon graph has vertex words
\[
  W_1=W_2=(1,2,3,4,5).
\]
The associated surface has genus two and one boundary
component, and $\AG(A_5)=[(2,8)]$.
The cycles obtained by traversing one edge and returning
along another have winding zero and span homology.
Thus $\sigma(\eta_{A_5})=0$, and
\cref{thm:noninterval-classification} gives
\begin{equation}
\label{eq:intro-counterexample}
  \GSpec(A_5)=\{4,6,8\}.
\end{equation}
This gives a negative answer to \cref{prob:interval}.

We mention that the same algebra also exhibits nontransitivity of
exceptional-sequence mutation.
Since its surface has genus two, one boundary component,
and two $\circ$-points,
\cite[Theorem~5.4]{ChangHaidenSchroll} shows that the action
of $\ZZ^5\rtimes B_5$ on full exceptional sequences in
$\Db(\modu A_5)$ has infinitely many orbits.
Here $\ZZ^5$ acts by independent shifts of the objects.

Its silting exchange graph is also disconnected.
To see this, let $\alpha$ be the autoequivalence induced
by interchanging $a_i$ and $b_i$, and set
$e=e_1+\cdots+e_4$ and $P_5=e_5A_5$.
The corner algebra $eA_5 e$ is the four-vertex
algebra in \cite[Proposition~4.4]{Dugas}.
Under the natural equivalence
\[
  \per(eA_5 e)\simeq\thick(eA_5),
\]
the construction in its proof gives a tilting object
$T$ of $\thick(eA_5)$ with $\alpha T\not\cong T$.
Since $\Hom(T,P_5[j])=0$ for every $j\in\ZZ$,
the object $T\oplus P_5[-N]$ is silting for $N\gg0$.
It is not $\alpha$-invariant, whereas every indecomposable
projective is $\alpha$-invariant.
By \cite[Corollary~2.2]{Dugas}, it therefore lies in a
different silting component from $A_5$.

More generally, for every $g\geq2$, the analogous algebra
$A_{2g+1}$ on $2g+1$ vertices satisfies
$\GSpec(A_{2g+1})=\{2g,2g+2,\ldots,4g\}$.
The braid group action on its full exceptional sequences
has infinitely many orbits, even modulo independent shifts,
and its silting exchange graph is disconnected, by the
same arguments.
\end{example}

\begin{remark}
The classification also covers the question of gaps if infinite global
dimension is allowed: finiteness of global dimension is derived
invariant, so such a class has the singleton extended spectrum
$\{\infty\}$.  It supplies no further noninterval example.
For finite global dimension, all internal gaps are precisely the
missing odd integers in \eqref{eq:classification-gap-spectrum}.
\end{remark}

\appendix
\section{Global dimension and silting mutation}
\label{sec:HRS-components}

We combine the spectrum classification with the known
global-dimension estimate for $t$-structures to obtain
consequences for mutation of tilting objects. We also
record the interpolation available within a fixed
connected component of the silting exchange graph.

Throughout this section, let $A$ be a nonzero
finite-dimensional algebra of finite global dimension,
and put $\cT:=\per(A)=\Db(\modu A)$.
A basic object $M\in\cT$ is \emph{silting} if
$\Hom_{\cT}(M,M[i])=0$ for every $i>0$ and
$\thick(M)=\cT$; it is \emph{tilting} if the same
vanishing holds for every $i\neq0$.
The silting exchange graph has isomorphism classes of basic silting objects
as vertices and irreducible silting mutations as edges
\cite{AiharaIyama}.

By the silting--$t$-structure correspondence
\cite{KoenigYang}, the subcategory
\[
 \cH_M:=\{X\in\cT\mid
 \Hom_{\cT}(M,X[i])=0\text{ for every }i\neq0\}
\]
is a bounded length heart with finitely many simple
objects. Following the global dimension of a heart
considered in \cite[Remark~5.5]{IkedaQiu}, set
\begin{equation}
\label{eq:delta-definition}
 \delta(M):=
 \sup\{n\in\NN\mid
 \Hom_{\cT}(X,Y[n])\neq0
 \text{ for some }X,Y\in\cH_M\}.
\end{equation}
We call this the \emph{relative global dimension} of
$\cH_M$ in $\cT$ to emphasize the ambient category.
A composition-series argument in both variables shows
that the supremum can be computed on pairs of simple
objects. It is finite, since there are finitely many
such pairs and each simple object belongs to
$\Kb(\proj A)$.

If $M$ is tilting, Rickard's equivalence identifies
$\cH_M$ with the standard heart for
$B_M:=\operatorname{End}_{\cT}(M)$, so
\[
 \delta(M)=\gldim B_M.
\]
Conversely, every value of $\GSpec(A)$ is attained in
this way by a tilting object \cite{Rickard}.

For bounded $t$-structures, the global dimension in
\cite[Definition~3.3]{ChenLinZhou} agrees with the
ambient-heart dimension used above. Their estimate
therefore gives the following consequence.

\begin{proposition}
\label{prop:silting-dimension-interval}
If $M$ and $N$ are adjacent vertices of the silting
exchange graph, then
\[
 |\delta(M)-\delta(N)|\leq1.
\]
Consequently, for every connected component $\mathscr C$,
the set
\[
 \{\delta(M)\mid[M]\in\mathscr C\}
\]
is an interval of nonnegative integers, either finite
or of the form $[a,\infty)$.
\end{proposition}

\begin{proof}
Under the Koenig--Yang correspondence, adjacent
silting objects have hearts related by a simple HRS
tilt. The associated $t$-structures have distance
at most one in the sense of
\cite[Definition~4.2]{ChenLinZhou}, so
\cite[Lemma~4.6(4)]{ChenLinZhou} gives the inequality.
Along a finite mutation path, the integer-valued
function $\delta$ therefore assumes every integer
between its endpoint values. Since any two vertices
in $\mathscr C$ are joined by such a path, its image
is an integer interval. The asserted forms follow
because this image is a nonempty subset of $\NN$.
\end{proof}

Let $\operatorname{TiltEx}(A)$ denote the full subgraph
of the silting exchange graph induced by tilting objects.
The spectrum constrains paths in this subgraph.

\begin{corollary}
\label{cor:intro-tilting-gap}
If $\GSpec(A)$ is not an interval, then
$\operatorname{TiltEx}(A)$ is disconnected.
More precisely, if
\[
 p<d<q,\qquad p,q\in\GSpec(A),\qquad d\notin\GSpec(A),
\]
then every silting mutation path from a tilting object
of relative global dimension $p$ to one of relative
global dimension $q$ passes through a non-tilting
object of relative global dimension $d$.
\end{corollary}

\begin{proof}
Proposition~\ref{prop:silting-dimension-interval}
forces the value $d$ to occur along every such path.
A vertex with this value cannot be tilting, since at
tilting objects $\delta$ belongs to $\GSpec(A)$.
Rickard's theorem supplies tilting objects attaining
$p$ and $q$, so they cannot be connected in
$\operatorname{TiltEx}(A)$.
\end{proof}

The classification gives a concrete consequence for
the exceptional gentle classes.

\begin{corollary}
\label{cor:classified-tilting-components}
Let $A$ satisfy the equivalent conditions of
\cref{thm:noninterval-classification}, with
$\AG(A)=[(2,4g)]$ and $g\geq2$.
Then $\delta$ is constant on every connected component
of $\operatorname{TiltEx}(A)$, and this graph has at
least $g+1$ connected components.
\end{corollary}

\begin{proof}
The classification gives
$\GSpec(A)=\{2g,2g+2,\ldots,4g\}$.
Thus adjacent tilting objects have even $\delta$-values
differing by at most one, and these values must be equal.
Each of the $g+1$ values is attained by a tilting object,
so they occur in distinct connected components.
\end{proof}

In particular, $\operatorname{TiltEx}(A_5)$ has at
least three connected components for the five-simple
algebra in \cref{thm:intro-counterexample}.

\begin{remark}
\label{rem:tilting-cotilting-existence}
We use tilting--cotilting equivalence for the equivalence relation
generated, up to Morita equivalence, by passing to endomorphism algebras
of classical tilting modules of projective dimension at most one or
classical cotilting modules of injective dimension at most one.
For algebras of finite global dimension, each such step changes
global dimension by at most one
\cite[Proposition~III.3.4]{Happel1988}; the cotilting case follows
by duality.

Let $A$ satisfy the conditions of
\cref{thm:noninterval-classification}, with $\AG(A)=[(2,4g)]$ and
$g\geq2$. By that theorem, there are algebras $B_0,\ldots,B_g$
derived equivalent to $A$ such that
\[
  \gldim B_j=2g+2j\qquad(0\leq j\leq g).
\]
Every algebra in a tilting--cotilting chain between two of these
algebras would remain derived equivalent to $A$, and hence would have
even global dimension. The one-step estimate therefore forces global
dimension to remain constant along the chain. Since the displayed
values are distinct, the algebras $B_0,\ldots,B_g$ are pairwise not
tilting--cotilting equivalent. Thus the derived-equivalence class of
$A$ contains at least $g+1$ tilting--cotilting equivalence classes.

In particular, since $\gldim A_{2g+1}=2g$, we may take
$B_0=A_{2g+1}$ and choose $B_1$ of global dimension $2g+2$.
For every odd number $2g+1\geq5$ of simple modules, this recovers
the existence phenomenon established by Ni~\cite[Theorem~1.2]{Ni2026}:
derived-equivalent gentle algebras need not be
tilting--cotilting equivalent.
\end{remark}

\begin{remark}
\label{rem:silting-dg-interpolation}
Let $M_0$ and $M_1$ be tilting objects in the same
silting component, with $\delta(M_0)=p\leq q=\delta(M_1)$.
For every integer $d\in[p,q]$, the proposition gives
a silting object $M_d$ in this component with
$\delta(M_d)=d$.
Its derived endomorphism dg algebra
$\Gamma_{M_d}:=\REnd_A(M_d)$ is proper and connective,
meaning that its total cohomology is finite-dimensional
and $H^{>0}(\Gamma_{M_d})=0$.
Under dg Morita equivalence, its canonical heart
corresponds to $\cH_{M_d}$
\cite{Keller1994,KellerDG,KoenigYang}.
If $d\notin\GSpec(A)$, then $M_d$ is not tilting and
\[
 H^{<0}(\Gamma_{M_d})\neq0.
\]
When $A$ is gentle, $\Gamma_{M_d}$ is quasi-isomorphic
to a nonpositively graded gentle algebra with zero
differential
\cite[Proposition~5.6 and Remark~5.7(2)]{JinSchrollWang}.
For a missing value $d$, every such model has a nonzero
negative-degree part. The common-component hypothesis
above is essential to this interpolation argument.
\end{remark}

\begin{example}
\label{ex:A2-relative-dimension}
Let $A=k(1\longrightarrow2)$, with right projectives
$P_i=e_iA$ and their simple tops $S_i$.
For $r\geq0$, put $M_r=P_1[r]\oplus P_2$.
These objects are silting, and
$\Hom_{\cT}(P_1[r],P_2)=0$ shows that the zero map
$P_1[r]\to0$ is a minimal left
$\operatorname{add}(P_2)$-approximation.
Thus mutation at $P_1[r]$ sends $M_r$ to $M_{r+1}$,
so all $M_r$ lie in the component containing $A$.

For $r\geq1$, the heart is the semisimple category
$\cH_{M_r}=\operatorname{add}\{S_1[r],S_2\}$.
Its intrinsic global dimension is zero, whereas
\[
 \Hom_{\cT}(S_1[r],S_2[r+1])
 \cong\Ext_A^1(S_1,S_2)\cong k.
\]
All other positive-degree Hom groups between these
simple objects vanish, so $\delta(M_r)=r+1$.
This also holds for $r=0$, since $M_0=A$.
Hence relative global dimensions are unbounded on this
silting component, while
$\AG(A)=[(3,1)]$ and
\cref{thm:intro-derived-interval} give
$\GSpec(A)=\{1\}$.
\end{example}
\section*{Acknowledgments}
The authors thank Martin Kalck for answering our questions concerning
\cref{prob:interval,prob:boundedness}.
The first author was supported by the Fundamental Research Funds for the
Central Universities (grant no.~GK202403003), Natural Science Basic Research Program of Shaanxi Province (grant no.~2026JC-YXQN-008) and the National Natural
Science Foundation of China (grant nos.~12271321 and 12671047). 

\section*{AI usage statement}
The authors acknowledge the assistance of ChatGPT with literature
searches, proofreading, computations, and figure preparation.
The computational assistance mainly involved the numerical
calculations in
Subsections~\ref{subsec:classification-higher-genus},
\ref{subsec:classification-winding-bound},
and~\ref{subsec:classification-genus-one}.
All AI-assisted content has been thoroughly reviewed and verified
by the authors, who take full responsibility for the paper,
including the correctness of its mathematical statements,
proofs, and figures.


\begin{thebibliography}{CJSW25}
\footnotesize
\raggedright
\raggedbottom

\bibitem[AI12]{AiharaIyama}
T.~Aihara and O.~Iyama,
\newblock Silting mutation in triangulated categories,
\newblock \emph{J. Lond. Math. Soc. (2)} \textbf{85} (2012), no.~3,
633--668,
\newblock \href{https://doi.org/10.1112/jlms/jdr055}
{doi:10.1112/jlms/jdr055}.

\bibitem[APS23]{AmiotPlamondonSchroll}
C.~Amiot, P.-G.~Plamondon, and S.~Schroll,
\newblock A complete derived invariant for gentle algebras via winding
numbers and Arf invariants,
\newblock \emph{Selecta Math. (N.S.)} \textbf{29} (2023), Paper No.~30,
36 pp.,
\newblock \href{https://doi.org/10.1007/s00029-022-00822-x}
{doi:10.1007/s00029-022-00822-x}.

\bibitem[AH81]{AssemHappel}
I.~Assem and D.~Happel,
\newblock Generalized tilted algebras of type $A_n$,
\newblock \emph{Comm. Algebra} \textbf{9} (1981), 2101--2125;
erratum, \emph{ibid.} \textbf{10} (1982), 1475,
\newblock \href{https://doi.org/10.1080/00927878108822697}
{doi:10.1080/00927878108822697}.


\bibitem[AS87]{AssemSkowronski}
I.~Assem and A.~Skowro{\'n}ski,
\newblock Iterated tilted algebras of type $\widetilde A_n$,
\newblock \emph{Math. Z.} \textbf{195} (1987), 269--290,
\newblock \href{https://doi.org/10.1007/BF01166463}
{doi:10.1007/BF01166463}.

\bibitem[AG08]{AvellaAlaminosGeiss}
D.~Avella-Alaminos and C.~Gei{\ss},
\newblock Combinatorial derived invariants for gentle algebras,
\newblock \emph{J. Pure Appl. Algebra} \textbf{212} (2008), no.~1,
228--243,
\newblock \href{https://doi.org/10.1016/j.jpaa.2007.05.014}
{doi:10.1016/j.jpaa.2007.05.014}.


\bibitem[BC21]{BaurCoelhoSimoes}
K.~Baur and R.~Coelho Sim\~oes,
\newblock A geometric model for the module category of a gentle algebra,
\newblock \emph{Int. Math. Res. Not. IMRN} \textbf{2021} (2021), no.~15,
11357--11392,
\newblock \href{https://doi.org/10.1093/imrn/rnz150}
{doi:10.1093/imrn/rnz150}.

\bibitem[C26]{Chang}
W.~Chang,
\newblock Geometric models for the algebraic hearts in the derived category
of a gentle algebra,
\newblock \emph{Israel J. Math.} \textbf{271} (2026), no.~2, 941--1000,
\newblock
\href{https://doi.org/10.1007/s11856-025-2801-7}
{doi:10.1007/s11856-025-2801-7}.


\bibitem[CHS25]{ChangHaidenSchroll}
W.~Chang, F.~Haiden, and S.~Schroll,
\newblock Braid group actions on branched coverings and full
exceptional sequences,
\newblock \emph{Adv. Math.} \textbf{472} (2025), 110284,
\newblock \href{https://doi.org/10.1016/j.aim.2025.110284}
{doi:10.1016/j.aim.2025.110284}.

\bibitem[CJS26]{ChangJinSchroll}
W.~Chang, H.~Jin, and S.~Schroll,
\newblock Recollements of partially wrapped Fukaya categories and surface cuts,
\newblock \emph{Int. Math. Res. Not. IMRN} \textbf{2026} (2026), no.~13,
Paper No.~rnag133,
\newblock \href{https://doi.org/10.1093/imrn/rnag133}
{doi:10.1093/imrn/rnag133}.


\bibitem[CLZ23]{ChenLinZhou}
X.-W.~Chen, Z.~Lin, and Y.~Zhou,
\newblock The extensions of $t$-structures,
\newblock \emph{Ark. Mat.} \textbf{61} (2023), 323--342,
\newblock
\href{https://doi.org/10.4310/ARKIV.2023.v61.n2.a4}
{doi:10.4310/ARKIV.2023.v61.n2.a4}.


\bibitem[D21]{Dugas}
A.~Dugas,
\newblock Some algebras that are not silting connected,
\newblock \emph{J. Algebra} \textbf{567} (2021), 371--385,
\newblock \href{https://doi.org/10.1016/j.jalgebra.2020.09.028}
{doi:10.1016/j.jalgebra.2020.09.028}.

\bibitem[GR92]{GabrielRoiter}
P.~Gabriel and A.~V.~Roiter,
\newblock \emph{Representations of Finite-Dimensional Algebras},
\newblock with a chapter by B.~Keller,
\newblock Encyclopaedia of Mathematical Sciences, vol.~73, Algebra VIII,
Springer-Verlag, Berlin, 1992.

\bibitem[GT87]{GrossTucker}
J.~L.~Gross and T.~W.~Tucker,
\newblock \emph{Topological Graph Theory},
\newblock Wiley-Interscience Series in Discrete Mathematics and Optimization,
John Wiley \& Sons, New York, 1987.

\bibitem[HKK17]{HaidenKatzarkovKontsevich}
F.~Haiden, L.~Katzarkov, and M.~Kontsevich,
\newblock Flat surfaces and stability structures,
\newblock \emph{Publ. Math. Inst. Hautes \'Etudes Sci.} \textbf{126}
(2017), 247--318,
\newblock \href{https://doi.org/10.1007/s10240-017-0095-y}
{doi:10.1007/s10240-017-0095-y}.


\bibitem[H87]{Happel1987}
D.~Happel,
\newblock On the derived category of a finite-dimensional algebra,
\newblock \emph{Comment. Math. Helv.} \textbf{62} (1987), 339--389,
\newblock \href{https://doi.org/10.1007/BF02564452}
{doi:10.1007/BF02564452}.

\bibitem[H88]{Happel1988}
D.~Happel,
\newblock \emph{Triangulated Categories in the Representation Theory of
Finite-Dimensional Algebras},
\newblock London Mathematical Society Lecture Note Series, vol.~119,
Cambridge University Press, Cambridge, 1988,
\newblock \href{https://doi.org/10.1017/CBO9780511629228}
{doi:10.1017/CBO9780511629228}.


\bibitem[HRS96]{HappelReitenSmalo}
D.~Happel, I.~Reiten, and S.~O.~Smal\o,
\newblock \emph{Tilting in Abelian Categories and Quasitilted Algebras},
\newblock Mem. Amer. Math. Soc. \textbf{120} (1996), no.~575,
\newblock \href{https://doi.org/10.1090/memo/0575}
{doi:10.1090/memo/0575}.

\bibitem[HR82]{HappelRingel}
D.~Happel and C.~M.~Ringel,
\newblock Tilted algebras,
\newblock \emph{Trans. Amer. Math. Soc.} \textbf{274} (1982), no.~2,
399--443,
\newblock \href{https://doi.org/10.1090/S0002-9947-1982-0675063-2}
{doi:10.1090/S0002-9947-1982-0675063-2}.


\bibitem[HKP19]{HKP2019}
A.~Hochenegger, M.~Kalck, and D.~Ploog,
Spherical subcategories in representation theory,
Math. Z. \textbf{291} (2019), no.~1--2, 113--147,
doi:10.1007/s00209-018-2075-4.

\bibitem[IQ23]{IkedaQiu}
A.~Ikeda and Y.~Qiu,
\newblock $q$-Stability conditions on
Calabi--Yau--$\mathbb X$ categories,
\newblock \emph{Compos. Math.} \textbf{159} (2023),
1347--1386,
\newblock
\href{https://doi.org/10.1112/S0010437X23007194}
{doi:10.1112/S0010437X23007194}.

\bibitem[JSW25]{JinSchrollWang}
H.~Jin, S.~Schroll, and Z.~Wang,
\newblock A complete derived invariant and silting theory for graded gentle
algebras,
\newblock arXiv:2303.17474v2, 2025,
\newblock \href{https://doi.org/10.48550/arXiv.2303.17474}
{doi:10.48550/arXiv.2303.17474}.




\bibitem[K15]{Kalck2015}
M.~Kalck,
Singularity categories of gentle algebras,
Bull. Lond. Math. Soc. \textbf{47} (2015), no.~1, 65--74,
doi:10.1112/blms/bdu093.

\bibitem[K17]{KalckDerivedUnique}
M.~Kalck,
\newblock Derived categories of quasi-hereditary algebras and their derived
composition series,
\newblock in \emph{Representation Theory---Current Trends and Perspectives},
EMS Ser. Congr. Rep., European Mathematical Society, Z\"urich, 2017,
pp.~269--308,
\newblock \href{https://doi.org/10.4171/171-1/11}
{doi:10.4171/171-1/11}.


\bibitem[K94]{Keller1994}
B.~Keller,
\newblock Deriving DG categories,
\newblock \emph{Ann. Sci. {\'E}cole Norm. Sup. (4)} \textbf{27} (1994),
63--102,
\newblock \href{https://doi.org/10.24033/asens.1689}
{doi:10.24033/asens.1689}.

\bibitem[K06]{KellerDG}
B.~Keller,
\newblock On differential graded categories,
\newblock in \emph{International Congress of Mathematicians. Vol. II},
European Mathematical Society, Z\"urich, 2006, pp.~151--190,
\newblock \href{https://doi.org/10.4171/022-2/8}
{doi:10.4171/022-2/8}.

\bibitem[KK20]{KellerKrause}
B.~Keller and H.~Krause,
\newblock Tilting preserves finite global dimension,
\newblock \emph{C. R. Math.} \textbf{358} (2020), no.~5, 563--570,
\newblock \href{https://doi.org/10.5802/crmath.72}
{doi:10.5802/crmath.72}.

\bibitem[KY14]{KoenigYang}
S.~Koenig and D.~Yang,
\newblock Silting objects, simple-minded collections, $t$-structures and
co-$t$-structures for finite-dimensional algebras,
\newblock \emph{Doc. Math.} \textbf{19} (2014), 403--438,
\newblock \href{https://doi.org/10.4171/DM/451}
{doi:10.4171/DM/451}.

\bibitem[L08]{Ladkani}
S.~Ladkani,
\newblock Bounds on the global dimension of certain piecewise hereditary
categories,
\newblock \emph{J. Pure Appl. Algebra} \textbf{212} (2008), 2140--2145,
\newblock \href{https://doi.org/10.1016/j.jpaa.2007.12.007}
{doi:10.1016/j.jpaa.2007.12.007}.

\bibitem[LP20]{LekiliPolishchuk}
Y.~Lekili and A.~Polishchuk,
\newblock Derived equivalences of gentle algebras via Fukaya categories,
\newblock \emph{Math. Ann.} \textbf{376} (2020), 187--225,
\newblock \href{https://doi.org/10.1007/s00208-019-01894-5}
{doi:10.1007/s00208-019-01894-5}.

\bibitem[LGH24]{LiuGaoHuang}
Y.-Z.~Liu, H.~Gao, and Z.~Huang,
\newblock Homological dimensions of gentle algebras via geometric models,
\newblock \emph{Sci. China Math.} \textbf{67} (2024), no.~4, 733--766,
\newblock \href{https://doi.org/10.1007/s11425-022-2120-8}
{doi:10.1007/s11425-022-2120-8}.

\bibitem[Ni26]{Ni2026}
Z.~Ni,
\newblock Derived equivalent gentle algebras which are not
tilting-cotilting equivalent,
\newblock arXiv:2608.22929, 2026,
\newblock \href{https://doi.org/10.48550/arXiv.2608.22929}
{doi:10.48550/arXiv.2608.22929}.

\bibitem[OPS18]{OpperPlamondonSchroll}
S.~Opper, P.-G.~Plamondon, and S.~Schroll,
\newblock A geometric model for the derived category of gentle algebras,
\newblock arXiv:1801.09659, 2018 (version~7, 2025),
\newblock \href{https://doi.org/10.48550/arXiv.1801.09659}
{doi:10.48550/arXiv.1801.09659}.

\bibitem[R89]{Rickard}
J.~Rickard,
\newblock Morita theory for derived categories,
\newblock \emph{J. Lond. Math. Soc. (2)} \textbf{39} (1989), 436--456,
\newblock \href{https://doi.org/10.1112/jlms/s2-39.3.436}
{doi:10.1112/jlms/s2-39.3.436}.



\bibitem[S99]{Schroer1999Modules}
J.~Schr\"oer,
Modules without self-extensions over gentle algebras,
J. Algebra \textbf{216} (1999), no.~1, 178--189,
doi:10.1006/jabr.1998.7696.


\bibitem[S23]{SchroerAtlas}
J.~Schr\"oer,
\newblock \emph{Atlas of Finite-Dimensional Algebras},
\newblock version of 29 December 2023, Section~13.5,
\href{https://www.math.uni-bonn.de/people/schroer/FD-Atlas.htmpl?language=en}
{available from the author's webpage}.

\bibitem[SZ03]{SchroerZimmermann}
J.~Schr\"oer and A.~Zimmermann,
\newblock Stable endomorphism algebras of modules over special biserial
algebras,
\newblock \emph{Math. Z.} \textbf{244} (2003), 515--530,
\newblock \href{https://doi.org/10.1007/s00209-003-0492-4}
{doi:10.1007/s00209-003-0492-4}.

\bibitem[S26]{SchrollSurvey}
S.~Schroll,
\newblock On geometric models in representation theory,
\newblock in \emph{Proceedings of the International Congress of
Mathematicians 2026. Vol.~3: Invited Lectures (Sections 1--4)},
SIAM, 2026, pp.~250--271,
\newblock \href{https://doi.org/10.1137/25M1810647}
{doi:10.1137/25M1810647}.

\bibitem[X79]{Xuong}
N.~H.~Xuong,
\newblock How to determine the maximum genus of a graph,
\newblock \emph{J. Combin. Theory Ser. B} \textbf{26} (1979), no.~2,
217--225,
\newblock \href{https://doi.org/10.1016/0095-8956(79)90058-3}
{doi:10.1016/0095-8956(79)90058-3}.

\end{thebibliography}
\end{document}